\documentclass[10pt]{article}

\usepackage{amsmath}
\usepackage{amssymb}
\usepackage[utf8]{inputenc}
\usepackage{tikz,pgfplots}
    \usetikzlibrary{patterns}
\usepackage[space,noadjust]{cite}
\usepackage{tabularx}
\usepackage{multirow}
\usepackage{rotating}
\usepackage{graphicx}
\usepackage{caption,subcaption}
\usepackage[algo2e,vlined,longend,linesnumbered,ruled]{algorithm2e}
\usepackage[software,hardware]{mymacros}
\usepackage{mathdots}

\pgfkeys{/pgf/number format/.cd,1000 sep={\,}}
\pgfplotscreateplotcyclelist{res}{%
{color=blue, thick,densely dashed, mark=*, mark options={solid,mark size=1.5pt}, mark repeat=1},
{color=blue, thick,densely dashed, mark=triangle*, mark options={solid,mark size=1.5pt}, mark repeat=1},
{color=mpigreen, thick,solid, mark=triangle*, mark options={solid,mark size=1.5pt}, mark repeat=1},
{color=red, thick,solid, mark=*, mark options={solid,mark size=1.5pt}, mark repeat=1},
{color=orange, dashdotted,thick, mark=diamond, mark options={fill=orange,solid}, mark repeat=1},
{color=red, densely dashed, mark=pentagon, mark options={solid}, mark repeat=1},
{color=mpigreen, densely dashed,solid, mark=square*, mark options={solid,mark size=1.5pt}, mark repeat=1},
{color=mpigreen, densely dotted,thick, mark=star, mark options={fill=orange,solid}, mark repeat=1},
{color=red, densely dashed, mark=triangle, mark options={solid}, mark repeat=1},
}

\pgfplotscreateplotcyclelist{resN}{%
{color=blue, thick,densely dashed, mark=*, mark options={solid,mark size=1.5pt}, mark repeat=5},
{color=blue, thick,densely dashed, mark=triangle*, mark options={solid,mark size=1.5pt}, mark repeat=5},
{color=mpigreen, thick,solid, mark=triangle*, mark options={solid,mark size=1.5pt}, mark repeat=5},
{color=red, thick,solid, mark=*, mark options={solid,mark size=1.5pt}, mark repeat=5},
{color=orange, dashdotted,thick, mark=diamond, mark options={fill=orange,solid}, mark repeat=5},
{color=red, densely dashed, mark=pentagon, mark options={solid}, mark repeat=5},
{color=mpigreen, densely dashed,solid, mark=square*, mark options={solid,mark size=1.5pt}, mark repeat=5},
{color=mpigreen, densely dotted,thick, mark=star, mark options={fill=orange,solid}, mark repeat=5},
{color=red, densely dashed, mark=triangle, mark options={solid}, mark repeat=5},
}

\newcommand{\change}[1]{\textcolor{black}{#1}}
\newcommand{\MMESS}{\textsf{M-M.E.S.S.}}
\newcommand{\trp}{\ast}

\newcommand{\jit}{{\ensuremath{j_{\text{it}}}}}
\newcommand{\jgal}{{\ensuremath{j_{\text{Gal}}}}}
\newcommand{\jitgek}{{\ensuremath{j^{\text{(G)EKSM}}_{\text{it}}}}}
\newcommand{\jitrk}{{\ensuremath{j^{\text{RKSM}}_{\text{it}}}}}

\newcommand{\jitir}{{\ensuremath{j^{\text{ILRSI}}_{\text{it}}}}}
\newcommand{\jitradi}{{\ensuremath{j^{\text{RADI}}_{\text{it}}}}}
\newcommand{\jitnkgp}{{\ensuremath{j^{\text{NK+GP}}_{\text{it}}}}}
\newcommand{\jitnkls}{{\ensuremath{j^{\text{iNK+LS}}_{\text{it}}}}}

\renewcommand{\Real}[1]{\ensuremath{\operatorname{Re} (#1)}}
\renewcommand{\Imag}[1]{\ensuremath{\operatorname{Im} (#1)}}

\definecolor{mpigreen}{HTML}{5C871D}
\newenvironment{sbm}{\left[\begin{smallmatrix}}{\end{smallmatrix}\right]}

\usepackage[colorlinks=true]{hyperref}
\usepackage{amsthm}

\newtheorem{example}{Example}
\newtheorem{remark}{Remark}

\usepackage{fullpage}

\title{A numerical comparison of solvers for large-scale, continuous-time algebraic Riccati equations \change{and LQR problems}}%
\author{Peter Benner\thanks{Max Planck Institute for Dynamics of Complex Technical Systems,
    Computational Methods in Systems and Control Theory, Magdeburg, Germany,
    {\tt{benner@mpi-magdeburg.mpg.de}, \tt{saak@mpi-magdeburg.mpg.de}}}
\and Zvonimir Bujanovi\'{c}\thanks{University of Zagreb,
    Faculty of Science, Department of Mathematics, Zagreb, Croatia,
    {\tt{zbujanov@math.hr}}}
\and Patrick K\"{u}rschner\thanks{Group Science,
    Engineering and Technology, KU Leuven Kulak, E. Sabbelaan 53, 8500 Kortrijk, Belgium and Department of Electrical Engineering ESAT/STADIUS, KU Leuven, Kasteelpark Arenberg 10, 3001 Leuven, Belgium,
    \tt{patrick.kurschner@kuleuven.be}}
\and Jens Saak\footnotemark[1]}

\hypersetup{%
pdftitle = {A numerical comparison of different solvers for large-scale, continuous-time algebraic Riccati equations and LQR problems}, %
pdfsubject = {Matrix Equations, \today}, %
pdfauthor = {Peter Benner, Zvonimir Bujanovi\'{c}, Patrick K\"urschner, Jens Saak}, %
pdfkeywords = {Algebraic Riccati equation, Lyapunov equation, alternating
direction implicit, rational Krylov subspaces, Newton's method}
}

\begin{document}
\maketitle

\begin{abstract}
  In this paper, we discuss numerical methods for solving large-scale
  continuous-time algebraic Riccati equations. These methods have been the focus
  of intensive research in recent years, and significant progress has been made
  in both the theoretical understanding and efficient implementation of various
  competing algorithms. There are several goals of this manuscript: first, to
  gather in one place an overview of different approaches for solving
  large-scale Riccati equations, and to point to the recent advances in each of
  them. Second, to analyze and compare the main computational ingredients of
  these algorithms, to detect their strong points and their potential
  bottlenecks. And finally, to compare the effective implementations of all
  methods on a set of relevant benchmark examples, giving an indication of their
  relative performance.
\end{abstract}

\section{Introduction}
Let $A,~M\in\Rnn$, $C\in\Rpn$, $B\in\Rnm$ be given matrices. Assuming that \(A,
M\) are sparse, $M$ is nonsingular and \(p,m\ll n\), we consider large-scale, generalized, continuous-time, algebraic Riccati equations (GCAREs)
\begin{align}\label{gcare}
  \cR(X)=A^{\trp}XM+M^{\trp}XA-M^{\trp}XBB^{\trp}XM+C^{\trp}C&=0.
\end{align}
Our goal is the fast and efficient computation of a low rank approximation of
a solution matrix $X\in\Rnn$.

For $M=I_n$,~\eqref{gcare} will be referred to as the standard continuous-time, algebraic Riccati equation (CARE).
If $B=0$,~\eqref{gcare} reduces to a generalized, continuous-time,
algebraic Lyapunov equation (GCALE).
GCAREs appear in various areas related to control theory, for
instance, linear-quadratic optimal regulator \change{(LQR)} problems~\cite{LanR95,Loc01}, $H_2$ and $H_\infty$
controller
design, nonlinear controller design via state-dependent Riccati equations~\cite{Cim08}, and balancing-related model
reduction~\cite{morJonS83,morAnt05,morBenH15}. Solving differential Riccati equations by implicit integration
schemes~\cite{Die92,ChoL90,LanMS15,BenM16} can also lead to GCAREs.

\subsection{Preliminaries, assumptions, motivation and goals}
Because of its nonlinear nature,~\eqref{gcare} can have several solutions.
We exclusively restrict ourselves to the usual situation where a stabilizing
solution is sought, i.e., our goal is to find $X=X^{\trp}\geq0$ such that the
spectrum of the closed loop matrix fulfills $\Lambda(A-BB^{\trp}XM,M)\subset\C_-$.
The stabilizing solution exists and is unique~\cite{LanR95,Bin14} when $(A,M,B)$
is stabilizable (i.e.~rank$[A-z M, \;B] = n$, for each value of $z$ in the
closed right half plane),
and $(A,M,C)$
is detectable ($(A^{\trp},M^{\trp},C^{\trp})$ stabilizable). These conditions are generically
fulfilled, and are assumed to hold in the remainder.
We also assume that the singular values of $X$ decay rapidly towards machine
precision. This enables us to numerically deal with large-scale GCAREs~\eqref{gcare} by computing low-rank approximations $X\approx X^{\text{lr}}=ZDZ^{\trp}$
with $Z\in\R^{n\times r}$, $D=D^{\trp}\in\Rrr$, $r\ll n$.
For GCALEs, the singular value decay has been investigated, e.g.,
in~\cite{Pen99,AntSZ02,Gra04,Sab07,Opm10,BakES14,morBenKS16}.
\change{The low-rank approximability of solutions of GCAREs is less
  well understood, see, e.g., ~\cite{BujB14,KreKM19,Opm15} for
  some results.  In particular, a rapid singular decay of the solution
  is often present when $\Lambda(A,M)\subset\C_-$ and $p,m\ll n$,
  i.e., $C^{\trp}C,~BB^{\trp}$ are of low rank, which we assume as
  well and which is common in large-scale applications.}
The computation of low-rank factors $Z$, $D$ is the common backbone
shared by all the discussed algorithms in this paper. An efficient
execution of these methods will also hinge on our assumption $p,m\ll n$.

The focus of this study is to compare the most prominent algorithms for
large-scale GCAREs.  The comparison is carried out with respect to the key
stages of the algorithms which consume the most numerical effort. In particular,
since all the considered methods are of an iterative nature often requiring different number of steps,
 we will analyze a
 single iteration step of each method.
The main work-intensive stages of the methods, that we are
especially concerned with, are:
\begin{itemize}
\item The numerical solution of the large-scale, real- or complex-valued linear
  systems of equations, usually of the form
  \begin{align*}
    (A+\alpha M)V=N, \quad V,N\in\C^{n\times s},
  \end{align*}
  with $A,M$ from~\eqref{gcare}. We assume that we are able to solve these by
  sparse-direct or iterative solvers at an approximately linear
  complexity. \change{Here, we focus on sparse direct
    solvers. Employing iterative solvers deserves further
    investigations, but this would be well beyond the scope of this
    paper due to the large amount of possible combinations of
    iterative schemes and preconditioners that could be used for this purpose}.
  The number of columns $s$ in the right hand side is of special interest. We employ
  the simplifying assumption that the numerical effort to solve a linear system
  scales linearly with $s$.
\item The construction of rectangular matrices with orthonormal columns spanning
  a basis for certain subspaces is a numerical task arising in methods that use
  a Galerkin projection framework. This is typically achieved by a stable
  implementation of a Gram-Schmidt process, e.g., (modified) Gram-Schmidt with
  iterative refinement. In the numerical costs for the orthogonalization we
  will include only the orthogonalization carried out for a Galerkin projection
  regarding~\eqref{gcare}, and not orthogonalization stages arising in other
  parts of the algorithms, such as the computation of shifts or the residual
  norm.
\item A Galerkin projection naturally leads to GCAREs~\eqref{gcare} defined by
  dense matrices of smaller dimension $\ell\ll n$.  Stable, numerical algorithms
  for small, dense GCAREs can be found, e.g., in~\cite[Chapter 1]{BenBKetal15a} and the
  references therein. If not stated otherwise, the \matlab{} routine
  \texttt{care} is used, whose costs are estimated as $\cO(\ell^3)$.
\item The majority of the considered methods relies on certain shift parameters,
  which are important for rapid convergence. The problem of generating and
  selecting high quality shift parameters might easily fill a survey article
  itself and is, therefore, beyond the scope of this study. We refer to the
  relevant literature~\cite{DruS11,SimSM13a,BenKS14,DruSZ14,Sim16a,BenBKetal15,Kue16} and provide only basic, necessary information where
  appropriate. For each particular GCARE solution method, we will employ the shift strategy
  providing the best results. With some exceptions, these are usually the more
  recent adaptive generation strategies.
\item All the considered algorithms are of iterative nature and, hence, require
  appropriate stopping criteria. For consistency, we terminate all
  algorithms when
  \begin{align}
    \|\cR(\tX)\|_2\leq \tau_{\cR}\|C^{\trp}C\|_2,\quad 0<\tau_{\cR}\ll 1,
  \end{align}
  for an approximate solution $\tX\approx X$, and usually choose
  \(\tau_{\cR}=10^{-8}\).
  We employ this residual norm based criterion because approaches based on
  relative changes, e.g., of the norm of the generated approximate feedback
  matrices $\tK=M^{\trp}\tX B$, have been proven unreliable~\cite{Wei16}. For
  large-scale matrix equations, computing or estimating the residual norm is not
  always an easy task and can yield significant numerical
  effort. \change{Recent developments reduce these costs for several
    of the considered methods to the level of
  being insignificant, though.} Details on how
  this is achieved will be included with the descriptions of the single
  algorithms.
\end{itemize}
Only the first and last points arise in all discussed methods.

The memory consumption of the methods is also considered and compared. Here, it is important to distinguish between two scenarios:
\begin{enumerate}
\item[\textbf{S1.}] An approximate solution of~\eqref{gcare} is sought, in the form of a
  low-rank approximate solution $ZDZ^{\trp}\approx X$ with solution factors
  $Z\in\R^{n\times r},~D\in\Rrr$ and $r\ll n$. This situation occurs especially
  in certain model order reduction techniques, e.g., LQG-, positive-real, and
  related balanced truncation type
  approaches~\cite{morJonS83,morDesP84,morOpdJ88,morAnt05,morBenH15}.  Note that
  in some of these model reduction approaches, the arising GCAREs can have
  positive definite or indefinite quadratic terms, e.g.~\cite{BenS14}. We do not pursue this issue
  any further here, but these can be treated with a low-rank version of the
  iteration proposed in~\cite{LanFA07}.
\item[\textbf{S2.}] Only an approximation of the stabilizing feedback matrix
  $K:=M^{\trp}XB\in\Rnm$ is sought. The optimal control of linear, time-invariant
  dynamical systems using LQR-, or LQG-feedback control approaches (e.g.~\cite{Saa09,Hei09}) is the
  prototypical application for this scenario. Since $K$ is a very thin
  $n\times m$ matrix, algorithms that are able to operate solely on
  approximations of $K$
  are considerably more efficient regarding the memory requirements, whenever
  \(m\) is significantly smaller than \(r\).
\end{enumerate}
\begin{remark}
 If one is only interested in the stabilizing feedback $K$ as in scenario \textbf{S2}, alternative approaches (for $M=I$) based on the Chandrasekar
differential equations
 \begin{subequations}\label{Chandr_ODE}
  \begin{align}
    \dot K(t) &= -{L(t)}^{\trp}(L(t)B),&K(0)=0, \\
    \dot L(t) &= -L(t)(A-B{K(t)}^{\trp}),&L(0)=C
 \end{align}
 \end{subequations}
can be used~\cite{BanI91,CosF97} by choosing a sufficiently large $t=t_f$ and solving the system of ordinary differential equations numerically backwards in
time. The motivation behind this is that $K=\lim\limits_{t\rightarrow -\infty} K(t)$ and $X=\lim\limits_{t\rightarrow -\infty}\intab{0}{t}
{L(s)}^{\trp}L(s)\mathrm{d}s$~\cite{Kai72,Cas77}.
Unfortunately, a numerical solution of~\eqref{Chandr_ODE} can be difficult, because the solution trajectories often exhibit a highly oscillatory behavior for
small times. An integrator may thus require a large number of tiny time steps, and it can take a long time for the above dynamical system to reach a
stationary phase~\cite{Rus79}. Therefore we do not pursue this approach in this study but leave it as an interesting future topic, especially since it negates the
need to work with an ARE at all in scenario \textbf{S2}. Similarly to the strategy in~\cite{BanI91}, one could solve~\eqref{Chandr_ODE} numerically until a
stabilizing $K$ is constructed, which then serves as an initial guess within an iterative method for the GCARE~\eqref{gcare}.
\end{remark}
\subsection{Outline}
In Section~\ref{sec:methods}, the considered methods are briefly described and categorized into three classes. We refrain from giving detailed derivations and
theoretical results for each single algorithm, as those can be found in the relevant literature. How the above computational stages arise in each method is
emphasized in that section.
Recent developments, which are important for the numerical
performance, are also mentioned\change{, including a few modifications
  so far not discussed in the literature such as, e.g., accessing the true GCARE residual
    norm in projection methods, defect correction strategies for occurring projected CAREs, and handling initial feedback approximations}. The comparison of the computational costs  is
given in Section~\ref{sec:compare}. This addresses a single iteration step of
each considered method respecting the subtasks mentioned above. The memory
requirements and consumption are addressed as well, especially regarding the scenarios
\textbf{S1} and \textbf{S2}.
A series of numerical experiments is carried out in Section~\ref{sec:numex}
comparing the computation time, memory consumption, and other performance
indicators of all methods. Section~\ref{sec:concl} summarizes our findings.
\subsection{Notation}
We use the following notation in this paper: $\C_-,\C_+$ are the open left and right half plane, $\Real{z},~\Imag{z}$, $\conj{z}=\Real{z}-\jmath\Imag{z}$,
$|z|$
are the real part, the imaginary part, the complex conjugate, and the absolute
value of a complex quantity $z$, respectively. For matrices, $A^{\trp}$ denotes
the transpose for real matrices, and the complex conjugate transpose for complex matrices.
If it exists, $A^{-1}$ is the inverse of $A$ and $A^{-\trp}={(A^\trp)}^{-1}$.
In most situations, expressions of the form $x=A^{-1}b$ are to be understood as
solving
the linear system of equations $Ax=b$ for $x$.
The relations $A\succ(\succeq)0$ and $A\prec(\preceq)0$ stand for the matrix $A$ being positive and negative (semi)-definite. Likewise, $A\succeq(\preceq) B$ refers to
$A-B\succeq(\preceq)0$. Unless stated otherwise, $\|\cdot\|$ is the Euclidean vector or subordinate matrix norm, and $\kappa(\cdot)$ is the associated condition
number. The Frobenius norm is denoted by $\|\cdot\|_F$.
For a given matrix $A \in \Rnn$, and a given (block) vector $v \in \Rnp$, the (block) Krylov subspace generated by $A$ and $v$ is denoted as
\[
    \cK_j(A, v) := \range{\left[v, \; Av, \; A^2v, \; \ldots, \; A^{j-1}v\right]}.
\]

\section{Classification and brief introduction of the considered
  methods}\label{sec:methods}
This section briefly introduces the methods compared in the later sections. We
divide the methods into three classes represented by the subsection herein. The
first class consists of all methods that use a certain subspace, to project~\eqref{gcare} to a much smaller dense representation, that can be solved by a
direct computation. Then the solution is lifted to the full coordinates
again. Direct iteration methods that successively approximate the solution,
without the need for projection, form the second class. The remaining class of
methods is made up out of the variants of Newtons method.
\subsection{Projection methods}\label{ssec:projectionmethods}
Let $\cQ$ be an $r$-dimensional subspace of $\Rn$ with $r\ll n$ and let the columns of $Q_r\in\Rnr$ form an orthonormal basis for $\cQ$.
We are looking for approximate solutions of~\eqref{gcare} in the space
\begin{align*}
 \cZ_r(Q_r):=\lbrace X_r=Q_{r}Y_{r}Q_{r}^{\trp}~:~Y_{r}^{\trp}\in\Rrr\rbrace.
\end{align*}
For the standard case $M=I_n$, in direct analogy to the case of Lyapunov equations~\cite{Saa90,morJaiK94}, we impose a Galerkin condition (using the Euclidean inner product) onto the CARE
residual:
\begin{align*}
 \cR(X_r)=A^{\trp}Q_{r}Y_{r}Q_{r}^{\trp}+Q_{r}Y_{r}Q_{r}^{\trp}A-Q_{r}Y_{r}Q_{r}^{\trp}BB^{\trp}Q_{r}Y_{r}Q_{r}^{\trp}+C^{\trp}C\quad\perp\quad\cZ_r(Q_r).
\end{align*}
This condition implies that $Y_r$ is the solution of the $r$-dimensional CARE
\begin{align}\label{proj_care}
 \tA_{r}Y_{r}+Y_{r}\tA^{\trp}_{r}+Y_{r}\tB_{r}\tB_{r}^{\trp}Y_{r}+\tC_{r}^{\trp}\tC_{r}=0
\end{align}
with $\tA_r:=Q_r^{\trp}A^{\trp}Q_r\in\Rrr$,~$\tB_r:=Q_r^{\trp}B\in\Rrm$,~$\tC_r:=CQ_r\in\Rpr$.
For $M\neq I_n$, this Galerkin projection is typically implicitly applied to an equivalent CARE defined, e.g., by $A_M:=AM^{-1}$, $B_M:=B$ and $C_M=CM^{-1}$ or,
if $0\prec M=L_{M}L_{M}^{\trp}$, to  $A_{M}:=L_{M}^{-1}AL_{M}^{-\trp}$, $B_{M}:=L_{M}^{-1}B$ and $C_{M}=CL_{M}^{-\trp}$.
If the resulting low-rank solution $X_{r}$ is not good enough, the subspace $\cQ_{r}$ is expanded by additional basis vectors.

Methods following this Galerkin projection principle mainly differ in the way the subspace $\cQ$ or, more precisely, the sequence
$\cQ_0\subseteq\cQ_1\subseteq\ldots\subseteq\cQ_j$ of subspaces, is constructed.
An intuitive choice would be the block Krylov subspace generated from $A_M^{\trp}$ and $C_M^{\trp}$~\cite{morJaiK94} --- however, because of the often
resulting slow convergence,
this approach has been superseded by the application of more general Krylov type subspaces.

The extended block Krylov subspace~\cite{DruK98,KniS10} is given by
\begin{align*}
 \cQ_j&=\cK^{\text{ext}}_j(A_M^{\trp},C_M^{\trp}):=\cK_j(A_M^{\trp},C_M^{\trp})\cup \cK_j(A_M^{-1},A_M^{-1}C_M^{\trp}).
\end{align*}
In each iteration step, the subspace is expanded by $2p$ new vectors, leading to\\dim$(\cQ_j)\leq2jp$.
The basis matrix $Q_j$ for $\cK^{\text{ext}}_j$ can be constructed by the extended block Arnoldi process which was used in~\cite{Sim07,KniS11} and~\cite{HeyJ09} to compute low-rank solutions of GCALEs and GCAREs, respectively. Where it does not lead to confusion, we will omit the prefix
``block'' in the remainder of the text, and refer to algorithms using
$\cK^{\text{ext}}_j$ as extended Krylov subspace methods (EKSM).

A further generalization is given by rational Krylov subspaces~\cite{Ruh84}
\begin{align*}
 \cQ_j&=\cK^{\text{rat}}_j(A_M^{\trp},C_M^{\trp},\boldsymbol{\alpha})\\
&:=\range{\left[{(A_M-\alpha_1I)}^{-\trp}C_{M}^{\trp},\ldots,\prod\limits_{i=1}
^j\left({(A_M-\alpha_{i}I)}^ { -\trp }\right) C_M^{ \trp }\right]}\\
&=\range{\left[{(A-\alpha_1M)}^{-\trp}M^{\trp}C_M^{\trp},\ldots,\prod\limits_{i=1}
^j\left({(A-\alpha_{i}M)}^ { -\trp }M^{\trp}\right) C_M^ { \trp } \right]} ,
\end{align*}
where $\boldsymbol{\alpha}:=\lbrace\alpha_1,\ldots,\alpha_j\rbrace\subset\C_+\cup\lbrace\infty\rbrace$
are shift
parameters whose selection is discussed later.
Following the well known moment matching interpretation~\cite{morAnt05}, we note
that $\cK^{\text{ext}}_j$ is a special case of $\cK^{\text{rat}}_{2j}$ with the shifts zero and infinity used in an alternating fashion.
The usage of rational Krylov subspace methods (RKSM) for GCALEs and GCAREs was
investigated in~\cite{DruS11,SimSM13a,DruKS11,Sim16a}.
The basic procedure of both EKSM and RKSM for~\eqref{gcare} is summarized in Algorithm~\ref{alg:rksm}.
Usually, the
rational Krylov subspace is generated so that $\range{C_M^{\trp}}\subseteq \cQ_j$~\cite{SimSM13a,DruKS11,DruS11}, which can be enforced by formally
setting $\alpha_1=\infty$.
  \LinesNumberedHidden{}
\begin{algorithm2e}[h]
  \SetEndCharOfAlgoLine{} \SetKwInOut{Input}{Input}\SetKwInOut{Output}{Output}
  \caption[Extended and Rational Krylov subspace method for GCAREs]{Extended and Rational Krylov Subspace
    Method for GCAREs}\label{alg:rksm}
  \Input{Matrices $A,~M,~B,~C$ defining~\eqref{gcare} and stopping tolerance $0<\tau_{\cR}\ll 1$.}
  \Output{$Q_j\in\C^{n\times k}$, $Y_j=Y_j^{\trp}\in\C^{k\times k}$ such that
    $Q_{j}Y_{j}Q_{j}^{\trp}\approx X$ with $k\ll n$ and $Q_{j}^{\trp} Q_{j}=I_k$, stabilizing feedback matrix
$K_j\in\Cnm$.}
\nlset{1} $q_0=C_M:=M^{-\trp}C^{\trp}$, $Q_0=\mathrm{orth}(C_M)$.\nllabel{rksmO1}\;
\nlset{2} EKSM:\@ $q_0^{(1)}=q_0$, solve $A^{\trp}q_0^{(2)}=M^{\trp}q_0^{(1)}$ for $q_0^{(2)}$, $Q_0=\mathrm{orth}([Q_0,q_0^{(2)}])$.\nllabel{rksmO2}\;
\nlset{3} \For{$j=1,2,\ldots,j_{\max}$}{
\nlset{4} RKSM:\@ Select shift  $\alpha_j \in \C_+$.\nllabel{rksmS}\;
\nlset{5}  Generate new basis vectors:

\nlset{5a}
   \hspace*{0.5cm} RKSM:~Solve ${(A-\alpha_{j}M)}^{\trp} q_j=M^{\trp}q_{j-1}$ for $q_j$.\nllabel{rksmL1}\;

\nlset{5b}
   \hspace*{0.5cm}
      EKSM:~Solve $M^{\trp}q_j^{(1)}=A^{\trp}q_{j-1}^{(1)}$, $A^{\trp}q_{j}^{(2)}=M^{\trp}q_{j-1}^{(2)}$ for $q_{j}^{(1)}$, $q_{j}^{(2)}$;
          set $q_j=[q_{j}^{(1)}$, $q_{j}^{(2)}]$.\nllabel{rksmL2}\;
\nlset{6} Orthogonally extend basis
    matrix $Q_{j-1}$:
$Q_{j}=\mathrm{orth}([Q_{j-1},~q_j])$.\nllabel{rksmO3}\;
\nlset{7} $\tA_j=Q_j^{\trp} M^{-\trp}A^{\trp}Q_j$, $\tB_j=Q_j^{\trp} B$, $\tC_j=C_{M}Q_j$.\nllabel{rksmAj}\;
\nlset{8} Solve projected CARE
$\tA_{j}Y_{j}+Y_{j}\tA_{j}^{\trp}-Y_{j}\tB_{j}\tB_{j}^{\trp} Y_j+\tC_j^{\trp}\tC_j=0$ for $Y_j$.\nllabel{rksmsmallARE}\;
\nlset{9}\lIf{$\|\cR(Q_{j}Y_{j}Q_{j}^{\trp})\|<\tau_{\cR}\|C^{\trp}C\|$}{$K_j=M^{\trp}Q_{j}Y_{j}\tB_{j}$, \textbf{stop}.}
}
\end{algorithm2e}
\LinesNumbered{}

We now give some remarks on the major steps of
Algorithm~\ref{alg:rksm}. In the Lines~\ref{rksmO1},~\ref{rksmO2},~\ref{rksmO3}, \texttt{orth} is to be understood as any stable (block) orthogonalization
routine, such as a repeated (block) modified Gram-Schmidt process~\cite{GolV13}, which we employ in this study. An efficient construction of the projected
matrices $\tA_j$, $\tC_j$ can be found in
the
respective literature on EKSM~\cite{Sim07} and RKSM~\cite{DruS11,SimSM13a,Gue13,BerG15}.  The small CARE in Line~\ref{rksmsmallARE} can
be solved by direct methods involving dense numerical linear
algebra, see e.g.~\cite{Bin14,BenBKetal15a,Ben99}, whose numerical complexity is cubic in the subspace dimension.

For RKSM, the choice of shift parameters in Line~\ref{rksmS} is crucial to achieve a fast
convergence. An overview of
different selection strategies can be found, e.g., in~\cite{Gue13}. The adaptive selection strategy proposed in~\cite{DruS11} and later improved
in~\cite{Sim16a}
turned out to be successful in the majority of cases. There, after
iteration step $j$ of RKSM, the next shift $\alpha_{j+1}$ is obtained by
minimizing a rational function over the convex hull of the eigenvalues of either
$\tA_j$~\cite{DruS11}, the projected closed loop matrix $\tA_j-\tB_j\tK_j^{\trp}$ with $\tK_j:=Y_j\tB_j$~\cite{Sim16a,LinSim15}, or of the matrix pair
$(Q_j^{\trp}AQ_j,~Q_j^{\trp}MQ_j)$~\cite{ElmR16}. Here, we restrict to the variant using $\tA_j-\tB_j\tK_j^{\trp}$.
It can happen that some of the generated shifts occur in complex conjugate pairs.
In order to reduce the amount of complex arithmetic operations, the basis matrix
$Q_j$, and therefore also most other quantities, can be kept real by applying the
real RKSM proposed in~\cite{Ruh94c}. Essentially, the real RKSM consists of augmenting $Q_{j-1}$ by $[\Real{q_j},~\Imag{q_j}]$ when $\alpha_j\in\C\setminus\R$,
which is equivalent to processing both shifts $\alpha_j$ and $\conj{\alpha_j}$ at once. With this observation, the only remaining complex operation that occurs
in the algorithm is solving the complex, sparse linear system in Line~\ref{rksmL1}.

To monitor the progress and to stop the iteration, the Euclidean or Frobenius norm of the residual $\cR_j:=\cR(Q_{j}Y_{j}Q_{j}^{\trp})$ can
also be computed efficiently without explicitly forming this large, dense $n\times n$ matrix~\cite{Sim07,HeyJ09,DruS11,SimSM13a,Sim16a}.
Following~\cite{DruS11}, the CARE residual in RKSM has the form (for $M=I$)
\begin{align*}
 \cR_j&=F_j+F_j^\trp,\quad F_j:=G_{j}S_{j}^\trp\in\C^{n\times n},\\
 G_j&:=q_{j+1}\alpha_j-(I-Q_{j}Q_{j}^{\trp})A^{\trp}q_{j+1},\quad
S_j:=Q_j^\trp Y_j H_{j}^{-\trp}E_{j}h_{j+1,j} ,\quad E_j=e_j\otimes I_p,
\end{align*}
where $H_j=[h_{i,k}]\in\C^{jp\times jp}$ is block upper Hessenberg, and  $h_{i,k}\in\C^{p\times p}$,
$i=1,\ldots,j+1$, $k=1,\ldots,j$ are the orthogonalization coefficients obtained in the Gram-Schmidt process. Since
$Q_j^\trp q_{j+1} = 0$, it also holds
$S_j^\trp G_j=0$ and, hence, we have the idempotence $F_j^2 = 0$. An easy calculation shows that
$\|\cR_j\| = \|F_j\| = \|\Psi h_{j+1,j} E_j^\trp H_{j}^{-1}Y_j\|$,
where $\Psi$ is the triangular factor of a thin QR factorization of $G_j$. The structure of $G_j$ and $S_j$ is slightly different in the real RKSM~\cite{Ruh94c}
and
in EKSM~\cite{Sim07,HeyJ09}, but the overall approach remains unchanged.
Note that for GCAREs ($M\neq I$) the projection method is implicitly applied to an equivalent CARE defined by matrices $A_M$, $B_M$
and
$C_M$, such that one would obtain the norm of an associated transformed  residual, e.g., $\cR_j^M:=M^{-\trp}\cR_{j}M^{-1}$.
In practice, $\|\cR_j\|$ and $\|\cR^M_j\|$ can differ
significantly. For our comparative study of different algorithms, we prefer to monitor the residual of the original
GCARE. This can be easily achieved by exploiting
\begin{align*}
 \cR_j&=M^\trp \cR_j^M M=M^\trp(F^M_j+{(F^M_j)}^\trp)M=N^M_j\begin{sbm}0_p&I_p\\I_p&0_p\end{sbm} {(N^M_j)}^\trp,\quad F^M_j:=G^M_{j}S_{j}^\trp,
\\
N_j&:=M^\trp\left[G^M_j,S_j\right]\in\C^{n\times 2p},~G^M_j:=q_{j+1}\alpha_j-(I-Q_{j}Q_{j}^{\trp})A_M^{\trp}q_{j+1},
\end{align*}
from which it follows $\|\cR_j\|=\|\Psi_M\begin{sbm}0_p&I_p\\I_p&0_p\end{sbm} \Psi_M^\trp\|$,
where $\Psi_M$ is the triangular factor of a thin QR decomposition of $N_j$.
In addition to the QR decomposition, $2p$ ($4p$ for complex shifts in real RKSM) matrix-vector multiplications with $M^\trp$ (or factors thereof) are
therefore required to compute $\|\cR_j\|$.
\subsubsection{Variants}
To slow down the growth of the column dimension of $Q_j$ in block RKSM for Lyapunov equations, a modification has been presented in~\cite{DruSZ14}. There,
tangential
rational Krylov subspaces
\begin{align*}
 \cQ_j
 	&=\cK^{\text{t-rat}}_j(A_M^{\trp},C_M^{\trp},\boldsymbol{\alpha}) \\
 	&:=\range{\left[{(A-\alpha_1M)}^{-\trp}M^{\trp}
		C_M^{\trp}d_1,\ldots,{(A-\alpha_{j}M)}^{-\trp}M^{\trp}C_M^{\trp}d_j\right]},~\alpha_i\neq\alpha_j,
\end{align*}
are used to generate the projection subspace, leading to the tangential RKSM (TRKSM). It is straightforward to apply this approach to GCAREs.
The tangential directions $d_j\in\C^{p\times p_j}$, $p_j\leq p$ are computed adaptively at every iteration step in conjunction with the shifts. Clearly,
dim$(\cK^{\text{t-rat}}_j)=\sum_{i=1}^{j}p_i\leq jp$, leading to a reduced effort for the orthogonalization and small-scale solution of the projected matrix
equation.
However, as pointed out in~\cite{DruSZ14} the decreased growth of dim$(Q_j)$ can be accompanied with a slower
convergence of the algorithm. Compared to the standard block RKSM, this may result in a larger number of linear systems with different coefficient matrices that
need to be solved. Some of the upcoming numerical experiments confirm this.

In~\cite{Jbi99,Jbi06a}, global standard and extended Krylov subspace methods for matrix equations are considered.
\change{These methods are based on a global orthogonalization process,
where the constructed basis matrix $Q_j=[q_1,\ldots,q_j]$ has blocks $q_i\in\Cnp$ satisfying $\langle q_i,q_k\rangle=\trace{q_k^{\trp}q_i}=\delta_{i,k}$, $\|q_i\|_F=1$. See also~\cite{FroLS17} for more details on these global processes.}
When applied for GCALEs, the projected quantities $\tA_j$, $\tC_j$ obtained in this way have special structures which can be exploited to gain savings when
solving the projected
small-scale problem~\cite{Jbi99,StyS12}. In global Krylov methods for GCAREs, however, such savings are not possible because $\tB_j$ does not reveal any
special structure one could exploit.
Moreover, global methods often converge slower than their
block counterparts, meaning that they require a higher number of iteration steps, which results in a larger subspace dimension.
In the GCARE case this leads to noticeably higher numerical costs
for the solution of the projected equation.
A detailed assessment and comparison of the costs of both block and global extended Krylov subspace methods for GCALEs is
given in~\cite{StyS12}, where further information on actual implementations can be found as well.
For the purpose of comparison in this paper, the global EKSM (GEKSM)~\cite{Jbi06a} is used in some of our numerical experiments.
Global rational Krylov subspace approaches are considered in~\cite{morBonFAetal15} in the context of model order reduction.
\subsection{Non-projective iterations}\label{sec:nonprojectiveiterations}
The second class of methods consists of iterative processes working directly on the Riccati equation, and circumventing the use of a Galerkin projection
framework or
Newton-scheme as the methods in the previous and next section, respectively.
The low-rank solution is built from direct relations that do not require the
solution of a projected CARE.\@ In particular, two representatives of this class
are considered, both of which implicitly work with the Hamiltonian matrix
\begin{align*}
 \cH=\begin{bmatrix}
      A     & -BB^{\trp} \\
      -C^{\trp}C & -A^{\trp}
     \end{bmatrix}\in\R^{2n\times 2n}
\end{align*}
associated to~\eqref{gcare} (here we restrain ourselves to the case $M=I$ for simplicity\change{, but the case $M\neq I$ can be developed in a straightforward manner by considering the associated Hamiltonian/skew-Hamiltonian matrix pencil}). It is well known that the solution of~\eqref{gcare} is given by
$X=QP^{-1}$,
where
${[P^{\trp},Q^{\trp}]}^{\trp}\in\R^{2n\times n}$ spans an $\cH$-invariant subspace with respect to the stable eigenvalues of $\cH$.
Not only is this relation the basis for many direct methods for
CAREs~\cite{Lau79,ChuLM06}, but also several works~\cite{AmoB10,LinSim15,BujB14,MasOR16,BenBKetal15}
investigate its application in a large-scale setting.

\subsubsection{The incremental low-rank subspace iteration}
In~\cite{LinSim15}, the Cayley iteration
\begin{align}
\label{eq::CayleyHamiltonIteration}
\begin{bmatrix} P_j\\Q_j\end{bmatrix}={(\cH+\alpha_j I_{2n})}^{-1}(\cH-\overline{\alpha}_j I_{2n})\begin{bmatrix}
I_n\\X_{j-1}\end{bmatrix}, \quad X^{\text{cay}}_j=Q_{j}P_{j}^{-1},
\end{align}
with shift parameters $\alpha_j\in\C_+$, is proposed for the Hamiltonian matrix $\cH$ to iteratively compute the desired $\cH$-invariant subspace.
For feasibility in the large-scale case, one can observe that if the initial approximation $X_1^{\text{cay}} \succeq 0$ has rank $p$, then $X_j^{\text{cay}}$
has rank $jp$, and
it is possible to rewrite the above iteration into a low-rank version~\cite{LinSim15}, incrementally generating a low-rank approximation $X\approx
X_j^{\text{cay}}=Z_{j}Y_{j}^{-1}Z_j^{\trp}$ with $Z_j\in\C^{n\times jp}$, $0\prec Y_j=Y_j^{\trp}\in\C^{jp\times jp}$. This method is referred to as Incremental
Low-Rank
Subspace Iteration (ILRSI) and illustrated in Algorithm~\ref{alg:ilrsi}.

\begin{algorithm2e}[t]
\SetEndCharOfAlgoLine{}
\SetKwInOut{Input}{Input}\SetKwInOut{Output}{Output}
\caption{Incremental Low-Rank Subspace Iteration (ILRSI)}\label{alg:ilrsi}
\Input{Matrices $A,~M,~B,~C$ defining~\eqref{gcare}}
\Output{$Z_j\in\C^{n\times jp}$, $Y_j\in\C^{jp\times jp}$ s.t. $Z_j
Y_j^{-1} Z_j^{\trp}\approx X$.}
    Determine shifts $\lbrace
\alpha_1,\dots,\alpha_{j_{\max}}\rbrace\subset\C_+$.\;
Solve  $(\alpha_j M^{\trp}-A^{\trp})V_1= C^{\trp}$ for $V_1$. \nllabel{ilrsiL1}\;
Update solution factors $Y_1=2\Real{\alpha_1}(I+V_1^{\trp} BB^{\trp} V_1)$, $Z_1=V_1$.\;
    \For{$j=2,3,\ldots,j_{\max}$}{
    Solve  $(\alpha_j M^{\trp}-A^{\trp})\tV=M^{\trp}V_{j-1}$ for $\tV$,
$V_j=V_{j-1}-(\alpha_j+\overline{\alpha_{j-1}})\tV$.\;
    Update low-rank solution factor $Z_j = [Z_{j-1} \;\; V_j]$.\;
$L_j=\begin{sbm} &I_{jp-1}\\1&\end{sbm}
\begin{sbm}1&\cdots&1\\&\ddots&\vdots\\&&1\end{sbm}
\begin{sbm}\tfrac{\overline{\alpha_1}+\alpha_j}{2\Real{\alpha_j}}&&\\
\tfrac{\alpha_1-\alpha_j}{2\Real{\alpha_j}}&\tfrac{\overline{\alpha_2}+\alpha_j}
{2\Real{\alpha_j}}&\\
 &\ddots&\ddots\\
 &\tfrac{\alpha_{j-1}-\alpha_j}{2\Real{\alpha_j}}&\tfrac{\overline{\alpha_j}+\alpha_j}{2\Real{\alpha_j}}
 \end{sbm}
 \otimes I_p$.\;
 $N_j=L_j^{-1}+\diag{I,0}$, $V_B=L_j^{-\trp}V_j^{\trp} B$.\;
Update middle low-rank factor $Y_j=N_j^{-\trp}\left(\begin{sbm}
Y_{j-1}&\\&2\Real{\alpha_j}\end{sbm}+\tfrac{1}{2\Real{\alpha_j}}V_{B}V_{B}^{\trp}\right)N_j^{-1}$.\;
\lIf{$\|\cR(Z_{j}Y_{j}^{-1}Z_j^{\trp})\|<\tau_{\cR}\|C^{\trp}C\|$}{\textbf{stop}. \nllabel{ilrsiL2}
}
 $j=j+1$\;
    }
\end{algorithm2e}

Interestingly, the same iteration has been derived independently in~\cite{MasOR16}, with a completely different approach by analyzing
an optimal control problem associated to~\eqref{gcare}. Only the setup of the middle factor $Y_j$ is done slightly differently in~\cite{MasOR16}.  If
$B=0$, ILRSI reduces to the low-rank alternating direction implicit (LR-ADI) iteration~\cite{Pen99,LiW02,BenLP08} for GCALEs. The shift parameters $\alpha_j$
are, once again, crucial for the speed of the convergence and may be selected in advance
or,
preferably, adaptively by strategies borrowed from the GCALE case, e.g.,~\cite{Wac13, BenKS14}. Approaches that take the GCAREs structure into account more will
be discussed later.
The relation to the LR-ADI iteration can be exploited to deal with pairs of complex conjugate shifts since~\cite[Theorem 1]{BenKS13} works for
Algorithm~\ref{alg:ilrsi} as well: For $\alpha_j\in\C_+,~\alpha_{j+1}=\overline{\alpha_j}$, it holds
that $V_{j+1}=\overline{V_j}-\beta_j\Imag{V_j}$ with $\beta_j:=2\tfrac{\Real{\alpha_j}}{\Imag{\alpha_j}}$ and
$Z_{j+1} = [Z_{j-1}, \Real{V_j}-\beta_j\Imag{V_j}, \sqrt{(\beta^2+1)}\Imag{V_j}]$. Hence, for each complex pair, a double iteration step can be performed by
solving
only one complex linear system.

The most costly numerical work in ILRSI appears to be the solution of the shifted linear systems with $p$ right hand sides in Line~\ref{ilrsiL1}.
However, the estimation of the residual norm $\|\cR(Z_{j}Y_{j}^{-1}Z_j^{\trp})\|$ in Line~\ref{ilrsiL2} is potentially expensive as well.
The reason is that,
unlike the LR-ADI iteration for GCALEs~\cite{BenKS13} or the other discussed methods for GCAREs, there is no factorization of the GCARE residual with fixed rank
known for ILRSI, e.g., \change{of the form $\cR(Z_{j}Y_{j}^{-1}Z_j^{\trp})=WDW^\trp$ with fixed sized $W,~D$ and coldim$(W)\ll n$.
Of course, one can always use
$\cR(Z_{j}Y_{j}^{-1}Z_j^{\trp})=\hat W\hat D\hat W^{\trp}$ with residual factors
\begin{align*}
\hat W&:=[A^{\trp}Z_j, M^{\trp}Z_j,C^{\trp}]\in\C^{n\times (2j+1)p},\\
\hat D&:=\diag{\left[\begin{smallmatrix}0 & Y_j^{-1}\\ Y_j^{-1}&Y_j^{-1}Z_j^{\trp}BB^{\trp}Z_jY_j^{-1}\end{smallmatrix}\right],I_p}\in\C^{(2j+1)p\times (2j+1)p},
\end{align*}
whose (column) dimensions obviously increase as the iteration proceeds. Computing the residual norm via a factorization like that (see, e.g.,~\cite{BenLP08}) would be quickly getting extraordinarily expensive both memory and computation wise.}
Currently, one of the most reasonable things one can do is to exploit that the spectral norm of a symmetric matrix coincides with its largest eigenvalue and
apply a Lanczos process implicitly to $\cR(Z_{j}Y_j^{-1}Z_j^{\trp})$.
Because of the way the middle low-rank factor $Y_j$ is constructed, computing an approximate feedback matrix $K_j$ requires the complete already computed
low-rank solution. Therefore, just as the projection methods of the previous section, ILRSI does not gain an advantage with respect to memory consumption in the
scenario~\textbf{S2}.
\subsubsection{The RADI iteration}
Assume that the equation~\eqref{gcare} has an approximate solution $X_1 = \Xi_1 \succeq 0$ with a positive semidefinite residual $\cR(X_1)=R_1R_1^{\trp}$.
Then the exact solution $X$ can be written as $X = X_1 + \tilde{X}_1$, where $\tilde{X}_1$ is
a unique stabilizing solution~\cite[Theorem 1]{BenBKetal15} of the residual Riccati equation
\begin{equation}
    \label{eq::residual}
   A_1^{\trp} \tilde{X}_1 + \tilde{X}_1A_1 + C_1^{\trp} C_1 - \tilde{X}_1 BB^{\trp} \tilde{X}_1 = 0,
\end{equation}
with $A_1 = A - BK_1^{\trp}$, $C_1 = R_1^{\trp}$, $K_1=X_1B$ (once more we
suppose $M=I$ for ease of presentation).
This is again a Riccati equation, so one can repeat the procedure, and find an approximate solution $\Xi_2 \succeq 0$ for~\eqref{eq::residual}
such that its residual with respect to this equation is positive semidefinite, and accumulate the approximation $X_2 = X_1 + \Xi_2$, etc.
This way, we obtain an increasing sequence of approximations $0 \preceq X_1 \preceq X_2 \preceq \ldots \preceq X$ to the
solution of the original equation, all of which will have positive semidefinite residuals.

It remains to explain how to construct an approximate solution with a positive semidefinite residual in the $j$th step of the above sequence.
One possibility is to carry out one step of the Cayley iteration~\eqref{eq::CayleyHamiltonIteration} for the residual equation
with the Hamiltonian matrix
\[
  \cH_j =
    \begin{bmatrix}
     A_j & -BB^{\trp} \\
      -C_j^{\trp} C_j & -A_j^{\trp}
    \end{bmatrix}
    =
    \begin{bmatrix}
      A - B K_j^{\trp}  & -BB^{\trp} \\
      -R_j R_j^{\trp} & -{(A - B K_j^{\trp})}^{\trp}
    \end{bmatrix};
\]
here $K_j = X_j B$.
Setting the initial approximation in~\eqref{eq::CayleyHamiltonIteration} to zero, after some calculation we obtain that
the first step approximation $\Xi_{j+1}$ in the Cayley-Hamilton iteration is given by:
\begin{align*}
  \Xi_{j+1} &= -2\Real{\alpha_j} V_{j+1} \tilde{Y}_{j+1}^{-1} V_{j+1}^{\trp}, \\
  V_{j+1} &= {(A -  BK_j^{\trp} + \alpha_j I)}^{-\trp} R_{j}, \\
  \tilde{Y}_{j+1} &= I + (V_{j+1}^{\trp} B){(V_{j+1}^{\trp} B)}^{\trp},
\end{align*}
and that the residual factor and the feedback matrix can be updated as
\begin{align*}
  R_{j+1} &= R_j - 2\Real{\alpha_j} V_{j+1} \tilde{Y}_{j+1}^{-1}, \\
  K_{j+1} &= K_j - 2\Real{\alpha_j} (V_{j+1} \tilde{Y}_{j+1}^{-1}) (V_{j+1}^{\trp} B).
\end{align*}
This describes one step of the RADI procedure. The whole method for generalized CAREs is displayed in Algorithm~\ref{alg:radi}. Note that, unlike the projection
methods, each RADI iteration takes the same amount of time, and does not get progressively slower as the algorithm proceeds: \change{the matrices $V_j$, $K_j$, $R_j$ remain of the same size throughout the iteration}.
\begin{algorithm2e}[t]
\SetEndCharOfAlgoLine{}
\SetKwInOut{Input}{Input}\SetKwInOut{Output}{Output}
\caption{The RADI Iteration for
GCAREs.}\label{alg:radi}
\Input{Matrices $A,~M,~B,~C$ defining~\eqref{gcare}, initial feedback $K_0$}
\Output{$Z_j\in\C^{n\times jp}$, $Y_j=Y_j^{\trp}\in\C^{jp\times jp}$ s.t. $Z_j
Y_j^{-1} Z_j^{\trp}\approx X$, stabilizing feedback matrix
$K_j\in\Cnm$.}
    $R_0 = C^{\trp}$,  $Y_0 = [\;]=[\;]$, $j=1$. 
    \;
    \For{$j=1,2,\ldots,j_{\max}$}{
    Determine shift $\alpha_j\in\C_-$,  $\gamma_j:=\sqrt{-2\Real{\alpha_j}}$.\nllabel{alg:radi_shifts}\;
         Solve  ${(A -  BK_{j-1}^{\trp} +\alpha_j M)}^{\trp} V_j= R_{j-1}$ for $V_j\in\Cnp$.\nllabel{alg:radi_LS}\;
        $V_B=V_j^{\trp} B$, $V_+=\gamma_j^2M^{\trp}V_j$, $\tilde{Y_j} = I +
V_{B}V_{B}^{\trp}$.\;
        If required, update low-rank solution factors $Z_j = [Z_{j-1} \;\; \gamma_{j}V_{j}]$, $Y_j = \begin{sbm} Y_{j-1} & \\ & \tilde{Y_j} \end{sbm}$.\;
        Update Riccati residual factor $R_j = R_{j-1} + (V_+\tilde{Y_j}^{-1})\in\Cnp$.\;
        Update  feedback matrix $K_j = K_{j-1} + (V_+ \tilde{Y_j}^{-1})V_B\in\Cnm$.\;
        \lIf{$\|R_j^{\trp} R_j\|<\tau_{\cR}\|C^{\trp}C\|$}{\textbf{stop}.}
    }
\end{algorithm2e}

The shift parameters $\alpha_j$ are again important for a rapid convergence. We mention and later on employ an approach proposed in~\cite{BujB14,BenBKetal15}:
Let $M=I$, $U_{\ell}\in\R^{n\times \ell p}$,  $U_{\ell}^{\trp}U_{\ell}=I_{\ell}$, $\ell p\ll jp$ with $\range{U_{\ell}}\subset\range{Z_j}$
and
\begin{align*}
 \cH_j^{\text{proj}} =
  \begin{sbm}
    U_{\ell}^{\trp}A_j U_{\ell}                                              & -(U_{\ell}^{\trp} B) {(U_{\ell}^{\trp} B)}^{\trp} \\
    -(U_{\ell}^{\trp} R_j^{\trp}) {(U_{\ell}^{\trp} R_j^{\trp})}^{\trp} & -U_{\ell}^{\trp}A^{\trp}_j U_{\ell}
  \end{sbm}\in\R^{2\ell\times 2\ell}
\end{align*}
be the Hamiltonian matrix associated to the residual-CARE projected onto $\range{U_{\ell}}$.
If $\theta_i, q_i=\begin{sbm} q^{(1)}_i \\ q^{(2)}_i\end{sbm}$, with $q^{(1,2)}\in\C^{\ell}$, are the eigenpairs of $\cH_j^{\text{proj}}$, then the
next shift $\alpha_{j+1}$ is
selected as the eigenvalue $\theta_i$ with the largest value of $\|q_i^{(2)}\|$. The motivation behind this strategy, along with further details, can be found
in~\cite{BujB14,BenBKetal15}. In the general case $M\neq I$, the eigenpairs of the corresponding (projected) Hamiltonian/skew-Hamiltonian pencil have to be
considered. The number of block columns $\ell$ taken into account is typically chosen very small, e.g., one takes the last $2p$ columns of $Z_j$.
This
shift selection strategy appears to be efficient and, judging by the numerical experiments in~\cite{BenBKetal15}, often superior to other approaches.
Pairs of complex conjugated shifts can be dealt with appropriately~\cite{BenBKetal15}, so that only an absolutely necessary amount of complex arithmetic
operations remains and, e.g., the $Z_j, Y_j, K_j$ will be real matrices.

In contrast to ILRSI, because of the block-diagonal structure of the constructed middle low-rank factor $Y_j$, it is possible to accumulate only the feedback
approximation $K_j$ in the RADI iteration without ever forming the low-rank factors $Z_j, Y_j$. This makes the method a good memory-efficient candidate for
scenario~\textbf{S2}.
\subsubsection{Equivalences and relations to other methods}
As shown in~\cite[Theorem~2]{BenBKetal15}, if executed with the same set of shift parameters, ILRSI and the RADI iteration both produce the same output and are,
in
fact, equivalent. Hence, both methods are simply different implementations of the same approximation sequence. Moreover, if $B=0$, both methods reduce to the
LR-ADI iteration for Lyapunov equations. Because of
this close relation, we will also execute ILRSI with the residual Hamiltonian shifts described above.
ILRSI and the RADI  iteration are in the same sense also equivalent to the quadratic ADI iteration~\cite{morWonB07}.
We will not incorporate the quadratic ADI iteration into this study since it has been shown to be inferior in terms of numerical efficiency to the
considered methods, see, e.g.,~\cite{BenBKetal15,BujB14}.

Although no Galerkin projection is used within ILRSI, RADI, it can still be shown that the computed low-rank solution factors $Z_j$ span a rational Krylov
subspace. In fact, since in ILRSI the construction of the columns in $Z_j$ is very close to the LR-ADI iteration for Lyapunov equations, it can be shown that
$\range{Z_j}=
\range{\left[{(A-\alpha_1M)}^{-\trp}C^{\trp},\ldots,\left[\prod\limits_{i=2}^j(A-\alpha_iM)^{-1}
\right] (A-\alpha_1 M)^{-1}C^{\trp}\right]}$,\\ see, e.g., \cite{morLi00,LiW02,morWol15}. Further connections of ILRSI to the rational Krylov framework are discussed
in~\cite{MasOR16}. By the above equivalence, those connections hold for RADI as well.

\subsection{Low-rank Newton Methods}
Motivated by the nonlinear nature of $\cR(X)$, Newton schemes are applied to solve~\eqref{gcare}.
With a given initial guess $X_0$, step $k$ of a Newton iteration for~\eqref{gcare} is,
following, e.g.,~\cite{ArnL84,LanR95}, given by
  \begin{align}\label{newton_ARE_step}
    \cR'|_{X_{k-1}}(N_{k})&=-\cR(X_{k-1}),\quad X_k=X_{k-1}+N_{k}.
  \end{align}
  Here, $\cR'|_{X}$ denotes the Fr\'{e}chet derivative of $\cR$ at $X$ and is given by
  \begin{align*}
    \cR'|_{X}~:~N \mapsto {(A-BB^{\trp}XM)}^{\trp}NM+M^{\trp}N(A-BB^{\trp}XM).
  \end{align*}
  Hence, the update $N_{k}$ is the solution of a generalized, algebraic,
  continuous-time Lyapunov equation (GCALE).
Kleinman~\cite{Kle68} showed that it is possible to reformulate the Newton step
to directly compute $X_k$ instead of working with the correction $N_k$. In that case, the iteration and the corresponding GCALE to be solved is
\begin{subequations}\label{newtonNK_ARE}
\begin{align}\label{newtonNK_ARE_GCALE}
  A_k^{\trp}X_{k}M+M^{\trp}X_{k}A_k&=-F_{k}F_k^{\trp},\quad F_k:=[C^{\trp},K_{k-1}],\quad A_k:=A-BK_{k-1}^{\trp}\\\label{newtonNK_ARE_rhs}
  \text{with}\quad K_{k-1}&:=M^{\trp}X_{k-1}B.
\end{align}
\end{subequations}
Under
the assumed properties of this paper and provided $X_0$ is a
stabilizing initial guess, the sequence $\lbrace X_k\rbrace_{k\geq1}$ converges quadratically towards $X$, the exact solution of~\eqref{gcare}.
Moreover, it converges monotonically decreasing with respect to the L\"owner ordering, i.e.~$X_1\succeq X_2\succeq\ldots\succeq X$, and for all $k>0$ it holds
that
$\Lambda(A-BB^{\trp}X_{k}M,M)\subset\C_-$.
Because of the size restriction on $B$, $C$ we set in the introduction, the right hand side
$F_k$ in~\eqref{newtonNK_ARE_GCALE} is (at most) of rank $p+m\ll n$.
We can, therefore, utilize low-rank algorithms for large-scale
GCALEs~\cite{Sim16,BenS13} that compute low-rank solution factors of the solution $X_k$ in~\eqref{newtonNK_ARE_GCALE}, e.g., $X_k\approx Z_k Z_k^{\trp}$. Since
these
algorithms are typically of an iterative nature, we end up with an inner-outer
iteration, consisting of the inner particular GCALE solver iteration, and the
outer Newton iteration.
To distinguish these two stages,
subscripts $_k$ and bracketed superscripts $^{(j)}$ will refer to
data associated to
the outer (\(_{k}\)) and inner (\(^{(j)}\)) iteration. In many papers~\cite{BenLP08,BenS08,FeiHS09,BenS10,BenS13,BenHSetal16}, the low-rank alternating
directions
implicit (LR-ADI) iteration~\cite{Pen99,LiW02,BenKS13a} is employed
for the purpose of solving the GCALE~\eqref{newtonNK_ARE_GCALE}. For a fixed $k\geq1$, the LR-ADI iteration in
the most recent formulation~\cite{BenKS13a,Kue16} proceeds for $j\geq1$ in the following form:
\begin{align}\label{glradi4NK}
\begin{split}
 V_k^{(j)}&={(A_k+\alpha^{(j)}_{k}M)}^{-\trp}W_k^{(j-1)},\quad
W_k^{(j)}=W_k^{(j-1)}-2\Real{\alpha^{(j)}_k}M^{\trp}V_k^{(j)},\\
 Z_k^{(j)}&=\left[Z_k^{(j-1)},~\sqrt{-2\Real{\alpha^{(j)}_k}}V_k^{(j)}\right],
 \end{split}
\end{align}
where $W_k^{(0)}:=F_k$, $Z_k^{(0)}=[~]$, and $\alpha^{(j)}_k\in\C_-$ are
shift parameters.
Using~\eqref{glradi4NK} as the inner Lyapunov solver within the Newton-Kleinman iteration~\eqref{newtonNK_ARE} leads to the low-rank Newton-Kleinman ADI
(NK-ADI) method, which is illustrated in Algorithm~\ref{alg:lrnm_adi}.

\begin{algorithm2e}[t]
\SetEndCharOfAlgoLine{}
\SetKwInOut{Input}{Input}\SetKwInOut{Output}{Output}
  \caption{Low-rank Newton-Kleinman ADI for GCAREs (NK-ADI)}\label{alg:lrnm_adi}
    \Input{Matrices $A,~M,~B,~C$ defining~\eqref{gcare}, initial feedback
$K_0$, and stopping tolerances $0<\tau_{\text{ADI}},\tau_{\cR}\ll 1$.}
    \Output{$Z_k\in\C^{n\times (m+p)j}$ such that
    $Z_{k}Z_{k}^{\trp}\approx X$, stabilizing feedback matrix
$K_k\in\Cnm$.}
\For{$k=1,\ldots,k_{\max}$}{%
$W_k^{(0)}=[C^{\trp},K_{k-1}]$,
$Z^{(0)}_k=[~]$, $K^{(0)}_k=0$,~$j=0$.\;
\While{$\|W_k^{(j)}\|^2>\tau_{\text{ADI}}\|W_k^{(0)}\|^2$}
{\nllabel{alg:lrnm_adi:stop_inner} %
$j=j+1$\;
Determine next shift $\alpha_k^{(j)}$.\nllabel{alg:lrnm_adi:shifts}\;
Solve ${(A_k+\alpha^{(j)}_{k}M)}^{\trp} V_k^{(j)}=W_k^{(j-1)}$
for $V_k^{(j)}$.\nllabel{alg:lrnm_adi:linsys}\;
$\gamma_k^{(j)}=\sqrt{-2\Real{\alpha_k^{(j)}}}$,
$V_+:={(\gamma_k^{(j)})}^2M^{\trp}V_k^{(j)}$.\;
Update factor of Lyapunov residual $W_k^{(j)}=W_k^{(j-1)}+V_+$.\;
Update low-rank solution factor
$Z_k^{(j)}=[Z_k^{(j-1)},~\gamma_k^{(j)}V_k^{(j)}]$ if
required.\nllabel{alg:lrnm_adi:lrf}\;
Implicit update of feedback matrix $K^{(j)}_k=K^{(j-1)}_k+V_{+}{(V^{(j)}_{k})}^{\trp}B$.
\nllabel{alg:lrnm_adi:feedbackR}\;
}
$K_k=K_k^{(j)}$, $Z_k=Z^{(j)}_k$.\;
\lIf{$\|\cR(Z_{k}Z_{k}^{\trp})\|<\tau_{\cR}\|C^{\trp}C\|$}{\textbf{stop}.}
}
\end{algorithm2e}
The shift parameters $\alpha^{(j)}_k$ in Line~\ref{alg:lrnm_adi:shifts} are,
similarly to those in RKSM, ILRSI, and RADI, crucial for a fast reduction of the error.
Here, we employ, without going into detail, the automatic shift
generation strategy proposed in~\cite{BenKS14}: first, several shifts are selected from the spectrum of the matrix pair
$(U_{\ell}^{\trp}A_{k}U_{\ell},~U_{\ell}^{\trp}MU_{\ell})$,
where $U_{\ell}\in\R^{n\times \ell p}$,  $U_{\ell}^{\trp}U_{\ell}=I_{\ell}$, $\ell< j$, with $\range{U_{\ell}}\subset\range{Z_k^{(j)}}$.
Once these shifts are depleted, the procedure is repeated.
This shift selection strategy usually leads to the
best performance of the LR-ADI iteration, both in terms of execution time and the
required number of iteration steps. Pairs of complex shifts can be handled as in the LR-ADI~\cite{BenKS13,BenKS13a}.
 It can be shown~\cite{BenS13,BenHSetal16,Kue16} that, inherited from the LR-ADI iteration~\cite{BenKS13a}, the GCARE residual at inner step $j$ and outer
step $k$ is given by
\begin{align*}
  \cR_k^{(j)}&:=\cR(Z_k^{(j)}{(Z_k^{(j)})}^{\trp})=[W_k^{(j)},~\Delta
  K_k^{(j)}]\diag{I_{p+m},-I_m}{[W_k^{(j)},~\Delta K_k^{(j)}]}^{\trp},
\end{align*}
where $\Delta K_k^{(j)}:=K_k^{(j)}-K_{k-1}\in\C^{n\times m}$ and $W_k^{(j)}\in\C^{n\times (p+m)}$ is the explicit low-rank factor of the current GCALE residual.
Hence, computing the GCARE residual norm requires only the computation of the spectral norm of a thin $n\times (p+2m)$ matrix.

Algorithm~\ref{alg:lrnm_adi} reveals another advantage of using the LR-ADI iteration to solve the GCALE~\eqref{newtonNK_ARE_GCALE}. Similar to the RADI
algorithm, it is not necessary to store the low-rank factor $Z_k^{(j)}$
(Line~\ref{alg:lrnm_adi:lrf}) if only the feedback matrix $K$ is of
interest, since $K_k^{(j)}$ can be updated incrementally (cf. Line~\ref{alg:lrnm_adi:feedbackR},~\cite{BenLP08}).
Thus, the NK-ADI method is also very memory efficient in scenario \textbf{S2}. Note that, in order to execute the proposed shift strategy, a
small number of additional columns needs to be stored to carry out the projection of $(A_k,~M)$.

Naturally, one could also use a projection method like EKSM or RKSM~\cite{Sim07,DruS11} for solving the Lyapunov equation~\eqref{newtonNK_ARE_GCALE}, which
leads to other interesting results~\cite{SimSM13a}, e.g., regarding the Riccati residual.
However,
doing so would sacrifice the advantage in scenario \textbf{S2} since, just as with the projection methods in Section~\ref{ssec:projectionmethods}, the matrices
$Z_k^{(j)}$ would need to be stored.

Regardless of the low-rank method employed to solve~\eqref{newtonNK_ARE_GCALE}, the low-rank Newton method in the presented form often
shows a comparatively slow convergence towards an approximate stabilizing solution. One intuitive explanation is that the method updates the feedback gain
approximations only after an outer iteration is completed.
In practice, one of the following improvements is therefore
usually mandatory to make the low-rank Newton method competitive.
\subsubsection{Galerkin acceleration of the outer iteration}
Similar to RADI, ILRSI, it holds that $\range{Z^{(j)}_k}=\cK^{\text{rat}}_j(A_k,M,[C^{\trp},K_{k-1}],\boldsymbol{\alpha}_k)$, see
\cite{morLi00,LiW02,morWol15}, so
in order to improve the convergence of the Newton iteration,~\cite{Saa09,BenS10}
suggest to perform a Galerkin projection onto $\range{Z_k}$ in each
outer iteration step, just after the GCALE~\eqref{newtonNK_ARE_GCALE} has been solved.
Assume the columns of $Q_k\in\C^{n\times j(p+m)}$ with $Q_k^{\trp}Q_k=I$ constitute an orthonormal basis for $\range{Z_k}$.
The projection is performed similar to the methods in Section~\ref{ssec:projectionmethods}:
a small-scale, at most $j(p+m)$ dimensional GCARE
\begin{align*}
 \tA_k^{\trp} Y_k\tM_k+\tM_k^{\trp} Y_k\tA_k^{\trp}-\tM_k^{\trp} Y_k\tB_k\tB_k^{\trp} Y_k\tM_k+\tC_k^{\trp}\tC_k=0
\end{align*}
with $\tA_k:=Q_k^{\trp} AQ_k$, $\tB_k=Q_k^{\trp} B$, $\tC_k=CQ_k$, $\tM_k:=Q_k^{\trp} MQ_k$
is solved for its stabilizing solution $Y_k$.
\change{The restriction $\tA_k$ can be constructed without matrix vector
multiplications involving $A$. The precise formula depend on the choice of the low-rank Lyapunov solver, e.g., the LR-ADI iteration~\cite{Kue16} or projection methods~\cite{DruS11,SimSM13a}.}
The approximate solution of~\eqref{gcare} associated to the Galerkin
acceleration is $X_k^{\text{pr}}=Q_{k}Y_{k}Q_k^{\trp}$. The next Newton-Kleinman
iteration step
 then continues with the updated feedback matrix $K^{\text{pr}}_{k+1}:=M^{\trp}X_k^{\text{pr}}B=M^{\trp}Q_{k}Y_{k}\tB_{k}$.
The GCARE residual $\cR(X^{\text{pr}})$ after this projection does not have the low-rank structure from above and, thus,
 $\|\cR(X^{\text{pr}}_k)\|$ has to be computed
differently, e.g., via applying a Lanczos process to
$\cR(X^{\text{pr}}_{k})$ similar to ILRSI.\@ Moreover, in order to perform this Galerkin projection, the whole low-rank
solution factor $Z_k$ has to be stored, so this variant of the NK-ADI method loses its ability to
operate solely on feedback approximations and the advantage for scenario \textbf{S2}.
However, numerical experiments~\cite{Saa09,BenS10,BenKS16,Kue16} show that this
projection tremendously
accelerates the outer Newton-Kleinman iteration. Often, it decreases
the number of required Newton steps down to one or two. We will abbreviate this approach by NK-ADI+GP.\@
\begin{remark}\label{rem:lradi_gcare_GPinner}
In~\cite{Saa09,BenS10} it is also suggested to accelerate the inner low-rank ADI
iteration of Algorithm~\ref{alg:lrnm_adi} by a Galerkin projection framework. In
view of the improvements in the LR-ADI formulation, we strongly refrain from
this idea. On the one hand, using a Galerkin projection
within the LR-ADI iteration destroys important structural properties of the method, e.g.,
the associated GCALE residual does no longer have  low-rank
factorization and it is then not clear what should be used as $W_k^{(j)}$
in~\eqref{glradi4NK}. On the other hand, one could have used a
projection method as inner iteration from the start as discussed in~\cite{SimSM13a}.
\end{remark}
 \subsubsection{Inexact GCALE solves and line-search}
Newton methods are, under some mild conditions, still able to converge to the desired solution if the Newton step is carried out
inexactly, see, e.g.,~\cite{DemES82,EisW96}. For the Newton-Kleinman iteration~\eqref{newtonNK_ARE} for GCAREs, this means that the
GCALE~\eqref{newtonNK_ARE_GCALE} is only solved inexactly, for example, the LR-ADI iteration might be terminated once it
satisfies $\|A_k^{\trp}X_{k}M+{M}^{\trp}X_{k}A_k+F_{k}F_k^{\trp}\|\leq\tau_{\text{ADI}}$ for an appropriately chosen $\tau_{\text{ADI}}>0$, e.g.,
$\tau_{\text{ADI}}=\eta\|\cR_{k-1}\|$, $\eta\in\mathopen{(}0,1\mathclose{)}$. These inexact solves, obviously, have the potential to drastically reduce the
number of required LR-ADI
steps. At the same time, convergence towards a stabilizing solution of the GCARE~\eqref{gcare} has to be maintained, which makes the analysis
difficult~\cite{FeiHS09,Hyl11,BenHSetal16,Wei16}.
In~\cite{BenHSetal16,Wei16}, a novel theoretical and numerical framework for the inexact
Newton-Kleinman iteration for~\eqref{gcare}
is proposed which additionally  incorporates a line search strategy~\cite{BenB98}.
At
Newton step $k$, let $K_k$ be a feedback approximation linked to the approximate or exact solution
$X_k^{(j)}=Z_k Z_k^{\trp}$ of the GCALE~\eqref{newtonNK_ARE_GCALE}.
Then the next Newton step is carried out with an improved feedback approximation
$\hK_k=(1-\beta)K_{k-1}+\beta K_k$ for an appropriately chosen \(\beta\in \mathopen{(}0,1\mathclose{]}\).
The numerical effort to carry out this line search strategy is
negligible and we refer to~\cite{BenHSetal16,Wei16} for details about its
practical implementation, including possible choices of $\eta,~\beta$.
Compared to the plain NK-ADI iteration (Algorithm~\ref{alg:lrnm_adi}), the
number of both inner and outer iteration steps, and consequently the number of
arising linear systems, can be drastically reduced by this approach.
The inexact NK-ADI method equipped with this line-search strategy (abbreviated
by iNK-ADI+LS) is still able to work only on the feedback approximations $K_k$,
preserving the advantage of NK-ADI in scenario~\textbf{S2}.
\subsection{Related and further methods}
Apart from the methods described so far a number of other methods for solving large-scale GCAREs can
also be found in the literature.

The structured doubling algorithm (SDA)~\cite{ChuFL05,Bin14} is a recent method
based on efficient computation of deflating subspaces by doubling. The original
formulation of the method is suited for small-scale dense CAREs, for which it
performs very well. The SDA has also been adapted to the large-scale
setting~\cite{ChuFL05,LiCLetal13}. However, our experience with a number of
numerical experiments indicates that the large-scale algorithm is not yet
competitive to the methods described in this paper. Just as an illustration, we
included the results obtained by the large-scale SDA algorithm in
Example~\ref{ex:toep} of Section~\ref{sec:numex}. We have therefore chosen to
omit this method in the detailed analysis.

Another class of methods are those that compute eigenvectors of the Hamiltonian matrix $\cH$ associated with the stable eigenvalues~\cite{AmoB10,BujB14}. As
mentioned in Section~\ref{sec:nonprojectiveiterations}, the exact solution of the CARE is $X = QP^{-1}$, where ${[P^{\trp},Q^{\trp}]}^{\trp}\in\R^{2n\times n}$
spans the $\cH$-invariant stable subspace. Computing the entire stable subspace of a large matrix is not feasible, so these methods target a small number $\ell
\ll n$ of eigenvectors, and use only these to approximate $X$ by means of carefully designed formulas~\cite{AmoB10}. The main issue of this
approach is that it is quite difficult to determine which eigenvectors should be targeted~\cite{BujB14}. This difficulty renders the eigenvector-based
methods non-competitive as well.
On the other hand, these methods have motivated the development of both ILRSI/RADI.\@ The original derivation
of RADI~\cite{BenBKetal15} uses the formula from~\cite{AmoB10} to find an approximate solution $\Xi_1$ to~\eqref{eq::residual}; the goal of ILRSI is
computation of the stable subspace of $\cH$ as well, although it is achieved by different means. For $p=1$, ILRSI/RADI can be considered a generalization of the
methods from~\cite{AmoB10,BujB14}, since it is equivalent to these methods when eigenvalues of $\cH$ are used as shifts~\cite{BenBKetal15}.
\subsection{Unstable GCAREs}\label{ssec:unstable}
In the literature regarding large-scale GCAREs, often only the stable case is
used to test the numerical methods. We briefly discuss the situation where the
GCARE~\eqref{gcare} is defined by an unstable matrix pencil $A-\lambda M$. Recall that
the overall system $(A,M,B)$ has to be stabilizable. One way of handling this
situation is to provide a stabilizing initial guess $X_0=X_0^{\trp}\in\Rnn$,
meaning that the pair $(A_{K_0}:=A-BB^{\trp}X_0^{\trp}M,M)$ is stable.
 Also, depending on how serious the
instability of $A-\lambda M$ is, some of the investigated methods might be able
to converge without an initial guess, see~\cite{Sim16a} for results regarding
projection methods.  Of course, providing an initial guess can also help to
speed up the whole method, independent of the stability of the pencil
$A-\lambda M$, see~\cite{SimSM13a} for some numerical tests on this.

Using an additive decomposition  $X=X_++X_0$ of the sought after stabilizing solution of~\eqref{gcare}, it is well known~\cite{MehT88,BenBKetal15a} that the
increment $X_+=X_+^{\trp}\in\Rnn$ is the stabilizing solution of the stable GCARE
\begin{align}\label{gcare_stabilized}
 A_{K_0}^{\trp}X_+M+M^{\trp}X_+A_{K_0}-M^{\trp}X_+BB^{\trp}X_+M+\cR(X_0)=0.%
\end{align}
The stabilized GCARE~\eqref{gcare_stabilized} can now be dealt with by using any of the Riccati methods described so far, providing
$A_{K_0}$, $M$, $B$, and $\cR(X_0)$ as the new input.
In the case of RKSM, EKSM and ILRSI, the altered structure of the input matrices
leads to an increase in the algorithm complexities, especially regarding the
linear systems solves.  In contrast to that, the RADI iteration and the low-rank
Newton methods can be employed to the original GCARE~\eqref{gcare} right away if
the initial feedback $K_0$ is given to them as well.

Strategies for computing an initial feedback $K_0$ with stabilizing property can
be found, e.g., in~\cite{Rus79, BanI91, LanR95, Sim96, Hei09, AmoB12, Bin14,
  BaeBSetal15, BenHSetal16, Wei16}. In the remainder, we only use the approach
from~\cite{AmoB12,BaeBSetal15}, which we briefly describe next.

For $M=I$, let $A^{\trp}Q_u=Q_{u}R_{u}$ be a partial real Schur decomposition of
$A^{\trp}$ such that $\Lambda(R_{u})=\Lambda(A)\cap \C_+$, i.e., the columns of
$Q_u$ form an orthonormal basis of the associated unstable invariant
subspace. An initial guess is then defined with $X_0:=Q_{u}S_{u}Q_{u}^{\trp}$
and $K_{0}:=X_{0}B$, where $S_{u}$ solves the algebraic Bernoulli equation
$R_{u}S_{u}+S_{u}R^{\trp}_{u}-S_{u}(Q_{u}^{\trp}B)(B^{\trp}Q_{u})S_{u}=0$. Assuming
that the number of unstable eigenvalues is small ($u\ll n$), this Bernoulli
equation can be solved by standard, dense methods~\cite{BarBQ07}. As an
alternative, one can, under some mild conditions~\cite{AmoB12}, use
$X_0:=Q_{u}T^{-1}_{u}Q_{u}^{\trp}$, $K_{0}=Q_{u}T^{-1}_{u}Q_{u}^{\trp}B$ with
$T_{u}$ being the solution of the small Lyapunov equation
$-R_{u}T_{u}-T_{u}R^{\trp}_{u}+(Q_{u}^{\trp}B)(B^{\trp}Q_{u})=0$. This equation
is potentially easier to solve than the Bernoulli equation.  As an important
side effect, the initial Riccati residual in~\eqref{gcare_stabilized} satisfies
$\cR(X_{0})=C^{\trp}C$ with both of these choices. It is noteworthy that with this selection of an initial guess, it is in scenario~\textbf{S2} sufficient
to provide just the feedback matrix $K_0:=M^{\trp}X_0B\in\Rnm$ instead of $X_0$.

\begin{remark}
  Originally, the left and right eigenvectors of $A$ corresponding to
  $\Lambda(A)\cap \C_+$ were used in the computation of
  $X_0$~\cite{Hei09,AmoB12}; this still leads to the same initial guess as
  above. We prefer the usage of invariant subspaces since computing a partial
  Schur decomposition of a large matrix $A$ is often easier than computing left
  and right eigenvectors.

  Nevertheless, computing the entire unstable invariant subspace and,
  especially, ensuring that no unstable eigenvalues are missed can in practice
  be a very demanding task. Large-scale eigenvalue methods for this purpose can
  be found e.g., in~\cite{BaiDDetal00} and the references therein.
\end{remark}
In the numerical experiments, whenever unstable CAREs are considered, we assume
for the sake of simplicity that the required initial feedback is provided and do
not consider the numerical effort and difficulties of its computation. The
later, however, is nevertheless a crucial subject to be addressed in future
research, where, e.g., strategies from general Newton and homotopy methods for an initial guess selection  might be adopted~\cite{All03}.

\section{Comparison of the Main Computational Stages}\label{sec:compare}
In this section, we analyze the major computational subtasks that arise in the
discussed algorithms: solving large linear systems of equations with multiple
right hand sides, computing shift parameters, building orthogonal bases and
solving projected CAREs, as well as computing, or estimating, the CARE residual
norm.
\subsection{Solving Linear Systems}\label{ssec:linsys}
The most prominent feature of all considered algorithms is the solution of a
large linear system of equations with possibly multiple right hand sides.  We,
here, restricted ourselves to sparse direct solution strategies, because, on the
one hand, this worked sufficiently well in our experiments.  On the other hand,
due of the sheer number of different iterative solvers and preconditioning
strategies, a thorough discussion of preconditioned iterative solves would
clutter the presentation and is, therefore, beyond the scope of this
study. Also, the occurrence of multiple right hand sides further increases the
number of
available iterative solution approaches.  Moreover, the effect of errors arising
in the solution of linear systems on the methods for large matrix equations needs to be considered in the case of iterative linear
solvers. Some research on this topic for large Lyapunov equations can be found
in~\cite{Sab07,Sun08,KueF18}.

An overview of the structure of the arising systems in stable and unstable
situations, as well as the number of columns in the right hand sides is given in
Table~\ref{tab:linsys}. In several GCARE methods, we notice the occurrence of
systems defined by the sum of a sparse matrix $A+\alpha_j M$ and a low-rank
term $BK_j^\trp$ of rank at most $m$. If sparse-direct solvers are to be applied, such systems can
be dealt with by the Sherman-Morrison-Woodbury (SMW) formula~\cite{GolV13} via
\begin{align*}
    & y_j={(A - BK_j^\trp  +\alpha_j M)}^{-\trp} b_j = w_j+g_j{(I_m-B^{\trp}
      g_j)}^{-1} B^{\trp}  w_j,\\
    &\text{where }[w_j,g_j]\text{ solves }{(A+\alpha_j
      M)}^{\trp}[w_j,g_j]=[b_j,K_j].
\end{align*}
Hence, solving a sparse-plus-low-rank system with $s$ right hand sides is
expressed as solving a sparse system with $A+\alpha_j M$ and $s+m$ right hand
sides.

Among the algorithms, EKSM has the advantage that the coefficients of the
arising linear systems do not change during the iteration. Hence, a single
pre-computed sparse LU factorization can be reused throughout the whole
iteration by means of often significantly cheaper sparse
triangular solves.
This fact can also be exploited in the presence of an initial feedback $K_0$, by
using the SMW formula once again:
\begin{align*}
  y_j={(A-BK_0^\trp )}^{-\trp} b_j&=w_j+g_0{(I_m-B^\trp  g_0)}^{-1}B^\trp  w_j,
\end{align*}
where $w_j$, $g_0$ solve $A^\trp w_j=b_j$, $A^\trp g_0=K_0$ and also
$g_0\in\Rnm$ has to be computed only once.  Thus, in the unstable situation,
except for this extra linear system for $g_0$, no major additional work is
required in EKSM in contrast to the other methods.

For generalized Riccati equations, Table~\ref{tab:linsys} considers only the
case when the equivalent CARE defined by $M^{-\trp} A^\trp $,
$M^{-\trp} C^\trp $ is used implicitly in the projection methods, as explained
in Section~\ref{ssec:projectionmethods}. This means that an extra initial linear
solve with $M^\trp $ is needed in RKSM to compute $M^{-\trp} C^\trp $, whereas
EKSM requires such an additional solve with $M^\trp $ in every iteration step. Using sparse
Cholesky factors of $M=M^\trp \succ 0$ would replace a solve with $M^\trp $ by
two sparse triangular solves.  Regarding the number of right hand sides in the
occurring linear system, Table~\ref{tab:linsys} reveals that low-rank NK
iterations have the largest number ($p+m$) among all algorithms. If the SMW
formula is used, linear systems with the matrix $A+\alpha_j M$ and $p+2m$ right
hand sides have to be solved in each iteration step, in contrast to $p+m$ or
only $p$ in the other methods. Consequently, among all considered methods,
solving linear systems is most expensive in low-rank NK methods.  If no initial feedback is given or required, RKSM and ILRSI have an advantage over RADI
since they do not need the SMW
formula to solve systems with the matrix $A+\alpha_j M$ and $p$ right hand sides
 (rhs).
\begin{table}[t]
  \centering
  \caption{Form of the occurring linear system in each iteration
    step. Subscripts $_j$ indicate quantities varying during the iterations.
        \change{We also provide the number of columns in the
                  right hand side of the original linear systems and
                  the SMW version.
        }
        }
  \begin{tabularx}{\linewidth}{|l|X|X|l|l|}
    \hline
    &\multicolumn{2}{c|}{Coefficients of arising linear system{(s)}}&\multicolumn{2}{c|}{columns $s$ of rhs}\\
    Method&stable (no initial $K_0$)&with initial $K_0$&original&SMW\\
    \hline\hline
    {(G)}EKSM&$M^\trp$, $A^\trp$&$M^\trp$, ${(A-BK_0^\trp )}^{\trp}$&$p$&$p$\\ \hline
    RKSM&\multirow{3}{*}{${(A+\alpha_j M)}^{\trp}$}&
    \multirow{3}{*}{${(A-BK_0^\trp +\alpha_j M)}^{\trp}$}&$p$&$p+m$\\
    TRKSM&&&$p_j\leq p$&$p_j+m$\\
    ILRSI&&&$p$&$p+m$\\\hline
    RADI&\multicolumn{2}{c|}{\multirow{2}{*}{${(A-BK_j^\trp +\alpha_j M)}^{\trp}$}}&$p$&$p+m$\\
    NK-ADI&\multicolumn{2}{c|}{}&$p+m$&$p+2m$\\
    \hline
  \end{tabularx}\label{tab:linsys}
\end{table}

\subsection{Shifts}
With the exception of EKSM, all considered methods require shift parameters to
achieve a rapid convergence.  Initially, these parameters were computed in
advance, prior to any iteration of the Riccati solver.  In the last years,
dynamic shift generation
strategies~\cite{DruS11,DruSZ14,BenKS14,Sim16a,BenBKetal15,Kue16} have
attracted increasing attention.  In these strategies, the shift needed in a
particular iteration step is computed by using all the data that is available to
the method in that step.  This typically allows the methods to achieve
convergence in a smaller number of iterations and with less user interaction,
but may increase the computational complexity of each step.

A substantial amount of work in all the adaptive shift generation strategies
considered here, is the solution of an eigenvalue problem, whenever new shifts
are needed. These eigenvalue problems differ in their size and structure,
depending on the GCARE method and the particular selection strategy. In
Table~\ref{tab:shift_evp} we summarize the defining matrices of the eigenvalue
problems, together with their sizes in the strategies used in this study. As
before,
$j$ indicates the current iteration number. Although the automatic shift
generation of RKSM involves the largest eigenvalue problem, the costs are often
still negligible since $jp\ll n$. The dimension of the eigenvalue problem of the
approaches in ILRSI, RADI, and NK-ADI depends on a pre-specified small number
$\ell$, for which typical values in the literature are in the range
$1\leq\ell \leq 6$.  In the projection methods (T)RKSM the involved matrices
$\tilde{A}_j$ are already constructed as a part of the main iteration loop,
while the other methods additionally need to construct them, as well as the
orthogonal bases $U$ for certain subspaces just for the purpose of generating
the shifts.  Note that the approach used in NK-ADI will return more than one
shift parameter, so it does not have to be executed in every iteration step.

\begin{table}[t]
  \centering
  \caption{Matrices of the arising eigenvalue problems for the dynamic shift
    generation.}
  \begin{tabularx}{\linewidth}{|l|X|l|}
\hline
    Method&matrices&size\\
    \hline
    RKSM&$\tA_j$&$jp$\\
    TRKSM&$\tA_j$&$\sum\limits_{i=1}^jp_i$\\
    ILRSI, RADI&$\cH_j^{\text{proj}},~\cM_j:=\diag{U^\trp M U,U^\trp M^\trp
                 U}$&$2\ell p$\\
    NK-ADI&$U^\trp A_k U,U^\trp M U$&$\ell (p+m)$\\
    \hline
  \end{tabularx}\label{tab:shift_evp}
\end{table}
\subsection{Building orthonormal bases of the projection spaces}\label{sec:orth}
The projection based methods {(G)}EKSM, {(T)}RKSM, and the Galerkin accelerated
NK-ADI iteration require the construction of orthonormal bases of the used
subspaces. Let $Q_{j-1}=[q_1,\ldots,q_{j-1}]$ denote the orthonormal subspace
basis after step $j-1$ in either of these methods, and let $q_+$ denote a block
of vectors that need to be added to the subspace in step $j$. In {(G)}EKSM,
$q_+$ and each of the $q_i$ have $2p$ columns, while in RKSM they have $p$
columns. In TRKSM, the column dimensions of the block $q_i$ is $p_i$, where
$p_i\leq p$ is the number of tangential directions that were used in iteration
$i$.  We employ the modified blockwise Gram-Schmidt process to expand the basis
with $q_+$:
\begin{align*}
  \gamma=q_i^{\trp}q_+,~h_j = h_j+ \gamma,~q_+ = q_+ - q_i\gamma,~i=1,\ldots,j-1,
\end{align*}
and this is done twice (initially, $0:=h_j\in\R^{d_{j-1}\times d_+}$). The new
orthonormal basis block $q_j$ is then obtained by a thin QR-factorization
$q_+=q_{j}h_+$, where $h_+\in\R^{d_+\times d_+}$ contains the orthogonalization
coefficients.  Note that GEKSM uses a different inner product in the
orthogonalization process.  If complex shifts and directions occur in (T)RKSM,
$Q_{j-1}$ is first orthogonally expanded by the real and then by the imaginary
part of the associated complex solution vectors of the linear system.  TRKSM in
the form proposed in~\cite{DruSZ14} requires the construction of a second
orthonormal basis for shift generation and residual norm estimation.  In the
\MMESS~implementation~\cite{SaaKB16-mmess-1.0.1} of NK-ADI+GP, the Galerkin
projection is performed only after a Newton step, and the basis for the
projection is orthonormalized using the \matlab{}  routine \texttt{orth}.
\subsection{Small-scale CARE solution}\label{sec:smallCARE}
The projection based methods {(G)}EKSM, {(T)}RKSM, and NK-ADI+GP have to solve
small, dense Riccati equations in some stages of the algorithms.  Assuming $M=I$
for simplicity, and denoting by $\tA_j$, $\tB_j$, $\tC_j$ the projected matrices
$A$, $B$, $C$, respectively, the existence of a stabilizing solution of this
small CARE is ensured if $(\tA_j$, $\tB_j)$ and $(\tA^{\trp}_j$, $\tC^{\trp}_j)$
are stabilizable and detectable, respectively. In~\cite[Proposition
3.3]{Sim16a}, a general condition on $A$ and $B$ is given which ensures
stabilizability of $(\tA_j$, $\tB_j)$ in the case $M=I$.  Unfortunately, in
practice it is difficult to check whether this condition holds; in the majority
of our experiments the projected CARE could be solved for a stabilizing
solution.  We observed, however, the occurrence of nearly imaginary eigenvalues
of the associated Hamiltonian matrix $
  \begin{sbm}
    \tA_j              & -\tB_j\tB_j^{\trp} \\
    -\tC_j^{\trp}\tC_j & -\tA_j^{\trp}
  \end{sbm}
$,
and these can cause numerical difficulties and accuracy losses in the methods
employed for solving the small-scale CAREs~\eqref{proj_care}. In such a
situation it is reasonable to improve the quality of $Y_j$ by defect correction
strategies, see e.g.~\cite{MehT88}.  In our experiments we, therefore, check the
quality of $Y_j$ using the norm of
$R_j:=\tA_j^\trp Y_j+Y_j\tA_j-Y_j\tB_j\tB_j^\trp Y_j+\tC_j^{\trp}\tC_j$ and, if
necessary, try to improve $Y_j$ by running at most two steps of a Newton scheme
for the CARE defined by $\tA_j$, $\tB_j\tB_j^\trp$, $R_j$.

The cost of solving a dense matrix equation is cubic in the dimension. Since the sizes of the small CAREs arising in each iteration of EKSM and RKSM are equal
to the dimension of the projection subspace, the small-scale CARE solution can become expensive in the later iterations. A basic strategy to reduce these costs
is
to avoid solving the projected matrix equation in each iteration step, and do it
only every $\jgal$ steps, where usually we take $\jgal=5$.
In the Galerkin accelerated NK-ADI iteration, a small-scale GCARE has to be
solved only at the end of each outer iteration step.
 \subsection{Computing or Estimating the GCARE Residual Norm}
\change{%
In Table~\ref{tab:residual} we summarize the major steps carried out
in each method for computing or estimating the GCARE residual norm
$\|\cR(\tX)\|$. For brevity we only list the parts that involve
$n$-dimensional vectors (matrix vector products with the sparse
coefficients $A,A^\trp,M,M^\trp$, inner products of $n$-vectors, thin
QR factorizations and norms of $n\times x$ block vectors) as those are
the most computationally intensive steps in the residual norm
calculation.
Since the residual norm computation is often tightly intertwined with
other parts of the methods (shift generation, solving projected
CAREs), the overview only provides computational steps that are exclusively used
for that distinct purpose.}

\change{The first observation is that in (G)EKSM, RKSM, RADI, and iNK-ADI+LS,
the computational effort remains constant during the iteration, and is
dominated by other operations.
That is the reason why the measured computation time for this task is usually negligible.
Following the approach in~\cite[Proposition~6.1]{DruSZ14}, the
residual norm computation in TRKSM requires the construction of an
additional array of size $n\times \bs$, $\mathbf{s}=p+\sum_{j=1}^\jit
p_j$, leading also to increasing costs.
In ILRSI and NK-ADI+GP, where $\|\cR(\tX)\|$ is estimated via a
Lanczos process applied to $\cR$, the effort naturally increases with the
iteration because the number of columns in the low-rank factors
increases. Note that in NK-ADI+GP, $\|\cR(\tX)\|$ only needs to be
estimated after each Newton step.
}
\begin{table}[t]
  \centering
  \caption{\change{Major computational steps involving $n$-vectors required for GCARE residual norm estimation in all methods.}}
  {\renewcommand{\arraystretch}{1.3}%
  \begin{tabularx}{\linewidth}{|l|X|}
\hline
    Method&computations\\
    \hline
    (G)E- \& RKSM&$2p$ matrix-vector products with $M$ (or factors thereof) and thin QR of $n\times 2p$ matrix\\
    TRKSM&$p_j$ matrix-vector products with $M$ (or factors thereof),  additional orthogonalization scheme to expand thin QR of size $n\times \sum_{i=1}^{j-1}p_{i-1}+p$ by $p_j$ vectors \\
    ILRSI&per Lanczos step on $\cR$: matrix-vector products with $A,A^\trp,M,M^\trp$, $(2\jit+1)p+m$ inner products\\
        RADI&norm of $n\times p$ matrix\\
        NK-ADI+GP&per Lanczos step on $\cR$: matrix-vector products with $A,A^\trp,M,M^\trp$, $(2\jit_{\text{inner}}+1)(p+m)$ inner products\\
    iNK-ADI+LS&norm of $n\times (p+2m)$ matrix\\
    \hline
  \end{tabularx}\label{tab:residual}}
\end{table}

\subsection{Memory Consumption}\label{sec:memory}
Next, we are interested in the memory consumption of the compared methods for,
both, scenarios \textbf{S1} and \textbf{S2}. For this purpose, we consider the
maximum number of $n$-dimensional vectors stored at any given time during a
single run of the algorithm.  The memory consumption of other quantities,
e.g.~the projected matrices, is neglected. We also neglect the memory required
for solving the linear systems, since this occurs in all methods, and is likely
to lead to similar amounts.

We only give some basic estimates on the required number of $n$-dimensional
vectors. These numbers might vary slightly depending on the actual
implementation of the algorithms. Assume $M\neq I$, that only real shifts are
used, and let $\jit$ denote the total number of iterations required to satisfy the termination criteria in a particular method.

Starting with the projection methods (T)RKSM and (G)EKSM, we recall
that they do not gain an advantage in scenario \textbf{S2}. The largest number
of $n$-dimensional vectors in both scenarios is used for storing the $\jit p$
basis vectors contained in the matrix $Q_{\jit}$. For RKSM, there are also $p$
auxiliary vectors arising when solving the linear systems, and the block matrix
$G_\jit$ containing $p$ vectors, which occurs during the residual computation as
described in Section~\ref{ssec:projectionmethods}. Generating the approximate
feedback $K_\jit$ requires an additional $m$ vectors. A similar count can be
done for block and global EKSM;\@ we only have to replace $\jit$ with $2\jit$
since each step of these methods expands the basis twice. For TRKSM, storing
the basis requires $\mathbf{s}=p+\sum_{j=1}^\jit p_j$ vectors, where
$p_j$ denotes the number of direction vectors in the $j$th iteration. In TRKSM~\cite[Proposition~6.1]{DruSZ14}, the computation of the shifts,
the tangential directions, and the residual norm requires additional
$\mathbf{s}$ vectors.

The amount of $n$-dimensional vectors in ILRSI can be easily read off from
Algorithm~\ref{alg:ilrsi} and equals $(\jit+1)p$, along with additional $m$
vectors for the feedback approximation. For the RADI iteration
(Algorithm~\ref{alg:radi}) a closer look at the scenarios \textbf{S1} and
\textbf{S2} is appropriate. Starting with \textbf{S2}, we require a total of
$2p+2m$ vectors: there are $m$ vectors needed for solving the linear systems via
the SMW formula, $2p$ vectors for storing $R_\jit$ and $V_+$, and finally $m$
vectors for storing $K_j$. In scenario \textbf{S1}, the low-rank
factor $Z_\jit$ naturally adds $\jit p$ vectors. Since we use the residual-Hamiltonian shifts in ILRSI/RADI,
additional $3\ell p$ basis vectors are required in for storing $U_{\ell}$,
$A^\trp U_{\ell}$, $M^\trp U_{\ell}$.

Finally, consider the low-rank NK-ADI iteration. Let $\jit$ now denote the largest number of inner (ADI) iteration steps reached in any of the outer
Newton-Kleinman iterations. In both scenarios, the Galerkin projected variant (NK-ADI+GP) requires $\jit (p+m)$ vectors to store the basis, $2p+m$ vectors
to
use the SMW formula, $2(p+m)$ vectors
to store $R_\jit$ and $V_+$, and $m$ vectors to store $K_\jit$. In scenario~\textbf{S1}, the inexact version with line-search (iNK-ADI+LS) needs
additional $p+2m$
vectors to carry out the line-search technique. On the other hand, in scenario~\textbf{S2}, iNK-ADI+LS does not require all $\jit (p+m)$ vectors to store
the low-rank factor $Z_{\jit}$, but only $2\ell(p+m)$ vectors in order to generate the projection based shifts.

These basic estimates for the amount of $n$-dimensional vectors are summarized in Table~\ref{tab:nvec}, where superscripts are added to $\jit$ to indicate the
particular method in question and to highlight that all the methods might need different numbers of steps.
\begin{table}[t]
\centering
\caption{Number of required $n$-dimensional vectors for all algorithms in scenarios \textbf{S1} and \textbf{S2}, $\jit^{\textit{xyz}}$ indicates
the executed number of steps of method \textit{xyz}.}
   \setlength{\tabcolsep}{0.3em}
 	{\renewcommand{\arraystretch}{1.3}%
\begin{tabularx}{\linewidth}{|X|l|l|}
\hline
Method&\textbf{S1}&\textbf{S2}\\
\hline
{(G)}EKSM&$2(\jitgek+1)p$&$2(\jitgek+1)p+m$\\
RKSM&$(\jitrk+2)p$&$(\jitrk+2)p+m$\\
TRKSM&$2 \mathbf{s} +p$&$2 \mathbf{s} +p+m$\\
ILRSI&$(\jitir+1+3\ell)p$&$(\jitir+1+3\ell)p+m$\\
RADI&$(\jitradi+2+3\ell) p+2m$&$(2+3\ell)p+2m$\\
NK-ADI+GP&$(\jitnkgp+2+3\ell) (p+m)+2m$&$(\jitnkgp+2+3\ell) (p+m)+2m$\\
iNK-ADI+LS&$(\jitnkls+3+3\ell) (p+m)+2m$&$(3\ell+3) (p+m)+2m$\\
\hline
\end{tabularx}\label{tab:nvec}}
 \end{table}

The memory consumption in scenario~\textbf{S1} is for all methods dominated by the number of required iteration steps $\jit$ times the number $p$.
RKSM appears
to be efficient in this situations as it only requires a small amount of auxiliary $n$-vectors compared to the other methods.
Regarding scenario~\textbf{S2}, we see that RADI and iNK-ADI+LS have an advantage, since their storage requirements are
independent on the number of taken iteration steps $\jit$. The number $\ell$ is a free parameter that affects the shift generation, and it is typically is
chosen very
small, e.g.,  $1\leq \ell\leq 6$.

\section{Numerical experiments}\label{sec:numex}
For all methods described in this paper, in the following comparison we use
\matlab{} implementations\footnote{The codes and examples data for
the experiments in this section are available at DOI:10.5281/zenodo.3662519.}. These incorporate the
latest advances of each method, in particular, dynamic shift generation techniques.

State-of-the-art implementations of the low-rank Newton-Kleinman iterations can
be found in the \MMESS{} package~\cite{SaaKB16-mmess-1.0.1}, which is a collection
of currently mostly low-rank ADI based algorithms for large, sparse GCALEs, GCAREs,
differential Lyapunov and Riccati equations, as well as routines for related
balancing based model order reduction approaches.
Our implementations of
the projection based methods of Section~\ref{ssec:projectionmethods} were
originally inspired by the source codes
available\footnote{\url{http://www.dm.unibo.it/~simoncin/software.html}} from
Valeria Simoncini's homepage. The major differences in our \matlab{} codes include
a
defect correction strategy (cf. Section~\ref{sec:smallCARE}) for the compressed CAREs
and the calculation of the true residual norm for GCAREs as discussed in
Section~\ref{ssec:projectionmethods}.  The linear solves with $A$ in {(G)}EKSM
are carried out using sparse LU factors of \(A\),
when $A\neq A^{\trp}$, and using sparse Cholesky factors of $-A$ when
$A=A^{\trp}\prec 0$.  As described before, these factorizations are computed only
once and reused thereafter.
In all methods except {(G)}EKSM, the arising shifted linear sparse systems were
solved by the \matlab{} backslash ``$\backslash$'' operator and by using the SWM
formula when needed.  If not stated otherwise, the small, dense GCAREs in
{(G)}EKSM and {(T)}RKSM are solved at every $\jgal=5$ iteration steps. The subspace
for the shift generation in RADI, IRLSI, as well as NK-ADI  is taken from the previous $\ell=2$
iteration steps, i.e., $U_{\ell}\in\R^{n\times 2p}$.

The methods were compared on a number of standard examples available from
benchmark collections and the literature on large-scale Riccati equations.
Our test environment consists of \matlab{} 8.0.0.783
(R2012b) running on an \intel{}\xeon{} CPU X5650 (2.67GHz) with 48GB RAM.\@ All
methods were stopped when
\begin{align*}
  \rho:=\|\cR(\tX)\|_2/\|C^{\trp}C\|_2\leq \tau_{\cR}\text{ with
  }\tau_{\cR}=10^{-8},
\end{align*}
where $\tilde{X}$ is the approximate solution produced by the algorithm. Let
$\cL_k$ indicate the continuous-time Lyapunov
operator~\eqref{newtonNK_ARE_GCALE} at iteration $k$ of the NK-ADI method. Then
the stopping criterion for the inner Lyapunov ADI iteration is
$\|\cL_k(\tX)\|<\tau_{\cL}\|W^{(0)}\|^2$ with $\tau_{\cL}=\tau_{\cR}/10$ in the
Galerkin accelerated version. In the inexact version it is chosen adaptively via
$\tau_{\cL_k}=\eta\rho_{k-1}^2$, where $0<\eta<1$ and $\rho_{k-1}$ is the GCARE
residual norm from the previous Newton step;
see~\cite{BenHSetal16,Wei16} for details.
\subsection{First test series}
For a selection of the examples, Table~\ref{tab:results1} summarizes the results
obtained with the methods introduced in Section~\ref{sec:methods}.  Next to some
basic setup information for the particular example, the table gives the column
dimension of the built up low-rank factor after termination, the numerical rank of the
approximate GCARE solution using the machine precision as tolerance, the final
scaled GCARE residual $\rho_j$, the computation times in seconds for subtasks
such as solving linear systems ($t_{\text{LS}}$), small-scale GCAREs
($t_{\text{care}}$), computing shift parameters ($t_{\text{shift}}$), the total
computation time ($t_{\text{total}}$), the times $t_{\text{rest}}$ of all remaining minor computations (e.g., estimating the residual norm and further minor
auxiliary routines), as well as the largest number of stored
$n$-vectors regarding scenarios \textbf{S1} (mem$_Z$) and \textbf{S2} (mem$_K$).
The orthogonalization costs for building the orthonormal bases in the projection
based methods is included in small-scale solutions times $t_{\text{care}}$.  If
an algorithm failed for a particular example, a brief comment about the reason
is provided. In the following, whenever some defining matrices of a test example
are generated randomly, this is to be understood as matrices with normally
distributed entries for which the random number generator of \matlab{} was at the start of each experiment
initialized by \texttt{randn{('state', 0)}}.

\begin{table}[!h]
  \centering
  \caption{Testing results. The columns show: dimension of the constructed
    subspace; rank of final low-rank solution; final relative residual
    norm~$\rho$; time (in seconds) spent solving linear systems
    $t_{\text{LS}}$, small-scale CAREs $t_{\text{care}}$, and computing shifts
    $t_{\text{shift}}$; total computation time $t_{\text{total}}$; time $t_{\text{rest}}$ for remaining minor subtasks; peak memory
    consumption in terms of the number of stored $n$ vectors for scenarios
    \textbf{S1} (mem$_Z$) and \textbf{S2} (mem$_K$).}
  \footnotesize
  \setlength{\tabcolsep}{0.6em}
  \begin{tabularx}{\textwidth}{|l|X|r|r|r|r|r|r|r|r|r|r|}
    \hline
    Ex.
    &method&dim.&rank&final res.&$t_{\text{LS}}$&$t_{\text{care}}$&$t_{\text{shift}}$&$t_{\text{rest}}$&$t_{\text{total}}$&mem$_Z$&mem$_K$\\
    \hline\hline
    \multirow{8}{*}{
    \begin{minipage}{0.125\linewidth}
      Chip,\\
      $n=20\,082$,\\
      $m=1$,\\
      $p=5$
    \end{minipage}
    }
    &EKSM&250&89&1.9e-10&6.2&2.6&--&0.3&9.1&260&261\\
&RKSM&130&91&1.6e-10&34.6&1.4&2.2&0.5&38.7&140&141\\
&TRKSM&114&90&5.5e-10&37.8&2.0&0.9&0.3&41.0&238&239\\
&GEKSM&350&96&7.4e-09&8.0&4.2&--&0.1&12.3&360&361\\
&ILRSI&180&92&7.2e-10&48.4&--&1.7&3.1&53.3&185&186\\
&RADI&185&91&2.4e-09&51.5&--&0.7&0.2&52.5&227&42\\
&NK-ADI+GP&180&92&2.0e-11&45.0&0.6&0.2&1.6&47.4&192&192\\
&iNK-ADI+LS&180&90&1.7e-09&153.5&--&0.4&1.1&155.0&201&33\\
    \hline
    \multirow{7}{*}{
    \begin{minipage}{0.125\linewidth}
      Filter3D,\\
      $n=106\,437$,\\
      $m=1$,\\
      $p=5$
    \end{minipage}
    }
    &EKSM&1\,300&218&2.7e-09&422.4&447.4&--&23.4&893.2&1\,310&1\,311\\
    &RKSM,GEKSM&\multicolumn{10}{c|}{no convergence in maximum iterations}\\
    &TRKSM&216&210&5.1e-10&138.6&79.6&5.7&5.4&229.3&442&443\\
    &ILRSI&220&210&6.8e-09&107.9&--&17.1&33.3&158.3&255&256\\
    &RADI&215&207&6.0e-09&97.8&--&6.2&2.2&106.2&257&42\\
    &NK-ADI+GP&\multicolumn{10}{c|}{failure at solving projected CARE}\\
    &iNK-ADI+LS&348&241&1.7e-09&667.7&--&10.8&27.7&706.2&381&51\\
    \hline
    \multirow{8}{*}{
    \begin{minipage}{0.125\linewidth}
      Rail,\\
      $n=317\,377$,
      $m=7$,\\$p=6$
    \end{minipage}
    }
    &EKSM&840&222&3.3e-09&199.3&368.6&--&24.3&592.2&852&859\\
&RKSM&210&201&3.0e-10&118.3&83.9&11.8&18.3&232.3&222&229\\
&TRKSM&186&186&7.6e-10&113.7&102.9&3.2&13.1&232.9&384&391\\
&GEKSM&1800&224&1.6e-09&383.0&2406.1&--&3.7&2792.8&1812&1817\\
&ILRSI&204&183&5.7e-09&88.6&--&60.7&30.9&180.2&260&267\\
&RADI&204&183&5.8e-09&101.1&--&27.5&9.5&138.1&264&62\\
&NK-ADI+GP&324&229&6.2e-14&137.3&42.0&7.4&13.6&200.4&350&350\\
&iNK-ADI+LS&624&219&2.4e-09&197.0&--&8.3&52.3&257.6&684&86\\
    \hline
    \multirow{8}{*}{
    \begin{minipage}{0.125\linewidth}
      Toeplitz,\\
      $n=100\,000$,
      $m=5$,\\$p=20$
    \end{minipage}
    }
    &EKSM&600&378&1.9e-10&3.5&28.0&--&0.8&32.3&640&650\\
    &RKSM&320&320&8.0e-11&1.3&13.4&0.8&4.8&20.3&365&365\\
    &TRKSM&260&260&6.3e-09&1.0&12.4&2.5&3.7&19.6&305&305\\
    &GEKSM&600&378&1.9e-10&3.3&20.9&--&0.5&24.6&640&650\\
    &ILRSI&340&340&4.1e-09&0.8&--&12.7&15.8&29.3&400&405\\ 
    &RADI&280&280&4.6e-09&1.0&--&2.9&1.0&4.9&450&170\\ 
    &NK-ADI+GP&300&300&2.0e-10&0.6&8.0&0.6&1.8&11.0&350&350\\
    &iNK-ADI+LS&350&297&4.7e-09&2.8&--&1.4&6.8&11.0&440&215\\
    \hline
    \multirow{6}{*}{
    \begin{minipage}{0.125\linewidth}
      Lung,\\
      $n=109\,460$,\\
      $m=10$,\\
      $p=10$
    \end{minipage}
    }
&EKSM&100&100&6.6e-09&1.1&1.6&--&0.3&3.0&110&120\\
&RKSM&120&120&8.4e-10&8.1&3.9&3.0&2.3&17.3&140&150\\
&TRKSM&\multicolumn{10}{c|}{failure at solving projected CARE}\\
&GEKSM&100&93&8.1e-09&1.1&1.4&--&0.1&2.6&110&120\\
&ILRSI&250&211&9.2e-09&18.2&--&1.4&50.5&70.1&260&270\\
&RADI&80&80&9.7e-09&7.8&--&1.4&0.4&9.6&180&100\\
&NK-ADI+GP&\multicolumn{10}{c|}{failure at solving projected CARE}\\
&iNK-ADI+LS&440&220&9.7e-09&163.6&--&4.7&14.9&183.2&1\,370&190\\
    \hline
        \multirow{6}{*}{%
    \begin{minipage}{0.125\linewidth}
      Stokes,\\
      $n=67\,199$,\\
      $m=p=5$
    \end{minipage}
    }
    &EKSM&250&104&6.8e-09&8.9&18.2&--&3.1&30.2&260&265\\
    &RKSM&125&125&7.2e-09&15.4&10.5&2.0&4.0&31.9&135&140\\
    &ILRSI&175&106&8.6e-10&19.9&--&2.7&68.4&91.0&180&185\\
    &RADI&110&103&1.5e-09&20.3&--&0.9&0.9&22.1&160&50\\
    &NK-ADI+GP&140&107&4.1e-09&16.1&0.3&0.4&5.8&22.6&160&160\\
    &iNK-ADI+LS&320&100&2.1e-09&69.2&--&2.0&4.2&75.4&365&65\\
    \hline
  \end{tabularx}\label{tab:results1}
\end{table}

We proceed by giving details on the examples and discussing the obtained results
for each example.

\begin{example}[Chip]\label{ex:CHIP}
  The \textit{Chip} model~\cite{morMooRGetal04} represents a finite element discretization of a
  cooling process of a microchip. This benchmark example from the
  model order reduction wiki (MorWiki)\footnote{Available at
\url{http://modelreduction.org/index.php/Convection}}~\cite{morWiki} (also part of the Oberwolfach benchmark collection) provides $A\neq A^{\trp}$, $M\succ 0$
diagonal with
  $n=20\,082$ and $m=1$, $p=5$.

  For this first example, all methods generate approximate GCARE solutions of
  more or less the same rank.  EKSM and its global version, GEKSM, achieve the
  smallest total computation times $t_{\text{total}}$. The reason is clearly
  that both methods only require a single factorization of $A$ which for this
  example could be computed very cheaply and, hence, EKSM and GEKSM also spend a
  comparatively very small amount of time ($t_{\text{LS}}$) in solving linear
  systems. However, the generated subspace dimensions are significantly larger
  compared to all other methods, with GEKSM even further surpassing EKSM.\@
  Consequently, a substantial amount (about 50\%) of the total computation time
  is spent for solving the reduced CARE ($t_{\text{care}}$) and the memory
  consumptions is higher than for the other methods.  This will be a typical
  observation also in almost all further experiments. Furthermore, there is also
  a large difference between the subspace dimension and the actual rank of the
  solution, indicating that the extended Krylov subspace approaches may lead to
  unnecessarily large projection subspaces for the purpose of solving GCAREs,
  resulting in a waste of computational effort in several cases as it will be
  evident in later examples.  Using a rational Krylov subspace appears to be a
  better option, as RKSM and its tangential version, TRKSM, require much smaller
  subspaces, which are actually close to the rank of the solutions. Thus,
  smaller storage requirements and small-scale solution times $t_{\text{care}}$
  can be expected. However, for the Chip example this is compensated by
  larger times $t_{\text{LS}}$ since the varying linear systems prevent a
  prefactorization of the system matrix. The setup of the subspaces in TRKSM
  leads to a further reduction of the subspace dimension compared to the
  standard block version in RKSM.\@

  Both ILRSI and RADI perform similarly, though slightly worse than {(T)}RKSM,
  since they need more iterations. This leads to higher subspace dimensions,
  linear system and total computation times.  Since $m=1$, the more complicated
  linear systems in RADI do not cause a significant difference in
  $t_{\text{LS}}$ compared to ILRSI.\@ RADI achieves the smallest memory
  requirement among all methods in scenario~\textbf{S2}, which is in line with
  the discussion in Section~\ref{sec:memory}.

  For the Chip example, the performance of NK-ADI+GP is overall similar to RADI
  and ILRSI.\@ The Galerkin acceleration led to a termination after the first
  Newton step. A very poor performance is exhibited by the iNK-ADI+LS iteration
  which requires many more linear systems solves and as a consequence, much
  higher times $t_{\text{LS}}$, $t_{\text{total}}$ compared to the other
  algorithms. It can only compete regarding the memory consumption in
  scenario~\textbf{S2}.

  It is also important to emphasize that the shift generation times
  $t_{\text{shift}}$ are always a very small fraction of the total times
  $t_{\text{total}}$ in all methods relying on shift parameters. 
\end{example}

\begin{example}[Filter3D]\label{ex:filter}
  Another MorWiki~\cite{morWiki} example represents a finite element model of a tunable optical
  filter\footnote{Available at
\url{http://modelreduction.org/index.php/Tunable_Optical_Filter}}~\cite[Chapter 15]{morBenMS05} with $A\prec 0$, $M\succ 0$, $n=106\,437$, and $m=1$, $p=5$.

  Here, EKSM requires far more iterations and therefore, a much larger
  subspace than the other methods.  This results in a large amount of
  time spent in solving the
  projected CAREs, so that any savings from the easier linear systems solves in
  EKSM are counterbalanced. RKSM and GEKSM are not able to satisfy the
  termination criterion in a maximum allowed number of 200 iterations.  In view
  of the subspace dimension, TRKSM, ILRSI, and RADI perform similarly. Because
  TRKSM has to solve more linear systems with a different coefficient matrix in
  each step, and a projected CARE every $\jgal$ steps (see
  Section~\ref{sec:memory}), its overall computation time is significantly
  larger than for ILRSI and RADI. The smallest memory consumption in
  scenario~\textbf{S2} is again obtained by the RADI iteration.  The higher
  shift generation times $t_{\text{shift}}$ in ILRSI compared to RADI are mainly caused
  by the employed residual Hamiltonian approach requiring to project the
  residual matrix $\cR(X_k)$ onto a low-dimensional subspace. Here, in contrast
  to RADI, there is no \change{effective} low-rank structure of $\cR(X_k)$ known for ILRSI. Even
  though $U_{\ell}^{\trp}\cR(X_k)U_{\ell}$ can be computed without explicitly
  forming $\cR(X_k)$ and by exploiting its symmetry, this construction is,
  nevertheless, noticeably more expensive than in RADI.

  The projected CARE after the first Newton step of the NK-ADI+GP iteration
  could not be solved by the employed \texttt{care} routine, and the method
  broke down.  The iNK-ADI+LS iteration was able to solve this problem but, as
  in the Chip example, required larger times $t_{\text{LS}}$, $t_{\text{total}}$
  caused by the much larger number of encountered linear systems compared to the
  other algorithms. The dimension of the generated subspace is also
  higher.
\end{example}
\begin{example}[Rail]\label{ex:rail}
  The steel profile cooling models\footnote{Available at
\url{http://modelreduction.org/index.php/Steel_Profile}} are also part of the MorWiki.  They represent spatial
  finite element discretizations of a two-dimensional heat transfer problem arising in the
  cooling of steel rail profiles~\cite[Chapter 19]{morBenMS05}. Different grid sizes and
  discretization levels result in five versions of the example, each having
  different dimension: \(n=1\,357\),
  \(5\,177\),
  \(20\,209\),
  \(79\,841\),
  \(317\,377\).\footnote{The largest version was obtained by a \fenics~based
  reimplementation available at
    \url{http://gitlab.mpi-magdeburg.mpg.de/models/fenicsrail/} of the
  discretization.}
  In all versions, $A\prec 0$, $M\succ 0$, and the provided matrices $B$, $C$
  have parameters $m=7$ and $p=6$.  Table~\ref{tab:results1} shows the results
  for the largest version only, while the other versions are examined later.

  EKSM, GEKSM, and the iNK-ADI+LS iteration generate substantially larger
  subspaces than all other algorithms.  For 
  EKSM and GEKSM this leads to a very large effort for solving the projected
  Riccati equation such that again, any savings gained by the simpler linear
  systems are completely lost.  Apparently, the global method performs far worse
  than the block method: GEKSM has the largest subspace dimensions, and the
  largest times $t_{\text{LS}}$, $t_{\text{care}}$, $t_{\text{total}}$. These
  results emphasize again that the extended Krylov subspaces are often not an
  adequate choice for solving large GCAREs and that, moreover, no computational
  advantages should be expected by the global over the block approach.  Because
  the iNK-ADI+LS iteration does not rely on solving projected GCAREs, its
  overall computation times are significantly smaller compared to {(G)}EKSM, but
  clearly larger than for {(T)}RKSM, ILRSI, and RADI. The latter
  four algorithms appear to be the winners for this example, having much smaller
  subspace dimensions that are, especially for {(T)}RKSM, very close to the
  actual solution rank. ILRSI and RADI need the smallest
  computation times, overall, since no projected GCARE has to be solved. The
  NK-ADI+GP iteration is able to compete for this example, with total
  computation times in between those of {(T)}RKSM and ILRSI, RADI, but with a
  significantly smaller final residual
  norm. This is due to the fact that NK-ADI+GP generates a roughly 1.5
  times larger subspace compared to {(T)}RKSM and ILRSI, RADI, such
  that the Galerkin projection results in a much smaller residual norm than actually desired.
\end{example}
\begin{example}[Toeplitz]\label{ex:toep}
  This is an artificial example defined by the Toeplitz matrix
  \begin{align*}
    A=\begin{sbm}
      2.8&1&1&1&0&\ldots&\\
      -1&2.8&1&1&1&0&\ddots\\
      0&-1&2.8&1&1&1&\ddots\\
      \vdots&\ddots&\ddots&\ddots&\ddots&\ddots&\ddots\\
      &&&&0&-1&2.8
     \end{sbm}, \quad M=I_n,
  \end{align*}
  and random $B$, $C$ with $\|B\|=1$. We set $n=100\,000$, $m=5$, and
  $p=20$. The idea for this setting comes
  from~\cite[Example~7.3]{LinSim15},~\cite[Example~5.5]{Sim16}, where slightly
  different entries and $m=p=1$ are used. We use this example because the
  arising linear systems are extremely cheap to solve even for large $n$, which
  helps to emphasize the required work in other stages of the algorithms.

  The data given in Table~\ref{tab:results1} confirm this because the timings
  $t_{\text{LS}}$ represent only a small fraction of $t_{\text{total}}$ for all
  methods. The timings $t_{\text{care}}$, and in some cases even the shift
  generation times $t_{\text{shift}}$, are higher than $t_{\text{LS}}$.  In
  fact, more than 50--75\% of time is spent in solving the small-scale CAREs.
  The methods {(G)}EKSM again build up larger subspaces than the
  other methods. All of the other approaches end up with significantly smaller
  subspace dimensions, which are all close to
  the actual solution rank. The approximate solution of the smallest rank is
  obtained by TRKSM which, having to solve projected CAREs, in the end does not
  achieve the smallest computation times $t_{\text{total}}$.  Because the RADI
  and iNK-ADI+LS methods do not have to solve such projected CAREs, they achieve
  the smallest total times $t_{\text{total}}$ with a substantial margin.  On par
  with iNK-ADI+LS regarding $t_{\text{total}}$ is the NK-ADI+GP method, which
  only requires one small CARE solve after the first Newton step. The cost for
  this is, however, a substantial portion of the overall cost, as reflected by
  $t_{\text{care}}$. Although ILRSI does not work with projected CAREs either, it
  ends up with the highest total time $t_{\text{total}}$ for two reasons: first,
  the shift generation is more costly (see observations made for Filter3D), and
  second, the costs for the residual norm estimation via a Lanczos process in
  this example are comparatively more expensive (around 50\% of the total
  computation time).
  Note that without the scaling $\|B\|=1$, ILRSI converges much slower because
  it encounters problems in generating good shift parameters.

  For this example, we also tried SDA~\cite{LiCLetal13}, but it was not able to
  compute a low-rank approximate solution of the desired accuracy.  We
  terminated SDA after 5 iteration steps which already took 247.8 seconds and
  led to a relative CARE residual norm
  $\|\cR(\tX)\|_2/\|C^{\trp}C\|_2\approx 0.998$. Interestingly, the built in
  residual norm estimation of SDA delivered a value $2.53\cdot10^{-13}$ and
  incorrectly indicated convergence of the method. The produced low-rank
  solution factors after these 5 iteration steps had reached the maximal allowed
  column dimension 600.\qed{}
\end{example}
\begin{example}[Lung]\label{num:ex:lung}
  The example \textsc{lung2} from the SuiteSparse Matrix
  Collection~\cite{DavH11} (formerly UF Sparse Matrix Collection) models
  temperature and water vapor transport in the human lung.  It provides matrices
  with leading dimension $n=109\,460$, where $A$ is nonsymmetric, $M=I$, and
  $B$, $C$ are generated as random matrices with $m=p=10$.

  For this example, {(G)}EKSM wins regarding the total computation times
  $t_{\text{total}}$ since, in contrast to the prior examples, they do not
  produce larger subspace dimension than the other methods. The extended Krylov method
  appears to be a viable choice for this example, and the cheaper linear solves,
  similar to the Chip example, pay off here.  RKSM has a larger time
  $t_{\text{LS}}$, because it, again, has to put more effort into solving the
  varying linear systems. TRKSM and NK-ADI+GP break down because the employed
  routine \texttt{care} for solving the projected CARE fails at some
  point. RADI achieves the smallest subspace dimension, but needs
  more total time $t_{\text{total}}$, due to the varying linear systems. For
  this example, ILRSI encounters problems in the dynamic shift generation and
  requires many more iterations compared to RADI.\@ We expect these issues to be
  the result of complex data in the basis used for projecting the residual
  Hamiltonian matrix.  Similar to prior observations, iNK-ADI+LS takes by far
  the last place regarding several measures, e.g., subspace dimension,
  $t_{\text{LS}}$, and $t_{\text{total}}$.\qed{}
\end{example}

\begin{example}[Stokes]\label{num:ex:stokes}
A spatial discretization of a two-dimensional Stokes equation from~\cite{morSch07} is included as test case of \MMESS{} consisting of matrices
 \begin{align*}
  A=\begin{sbm}
         A_1&G\\G^{\trp}&0
         \end{sbm},\quad M:=\begin{sbm}
         M_1&0\\0&0
         \end{sbm}\in\Rnn,\quad B=\begin{sbm}
         B_1\\0
         \end{sbm}\in\Rnm,\quad C=\begin{sbm}
         C_1&0
         \end{sbm}\in\Rpn,
 \end{align*}
 with $n=n_1+n_2$, $n_1>n_2$, $A_1,M_1\in\R^{n_1\times n_1}$ sparse, $G\in\R^{n_1\times n_2}$, $B_1\in\R^{n_1\times m}$, $C_1\in\R^{p\times n_1}$.
 This results in a descriptor system of index-2, such that using a Riccati based feedback stabilization requires some additional steps.
 For the sake of brevity, we only give the absolutely necessary information regarding the handling of those structured descriptor systems, especially with
respect to the
numerical steps in the GCARE algorithms. With rank$(G)=n_2$ we can associate the projector
\begin{align*}
 \Pi:=I_{n_1}-G(G^{\trp}M_1G)^{-1}G^{\trp}M_1
\end{align*}
to the differential algebraic system. The purpose of $\Pi$ is to ensure that algebraic constraints of the descriptor system are satisfied, i.e., the
calculations happen in the correct hidden
manifold.
Following~\cite{morHeiSS08, BaeB12, BaeBSetal15, Wei16, BenHSetal19}, the GCARE to be solved is defined by the large, dense matrices
$\cA=\Theta_1^{\trp}A_1\Theta_1$, $\cM=\Theta_1^{\trp}M_1\Theta_1$, $\cB=\Theta_1^{\trp}B_1$,  $\cC=C_1\Theta_1$, where
$\Theta_1,\Theta_2\in\R^{n_1\times (n_1-n_2)}$ are factors of $\Pi$: $\Pi=\Theta_1\Theta_2^{\trp}\in\R^{n_1\times n_1}$.

It has been shown in~\cite{BaeBSetal15, Wei16, BenHSetal19} that the main numerical subtasks in GCARE methods based on the low-rank NK-ADI framework
can be
implemented without the explicit projection $\Pi$ or its factors $\Theta_1, \Theta_2$, such that working with the original matrices $A,M,B$ is sufficient
and numerically desirable since $A,M$
are sparse. Only the transformed right hand side factor $C_{\Pi}:=(\Pi C^{\trp}_1)^{\trp}$ is required to start the iteration. It is straightforward to carry
these ideas over to the direct iterations (Section~\ref{sec:nonprojectiveiterations}) and the projection methods (Section~\ref{ssec:projectionmethods}).
An application of the projection $\Pi$ requires solving a
symmetric indefinite linear system in saddle point form defined by $\hM:=\begin{sbm}
         M_1&F\\F^{\trp}&0
         \end{sbm}$.
In the projection based methods, it is possible that the algorithms drift off the hidden manifold~\cite{StyS12}, especially after the orthogonalization
scheme for expanding the orthonormal basis. Hence, it is wise to apply $\Pi$ also to the outcome of the orthogonalization routine (repeated modified
Gram-Schmidt in this exposition). Moreover, estimating the GCARE residual norm via a Lanczos process, as it is done in ILRSI, NK-ADI+GP, also requires
applications of $\Pi$ in each Lanczos step~\cite[Chapter 4.3]{Wei16}.
Consequently, the occurrence of $\Pi$-application adds an additional
source of numerical cost to each method not present in the examples before.
How much extra effort this introduces depends on how often $\Pi$ needs to be applied: once in RADI, iNK+ADI+LS, $1+\jit$ times in EKSM, (T)RKSM,
$j_{\text{Lan}}$ times after each Galerkin projection in NK-ADI+GP, and $1+\jit j_{\text{Lan}}$ times in ILRSI, where $j_{\text{Lan}}$ indicates the number
of executed Lanczos steps for estimating the GCARE residual norm.  The computation of the residual Hamiltonian shifts in ILRSI also
relies explicitly on $\cR$, s.t. further $\Pi$-applications are needed (cf.
discussion in Example~\ref{ex:filter}).

Here, we use a discretization with $150$ grid points in each spatial dimension resulting in $n=67\,199$, $n_1=44\,700$, $n_2=22\,201$, $m=p=5$. These
sizes still allow to employ sparse-direct techniques to handle the linear systems defined by $\hM$ as well as all further occurring
linear systems in the GCARE methods. In particular, we use a sparse $LDL^{\trp}$ factorization of $\hM$ which is computed once before the
iterations.

The results are given at the bottom of Table~\ref{tab:results1}. The winner for this example in terms of computation time, subspace dimension, and memory
requirements is RADI, closely followed by the projection methods  EKSM, RKSM and NK-ADI+GP, although EKSM again produced larger subspaces. Similar as before,
iNK+ADI+LS required much more LR-ADI steps leading to larger times $t_{\text{LS}}$, $t_{\text{total}}$, and larger subspace dimensions.
The last place takes
ILRSI because of the significantly more expensive residual norm estimations leading to large $t_{\text{rest}}$ timings. This is similar as in the
\textit{lung} example, but here, the Lanczos process additionally requires $\Pi$-applications.
\end{example}

This concludes this first series of examples. We point out that the amount of
work for adaptively generating shifts is, with minor exceptions, only a small
fraction of the overall computational effort due to the advances in this
topic~\cite{DruS11,DruSZ14,BenKS14,Sim16a,Kue16,BenBKetal15} in the recent
years. If cleverly implemented, precomputed shifts could also be generated
efficiently, but in our experience, the performance of the low-rank GCARE (and GCALE)
methods typically lags behind compared to dynamic shift selection.

Although not the topic of this study, we expect a similar conclusion
regarding the related methods for large, sparse Lyapunov and Sylvester
equations.  We did not put timings for the residual norm computation in
Table~\ref{tab:results1} because these were only a tiny fraction of the overall
computation time, mostly because of low-rank or otherwise exploitable structures
of the GCARE residual matrix which made this task significantly better
manageable than with earlier approaches. The clear exception is here ILRSI
which still relies on a Lanczos process for this task and, hence, could benefit
from low-rank expressions \change{with fixed-sized factors} of the GCARE residual. This would further improve the performance of the employed residual Hamiltonian shift generation in ILRSI.
\subsection{Scaling of the performance with respect to the problem dimensions}
Now we examine the behavior of the algorithms when the leading dimension $n$ is
increased. For this, we use all five versions, described above, for the Rail
example, which correspond to different coarseness of finite element meshes.

\begin{figure}[!t]
  \caption{Scaling with respect to the problem dimension of \textit{Rail}
    examples: plotted are the generated subspace dimension (left), the ranks of the approximate solutions (middle), and the total
    computation time (right) versus the leading dimension
    $n$.}\label{fig:rail_scale_n}
  \centering
  \begin{tikzpicture}

\begin{axis}[%
width=0.365\linewidth,
name=raildim,
xmode=log,
xmin=1000,
xmax= 320000,
xminorticks=false,
ymin=100,
ymax=1800,
cycle list name=res,
xlabel={$n$},
ylabel={coldim($Z$)},
legend entries={EKSM,RKSM,TRKSM,GEKSM,ILRSI,RADI,NK-ADI+GP,iNK-ADI+LS},
legend style={draw=black,fill=white,legend pos=north west,below right,
at={(0.26,1.5)},
legend columns=4,
font=\small,
legend cell
align=left,rounded
corners=3pt},
every axis x label/.style={at={(ticklabel cs:0.5)},font=\scriptsize,rotate=0,anchor=center,yshift=-.2em},
every axis y label/.style={at={(ticklabel cs:0.5)},font=\scriptsize,rotate=90,anchor=center,yshift=.2em},
every  y tick label/.style={font=\scriptsize},
every  x tick label/.style={font=\scriptsize},
]
\addplot 
  table[row sep=crcr]{1357	300\\
5177	420\\
20209	540\\
79841	660\\
317377	840\\
};

\addplot 
  table[row sep=crcr]{1357	150\\
5177	150\\
20209	186\\
79841	180\\
317377	210\\
};

\addplot 
  table[row sep=crcr]{1357	127\\
5177	157\\
20209	159\\
79841	160\\
317377	186\\
};

\addplot 
  table[row sep=crcr]{1357	660\\
5177	960\\
20209	1260\\
79841	1620\\
317377	1800\\
};

\addplot 
  table[row sep=crcr]{1357	156\\
5177	168\\
20209	180\\
79841	192\\
317377	204\\
};

\addplot 
  table[row sep=crcr]{1357	156\\
5177	168\\
20209	180\\
79841	210\\
317377	204\\
};

\addplot 
  table[row sep=crcr]{1357	210\\
5177	300\\
20209	288\\
79841	420\\
317377	324\\
};

\addplot 
  table[row sep=crcr]{1357	442\\
5177	507\\
20209	676\\
79841	676\\
317377	624\\
};
\end{axis}

\begin{axis}[%
width=0.365\linewidth,
name=railrank,
at={(raildim.north east)},
xshift=3.25em,
anchor=north west,
xmode=log,
xmin=0,
xmax= 320000,
xminorticks=false,
ymin=120,
ymax=240,
cycle list name=res,
xlabel={$n$},
ylabel={rank($\tX$)},
every axis x label/.style={at={(ticklabel cs:0.5)},font=\scriptsize,rotate=0,anchor=center,yshift=-.2em},
every axis y label/.style={at={(ticklabel cs:0.5)},font=\scriptsize,rotate=90,anchor=center,yshift=.2em},
every  y tick label/.style={font=\scriptsize},
every  x tick label/.style={font=\scriptsize},
]
\addplot 
  table[row sep=crcr]{1357	140\\
5177	161\\
20209	185\\
79841	207\\
317377	222\\
};

\addplot 
  table[row sep=crcr]{1357	138\\
5177	148\\
20209	181\\
79841	177\\
317377	201\\
};

\addplot 
  table[row sep=crcr]{1357	127\\
5177	156\\
20209	159\\
79841	160\\
317377	186\\
};

\addplot 
  table[row sep=crcr]{1357	147\\
5177	169\\
20209	193\\
79841	214\\
317377	224\\
};

\addplot 
  table[row sep=crcr]{1357	127\\
5177	141\\
20209	158\\
79841	170\\
317377	183\\
};

\addplot 
  table[row sep=crcr]{1357	124\\
5177	140\\
20209	157\\
79841	179\\
317377	183\\
};

\addplot 
  table[row sep=crcr]{1357	140\\
5177	165\\
20209	189\\
79841	214\\
317377	229\\
};

\addplot 
  table[row sep=crcr]{1357	134\\
5177	159\\
20209	178\\
79841	199\\
317377	219\\
};
\end{axis}
%

\begin{axis}[%
width=0.365\linewidth,
at={(railrank.north east)},
xshift=3.25em,
anchor=north west,
xmode=log,
xmin=1000,
xmax=320000,
xminorticks=false,
ymode=log,
ymin=0,
ymax=5000,
yminorticks=false,
ytick={10,100,1000},
cycle list name=res,
xlabel={$n$},
ylabel={time},
every axis x label/.style={at={(ticklabel cs:0.5)},font=\scriptsize,rotate=0,anchor=center,yshift=-.2em},
every axis y label/.style={at={(ticklabel cs:0.5)},font=\scriptsize,rotate=90,anchor=center,yshift=.2em},
every  y tick label/.style={font=\scriptsize},
every  x tick label/.style={font=\scriptsize},
]

\addplot 
  table[row sep=crcr]{1357	1.671842\\
5177	5.126833\\
20209	18.125566\\
79841	75.395434\\
317377	592.222277\\
};

\addplot 
  table[row sep=crcr]{1357	1.83152\\
5177	2.647356\\
20209	7.553109\\
79841	31.340926\\
317377	232.281323\\
};

\addplot 
  table[row sep=crcr]{1357	1.411588\\
5177	3.115928\\
20209	7.611653\\
79841	31.95364\\
317377	232.897317\\
};

\addplot 
  table[row sep=crcr]{1357	19.093735\\
5177	84.331566\\
20209	269.113771\\
79841	906.724853\\
317377	2792.753713\\
};

\addplot 
  table[row sep=crcr]{1357	0.557835295169946\\
5177	1.325649\\
20209	5.142777\\
79841	30.219349\\
317377	180.246787\\
};

\addplot 
  table[row sep=crcr]{1357	0.386693609254499\\
5177	0.934955048076923\\
20209	3.758528\\
79841	19.6063550328947\\
317377	138.063078\\
};

\addplot 
  table[row sep=crcr]{1357	1.36518479035793\\
5177	3.95192301947888\\
20209	7.81234208420783\\
79841	43.375312\\
317377	200.396930000349\\
};

\addplot 
  table[row sep=crcr]{1357	1.750561\\
5177	4.728969\\
20209	32.42917\\
79841	149.678486003434\\
317377	257.592788\\
};

\end{axis}
\end{tikzpicture}%
\end{figure}
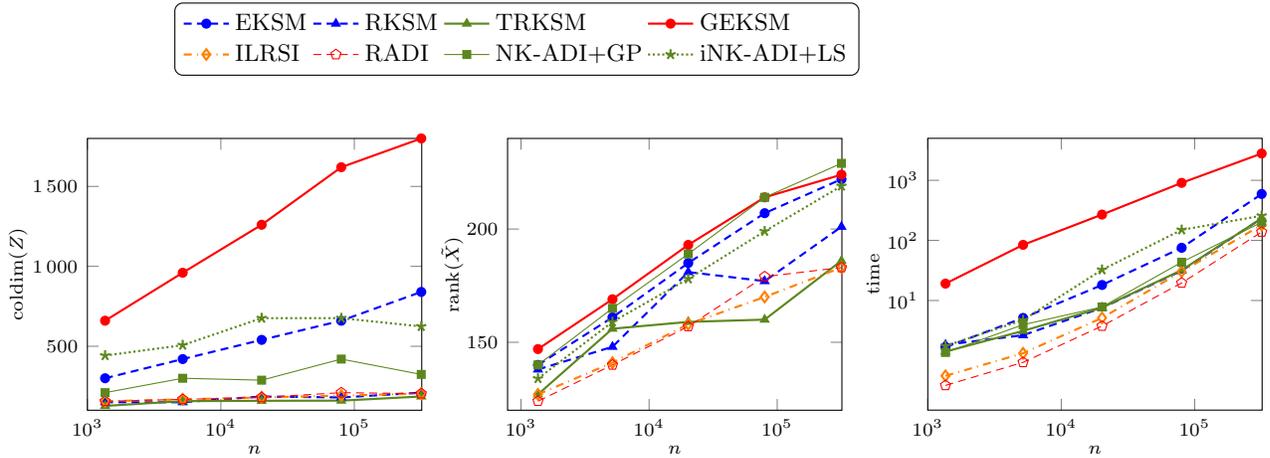

The subspace dimensions, the solution ranks, and total computation times are plotted against the
different values of $n$ in Figure~\ref{fig:rail_scale_n}.  For this particular
example, the increasing dimensions do seem to only lead to moderately larger ranks
of the approximate GCARE solutions required to satisfy the stopping criterion, see middle plot of~Figure~\ref{fig:rail_scale_n}.
Hence,
if the algorithms manage to produce a subspace of dimension not far from the
solution rank, the numerical effort should essentially only increase with $n$ because the
linear systems are increasingly expensive to solve.  In the left plot of
Figure~\ref{fig:rail_scale_n} we see that indeed, the majority of methods ends up
with subspace dimensions that remain approximately unchanged for increasing
$n$. The striking exceptions are {(G)}EKSM whose subspace dimensions are again
much larger than those of the other methods and do also clearly increase with
$n$. Similar to prior observation, GEKSM performs worse than EKSM.\@ For the
largest example, GEKSM requires a subspace dimension more than three times as
large as for the other methods.
The NK-ADI methods produce somewhat larger, and slightly increasing subspace dimensions
compared to {(T)}RKSM, ILRSI, and RADI.

The total computation times illustrated in the right plot of
Figure~\ref{fig:rail_scale_n} indicate that, as expected, increasing system
dimensions $n$ lead to increasing computation times because of the larger
linear systems. The trend in which $t_{\text{total}}$ increases is similar for
{(T)}RKSM, ILRSI, RADI, NK-ADI+GP, especially for dimensions $n>10^4$. The
timings for {(G)}EKSM and the iNK-ADI+LS iteration are somewhat larger than for the other algorithms. For the iNK-ADI+LS
iteration, the difference in $t_{\text{total}}$ to the other methods except GEKSM appears to decrease for increasing $n$.

\smallskip

The purpose of the next examples is to study how the change in the dimensions $m$
and $p$ affects the performance of the algorithms.  We keep the matrices $A$ and
$M$ fixed, but alter the number of columns in $B$ and the number of rows in
$C$. For most methods, this will result in a different number of columns in the
right hand side when solving linear systems, see
Table~\ref{tab:linsys}. Furthermore, for the projection based methods,
the size of $p$ also dictates the
growth of the subspace dimension and, thus, significantly influences the cost for solving the Galerkin systems.
Another way of manipulating the influence of linear systems during the algorithm
runs is to provide an initial guess $X_0$, e.g., when solving an unstable Riccati equation. Changing $m$ changes in this case the number of right hand sides in
the linear system if the SMW formula is employed (cf. Table~\ref{tab:linsys}); we study this effect as well.

\begin{example}[CUBE-FD]\label{ex:fdm3d}
  Consider a centered finite difference discretization of the differential
  equation
  \begin{align*}
    \partial_t f(\xi,t) = \Delta f(\xi,t) - 10 \xi_1 \partial_{\xi_1} f(\xi,t) - 1\,000 \xi_2 \partial_{\xi_2} f(\xi,t) - 10 \partial_{\xi_3} f(\xi,t)+
b(\xi)u(t),
  \end{align*}
  on a unit cube with \(\xi=(\xi_{1},\xi_{2},\xi_{3})\). Using Dirichlet boundary conditions and
  $n_0=32$ nodes in each spatial direction yields a nonsymmetric matrix $A$ of
  order $n=n_0^3$ and $M=I_n$. The matrices $B\in\Rnm$, $C\in \Rpn$ are
  generated randomly.  For $n_0=22$, $C=B^{\trp}$, and $m=10$, this is
  exactly~\cite[Example 5.2]{Sim07}.
   \end{example}
   We will also use the last example to evaluate the influence of a nonzero initial
  feedback $K_0$. However, for several unstable variations of this examples (e.g. using shifted matrices $A+\psi I$, $\psi>\vert \max
\Real{\lambda(A)}\vert$), we
  encountered severe problems with the initial feedback generation explained in
  Section~\ref{ssec:unstable}. Either the required eigen- or Schur vectors could
  not be computed in a stable manner, or the small matrix equations (Bernoulli or
  Lyapunov) could not be solved satisfactorily.
  To avoid these problems, we modify Example~\ref{ex:fdm3d} as follows.
\begin{example}[CUBE-FD-unstable]\label{ex:fdm3dunstable}
  Consider
  $A_u=\diag{A,A_+}$, $B_u^{\trp}=[B^{\trp},B_+^{\trp}]$, where $A,~B$ are as in Example~\ref{ex:fdm3d},
  $A_+\in\R^{u\times u}$ is an
  artificially generated matrix with eigenvalues in $\C_+$, and
  $B_+\in\R^{u\times m}$, $C\in\R^{p\times u+n}$ are given and randomly generated. We used \(u=5\) in our experiments.
  The required anti-stable invariant subspace for the above $A_u$
  is trivially spanned by $Q_u={[0,I_u]}^{\trp}\in\R^{(n+u)\times u}$ and it
  holds that $Q_u^{\trp}B_u=B_+$, so we set $K_0={[0,B^{\trp}_+]}^{\trp}$.  To
  generate the anti-stable part $A_+$, we enforce the solution of the Bernoulli
  equation $A_+S_u+S_{u}A^{\trp}_+-S_u(Q_u^{\trp}B)(B^{\trp}Q_u)S_u=0$ to be
  $S_u=I_u$.  Therefore, $A_++A_+^{\trp}=B_+B_+^{\trp}$, and we simply take
  $A_+=\half B_+B_+^{\trp}$.  We again point out that this construction of an
  unstable system is done entirely for demonstrating the effects of an initial
  feedback $K_0$.  From a practical point of view, the way
  $A_u,B_u,C_u$ are built makes this setting viable for partial stabilization
  approaches~\cite{HeM94,Var95c,Ben11}.  Moreover, as we mentioned earlier, projection
  based methods might, under certain conditions~\cite[Proposition 3.3]{Sim16a},
  not need an initial guess to converge.
  \end{example}

  \begin{table}[!t]
    \centering
    \caption{Testing results for \textit{CUBE-FD} and \textit{CUBE-FD-unstable} with $n=32\,768$,
      $\tau_{\cR}=10^{-8}$  and different values of $m,p$.}

    \footnotesize
    \setlength{\tabcolsep}{0.7em}

    \begin{tabularx}{\linewidth}{|l|X|r|r|r|r|r|r|r|r|r|r|}
      \hline
      Setup
      &method&dim.&rank&final res.&$t_{\text{LS}}$&$t_{\text{care}}$&$t_{\text{shift}}$&$t_{\text{rest}}$&$t_{\text{total}}$&mem$_Z$&mem$_K$\\
      \hline\hline
      \multirow{8}{*}{\begin{minipage}{0.05\linewidth}
               \rotatebox{90}{\mbox{$m=10$, $p=10$}}
                      \end{minipage}}%
      &EKSM&1\,200&568&6.9e-09&27.7&142.7&--&1.9&172.4&1\,220&1\,230\\
      &RKSM&460&460&7.4e-09&114.3&17.3&6.3&2.6&140.4&480&490\\
      &TRKSM&470&470&5.4e-09&114.4&18.4&26.7&1.9&161.5&960&970\\
      &GEKSM&1\,300&558&5.7e-09&29.8&157.8&--&0.3&187.9&1\,320&1\,330\\
      &ILRSI&520&490&9.6e-09&142.5&--&18.6&19.5&180.6&570&580\\
      &RADI&520&471&8.9e-09&186.8&--&5.8&1.1&193.6&620&100\\
      &NK-ADI+GP&550&550&1.5e-10&134.7&9.3&0.7&1.1&145.8&590&590\\
      &iNK-ADI+LS&1\,120&546&9.4e-09&477.6&--&3.3&3.3&484.9&1\,210&190\\
      \hline
      \multirow{8}{*}{\rotatebox{90}{\mbox{$m=10$, $p=20$}}}%
      &EKSM&2\,400&1\,099&4.4e-09&56.8&944.3&--&4.1&1\,005.3&2\,440&2\,450\\
      &RKSM&1\,120&1\,100&3.2e-10&200.4&114.7&24.9&6.7&346.6&1\,160&1\,170\\
      &TRKSM&1\,060&1\,060&4.6e-10&185.0&92.7&131.4&5.3&414.5&2\,160&2\,170\\
      &GEKSM&2\,800&1\,096&3.0e-09&67.8&1\,553.4&--&0.6&1\,621.8&2\,840&2\,850\\
      &ILRSI&1\,040&978&6.0e-09&199.4&--&52.2&39.8&291.3&1\,140&1\,150\\
      &RADI&1\,040&942&7.0e-09&246.6&--&18.7&2.3&267.7&1\,220&180\\
      &NK-ADI+GP&1\,140&1\,120&4.0e-11&204.9&50.5&2.0&2.3&259.6&1\,200&1\,200\\
      &iNK-ADI+LS&1\,740&1\,081&6.6e-09&616.6&--&6.6&6.6&629.4&1\,860&270\\
      \hline
      \multirow{8}{*}{\rotatebox{90}{\mbox{$m=20$, $p=10$}}}%
      &EKSM&1\,200&568&8.2e-09&25.6&148.0&--&1.9&175.5&1\,220&1\,240\\
      &RKSM&460&460&7.7e-09&126.9&18.6&6.7&2.6&154.8&480&500\\
      &TRKSM&480&480&6.0e-09&129.5&21.5&29.3&2.0&182.3&980&1\,000\\
      &GEKSM&1\,400&569&2.8e-09&29.1&210.6&--&0.2&239.9&1\,420&1\,240\\
      &ILRSI&530&509&7.3e-09&173.6&--&20.5&22.4&216.5&580&600\\
      &RADI&540&499&2.1e-09&269.4&--&8.2&1.5&279.1&660&120\\
      &NK-ADI+GP&570&570&6.2e-11&156.2&10.1&0.7&1.2&168.2&630&630\\
      &iNK-ADI+LS&1\,680&534&9.7e-09&692.1&--&5.8&16.0&703.9&1\,830&300\\
      \hline
      \multirow{6}{*}{\rotatebox{90}{\parbox{7em}{\mbox{$m=20$, $p=10$}\\unstable}}}%
      &EKSM&1\,200&571&8.0e-09&29.4&139.2&--&0.4&169.0&1\,220&1\,240\\
      &RKSM&570&545&2.5e-09&168.5&28.1&11.0&3.3&210.9&590&610\\
      &TRKSM&464&464&5.5e-09&136.7&21.9&50.6&2.1&211.3&948&968\\
      &ILRSI&600&504&5.5e-09&187.7&--&4.9&108.8&301.4&670&680\\
      &RADI&600&490&5.5e-09&280.5&--&4.9&1.8&287.2&720&120\\
      &NK-ADI+GP&800&572&9.2e-11&267.3&37.5&1.9&4.8&311.5&860&860\\
      &iNK-ADI+LS&2\,190&585&9.0e-09&1\,467.1&--&11.8&15.5&1\,494.4&2\,550&300\\
      \hline
    \end{tabularx}\label{tab:cube_scale_pm}
  \end{table}

  \begin{figure}[!t]
    \caption{\textit{CUBE-FD(-unstable)}: Results with respect to different values of $m$
      and $p$: generated subspace dimension (left), the ranks of the approximate solutions (middle), and the total
    computation time (right) vs $p$ or $m$.
      Top row of plots refers to $m=10$, $p=1, 10, 20, 30, 40$, the middle row refer to $p=10$,
      $m=1, 10, 20, 30, 40$, and in the bottom row the results for the unstable case \textit{CUBE-FD-unstable} with
      $p=10$ and varying $m=5, 10, 20, 30, 40$ are
      plotted.}\label{fig:cubefd_scale}
    \centering
    \begin{tikzpicture}

\begin{axis}[%
width=0.365\linewidth,
name=cubefd_dim_p,
xmin=1,
xmax=40,
ymin=0,
ymax=6000,
cycle list name=res,
xlabel={$p$},
ylabel={coldim($Z$)},
legend entries={EKSM,RKSM,TRKSM,GEKSM,ILRSI,RADI,NK-ADI+GP,iNK-ADI+LS},
legend style={draw=black,fill=white,legend pos=north west,below right,
at={(0.26,1.5)},
legend columns=4,
font=\small,
legend cell
align=left,rounded
corners=3pt},
every axis x label/.style={at={(ticklabel cs:0.5)},font=\scriptsize,rotate=0,anchor=center,yshift=-.2em},
every axis y label/.style={at={(ticklabel cs:0.5)},font=\scriptsize,rotate=90,anchor=center,yshift=.2em},
every  y tick label/.style={font=\scriptsize},
every  x tick label/.style={font=\scriptsize},
]
\addplot
  table[row sep=crcr]{1	130\\
10	1200\\
20	2400\\
30	3300\\
40	4400\\
};

\addplot
  table[row sep=crcr]{1	46\\
10	460\\
20	1120\\
30	1530\\
40	2680\\
};

\addplot
  table[row sep=crcr]{1	46\\
10	470\\
20	1060\\
30	1410\\
40	2120\\
};

\addplot
  table[row sep=crcr]{1	120\\
10	1300\\
20	2800\\
30	4200\\
40	5600\\
};

\addplot
  table[row sep=crcr]{1	49\\
10	520\\
20	1040\\
30	1530\\
40	2040\\
};

\addplot
  table[row sep=crcr]{1	50\\
10	520\\
20	1040\\
30	1530\\
40	2040\\
};

\addplot
  table[row sep=crcr]{1	56\\
10	550\\
20	1140\\
30	1770\\
40	2280\\
};

\addplot
  table[row sep=crcr]{1	715\\
10	1120\\
20	1740\\
30	2280\\
40	2800\\
};

\end{axis}

\begin{axis}[%
width=0.365\linewidth,
name=cubefd_rank_p,
at={(cubefd_dim_p.north east)},
xshift=3.25em,
anchor=north west,
xmin=1,
xmax=40,
ymin=1,
ymax=2500,
cycle list name=res,
xlabel={$p$},
ylabel={rank($\tX$)},
every axis x label/.style={at={(ticklabel cs:0.5)},font=\scriptsize,rotate=0,anchor=center,yshift=-.2em},
every axis y label/.style={at={(ticklabel cs:0.5)},font=\scriptsize,rotate=90,anchor=center,yshift=.2em},
every  y tick label/.style={font=\scriptsize},
every  x tick label/.style={font=\scriptsize},
]
\addplot
  table[row sep=crcr]{1	62\\
10	568\\
20	1099\\
30	1543\\
40	2016\\
};

\addplot
  table[row sep=crcr]{1	46\\
10	460\\
20	1100\\
30	1518\\
40	2134\\
};

\addplot
  table[row sep=crcr]{1	0\\
10	470\\
20	1060\\
30	1410\\
40	2098\\
};

\addplot
  table[row sep=crcr]{1	60\\
10	558\\
20	1096\\
30	1598\\
40	2079\\
};

\addplot
  table[row sep=crcr]{1	49\\
10	490\\
20	978\\
30	1426\\
40	1857\\
};

\addplot
  table[row sep=crcr]{1	50\\
10	471\\
20	942\\
30	1419\\
40	1842\\
};

\addplot
  table[row sep=crcr]{1	56\\
10	550\\
20	1120\\
30	1642\\
40	2130\\
};

\addplot
  table[row sep=crcr]{1	62\\
10	546\\
20	1081\\
30	1558\\
40	2020\\
};

\end{axis}

\begin{axis}[%
width=0.365\linewidth,
name=cudefd_time_p,
at={(cubefd_rank_p.north east)},
xshift=3.25em,
anchor=north west,
xmin=1,
xmax=40,
ymin=0,
ymax=13500,
cycle list name=res,
xlabel={$p$},
ylabel={time},
ymode=log,
yminorticks=false,
every axis x label/.style={at={(ticklabel cs:0.5)},font=\scriptsize,rotate=0,anchor=center,yshift=-.2em},
every axis y label/.style={at={(ticklabel cs:0.5)},font=\scriptsize,rotate=90,anchor=center,yshift=.2em},
every  y tick label/.style={font=\scriptsize},
every  x tick label/.style={font=\scriptsize},
]
\addplot
  table[row sep=crcr]{1	8.883157\\
10	172.377463\\
20	1005.263716\\
30	2363.756825\\
40	5628.958706\\
};

\addplot
  table[row sep=crcr]{1	73.532385\\
10	140.408017\\
20	346.618754\\
30	517.408029\\
40	1917.822628\\
};

\addplot
  table[row sep=crcr]{1	0\\
10	161.487314\\
20	414.4899\\
30	644.987428\\
40	1661.728385\\
};

\addplot
  table[row sep=crcr]{1	7.613466\\
10	187.899509\\
20	1621.755733\\
30	5402.912363\\
40	13264.715996\\
};

\addplot
  table[row sep=crcr]{1	90.584298\\
10	180.612392003366\\
20	291.31353\\
30	411.42859\\
40	556.041756\\
};

\addplot
  table[row sep=crcr]{1	128.558308005383\\
10	193.614397\\
20	267.651071\\
30	335.370913000259\\
40	417.927268999356\\
};

\addplot
  table[row sep=crcr]{1	83.0489280008669\\
10	145.800274004979\\
20	259.576480002003\\
30	458.434375\\
40	713.469168002332\\
};

\addplot
  table[row sep=crcr]{1	408.159764\\
10	484.898113999765\\
20	629.447163\\
30	727.370429999409\\
40	847.555541000542\\
};

\end{axis}
\end{tikzpicture}%

%
%
%
\definecolor{mycolor1}{rgb}{1.00000,0.00000,1.00000}%
\begin{tikzpicture}

\begin{axis}[%
width=0.365\linewidth,
name=cubefd_dim_m,
at={(cubefd_dim_p.south west)},
xmin=1,
xmax=40,
ymin=0,
ymax=4000,
cycle list name=res,
xlabel={$m$},
ylabel={coldim($Z$)},
every axis x label/.style={at={(ticklabel cs:0.5)},font=\scriptsize,rotate=0,anchor=center,yshift=-.2em},
every axis y label/.style={at={(ticklabel cs:0.5)},font=\scriptsize,rotate=90,anchor=center,yshift=.2em},
every  y tick label/.style={font=\scriptsize},
every  x tick label/.style={font=\scriptsize},
]
\addplot
  table[row sep=crcr]{1	1200\\
10	1200\\
20	1200\\
30	1200\\
40	1200\\
};

\addplot
  table[row sep=crcr]{1	560\\
10	460\\
20	460\\
30	570\\
40	470\\
};

\addplot
  table[row sep=crcr]{1	480\\
10	470\\
20	480\\
30	480\\
40	480\\
};

\addplot
  table[row sep=crcr]{1	530\\
10	530\\
20	530\\
30	530\\
40	530\\
};

\addplot
  table[row sep=crcr]{1	540\\
10	540\\
20	540\\
30	540\\
40	540\\
};

\addplot
  table[row sep=crcr]{1	540\\
10	540\\
20	540\\
30	540\\
40	540\\
};

\addplot
  table[row sep=crcr]{1	570\\
10	570\\
20	570\\
30	570\\
40	570\\
};

\addplot
  table[row sep=crcr]{1	638\\
10	1200\\
20	1680\\
30	2520\\
40	3150\\
};

\end{axis}

\begin{axis}[%
width=0.365\linewidth,
name=cubefd_rank_m,
at={(cubefd_dim_m.north east)},
xshift=3.25em,
anchor=north west,
xmin=1,
xmax=40,
ymin=400,
ymax=600,
cycle list name=res,
xlabel={$m$},
ylabel={rank($\tX$)},
every axis x label/.style={at={(ticklabel cs:0.5)},font=\scriptsize,rotate=0,anchor=center,yshift=-.2em},
every axis y label/.style={at={(ticklabel cs:0.5)},font=\scriptsize,rotate=90,anchor=center,yshift=.2em},
every  y tick label/.style={font=\scriptsize},
every  x tick label/.style={font=\scriptsize},
]
\addplot
  table[row sep=crcr]{1	568\\
10	568\\
20	568\\
30	568\\
40	568\\
};

\addplot
  table[row sep=crcr]{1	560\\
10	460\\
20	460\\
30	570\\
40	470\\
};

\addplot
  table[row sep=crcr]{1	480\\
10	470\\
20	480\\
30	480\\
40	480\\
};

\addplot
  table[row sep=crcr]{1	569\\
10	569\\
20	569\\
30	569\\
40	569\\
};

\addplot
  table[row sep=crcr]{1	509\\
10	509\\
20	509\\
30	509\\
40	509\\
};

\addplot
  table[row sep=crcr]{1	499\\
10	499\\
20	499\\
30	499\\
40	505\\
};

\addplot
  table[row sep=crcr]{1	570\\
10	570\\
20	570\\
30	570\\
40	570\\
};

\addplot
  table[row sep=crcr]{1	562\\
10	570\\
20	534\\
30	586\\
40	585\\
};

\end{axis}
%
%
%
%

\begin{axis}[%
width=0.365\linewidth,
name=cubefd_time_m,
at={(cubefd_rank_m.north east)},
xshift=3.25em,
anchor=north west,
xmin=1,
xmax=40,
ymin=100,
ymax=1400,
ymode=log,
cycle list name=res,
xlabel={$m$},
ylabel={time},
every axis x label/.style={at={(ticklabel cs:0.5)},font=\scriptsize,rotate=0,anchor=center,yshift=-.2em},
every axis y label/.style={at={(ticklabel cs:0.5)},font=\scriptsize,rotate=90,anchor=center,yshift=.2em},
every  y tick label/.style={font=\scriptsize},
every  x tick label/.style={font=\scriptsize},
]
\addplot
  table[row sep=crcr]{1	171.139408\\
10	173.742786\\
20	175.506996\\
30	176.737327\\
40	173.241497\\
};

\addplot
  table[row sep=crcr]{1	194.833863\\
10	153.736581\\
20	154.793037\\
30	201.948465\\
40	154.942238\\
};

\addplot
  table[row sep=crcr]{1	174.524462\\
10	170.398898\\
20	182.325646\\
30	176.455232\\
40	176.783705\\
};

\addplot
  table[row sep=crcr]{1	247.168935\\
10	242.774016\\
20	239.921551\\
30	240.280448\\
40	244.836689\\
};

\addplot
  table[row sep=crcr]{1	187.270239978215\\
10	202.063634\\
20	216.478819\\
30	228.546234003352\\
40	242.363747\\
};

\addplot
  table[row sep=crcr]{1	162.845270973376\\
10	216.457394\\
20	279.089934000236\\
30	335.849054\\
40	397.912304\\
};

\addplot
  table[row sep=crcr]{1	151.148869032926\\
10	171.793109\\
20	168.179373\\
30	170.156940000353\\
40	168.932409\\
};

\addplot
  table[row sep=crcr]{1	326.211224002379\\
10	547.724820000329\\
20	703.904306999128\\
30	974.68358100129\\
40	1165.738036\\
};

\end{axis}
\end{tikzpicture}%

%
%
\begin{tikzpicture}

\begin{axis}[%
width=0.365\linewidth,
name=cubefd_unst_dim_m,
at={(cubefd_dim_m.south west)},
xmin=5,
xmax=40,
ymin=0,
ymax=1500,
cycle list name=res,
xlabel={$m$},
ylabel={coldim($Z$)},
every axis x label/.style={at={(ticklabel cs:0.5)},font=\scriptsize,rotate=0,anchor=center,yshift=-.2em},
every axis y label/.style={at={(ticklabel cs:0.5)},font=\scriptsize,rotate=90,anchor=center,yshift=.2em},
every  y tick label/.style={font=\scriptsize},
every  x tick label/.style={font=\scriptsize},
]
\addplot 
  table[row sep=crcr]{5	1200\\
10	1200\\
20	1200\\
30	1200\\
40	1200\\
};

\addplot 
  table[row sep=crcr]{5	570\\
10	570\\
20	570\\
30	570\\
40	570\\
};

\addplot 
  table[row sep=crcr]{5	455\\
10	475\\
20	464\\
30	464\\
40	523\\
};
\pgfplotsset{cycle list shift=1}
\addplot
  table[row sep=crcr]{5	630\\
10	600\\
20	600\\
30	590\\
40	580\\
};

\addplot 
  table[row sep=crcr]{5	630\\
10	600\\
20	600\\
30	590\\
40	580\\
};

\addplot 
  table[row sep=crcr]{5	920\\
10	770\\
20	800\\
30	910\\
40	830\\
};

\addplot 
  table[row sep=crcr]{5	1020\\
10	1360\\
20	2190\\
30	2920\\
40	3950\\
};

\end{axis}

\begin{axis}[%
width=0.365\linewidth,
name=cubefd_unst_rank_m,
at={(cubefd_unst_dim_m.north east)},
xshift=3.25em,
anchor=north west,
xmin=5,
xmax=40,
ymin=440,
ymax=620,
cycle list name=res,
xlabel={$m$},
ylabel={rank($\tX$)},
every axis x label/.style={at={(ticklabel cs:0.5)},font=\scriptsize,rotate=0,anchor=center,yshift=-.2em},
every axis y label/.style={at={(ticklabel cs:0.5)},font=\scriptsize,rotate=90,anchor=center,yshift=.2em},
every  y tick label/.style={font=\scriptsize},
every  x tick label/.style={font=\scriptsize},
]
\addplot
  table[row sep=crcr]{5	550\\
10	571\\
20	571\\
30	571\\
40	571\\
};

\addplot
  table[row sep=crcr]{5	534\\
10	545\\
20	545\\
30	545\\
40	546\\
};

\addplot
  table[row sep=crcr]{5	455\\
10	475\\
20	464\\
30	464\\
40	523\\
};
\pgfplotsset{cycle list shift=1}

\addplot
  table[row sep=crcr]{5	569\\
10	569\\
20	569\\
30	569\\
40	569\\
};

\addplot
  table[row sep=crcr]{5	505\\
10	514\\
20	504\\
30	505\\
40	505\\
};

\addplot
  table[row sep=crcr]{5	501\\
10	495\\
20	490\\
30	491\\
40	493\\
};

\addplot
  table[row sep=crcr]{5	566\\
10	571\\
20	572\\
30	615\\
40	574\\
};

\addplot
  table[row sep=crcr]{5	541\\
10	553\\
20	585\\
30	589\\
40	585\\
};

\end{axis}
%
%

\begin{axis}[%
width=0.365\linewidth,
name=cubefd_unst_time_m,
at={(cubefd_unst_rank_m.north east)},
xshift=3.25em,
anchor=north west,
xmin=5,
xmax=40,
ymin=100,
ymax=5e2,
cycle list name=res,
xlabel={$m$},
ylabel={time},
every axis x label/.style={at={(ticklabel cs:0.5)},font=\scriptsize,rotate=0,anchor=center,yshift=-.2em},
every axis y label/.style={at={(ticklabel cs:0.5)},font=\scriptsize,rotate=90,anchor=center,yshift=.2em},
every  y tick label/.style={font=\scriptsize},
every  x tick label/.style={font=\scriptsize},
]
\addplot 
  table[row sep=crcr]{5	165.548201\\
10	170.160433\\
20	169.041926\\
30	175.545303\\
40	171.866834\\
};

\addplot
  table[row sep=crcr]{5	192.661764\\
10	195.895804\\
20	210.909547\\
30	224.532116\\
40	238.181356\\
};

\addplot 
  table[row sep=crcr]{5	182.214424\\
10	202.127348\\
20	211.276616\\
30	224.391023\\
40	286.153832\\
};
\pgfplotsset{cycle list shift=1}

\addplot 
  table[row sep=crcr]{5	285.135896700968\\
10	279.331029419715\\
20	301.40915224812\\
30	309.629074968528\\
40	320.422600583748\\
};
\addplot 
  table[row sep=crcr]{5	195.342289\\
10	220.472162994825\\
20	287.169072998221\\
30	345.352422992651\\
40	399.773547991775\\
};
%
\addplot 
  table[row sep=crcr]{5	268.55391300162\\
10	247.871197\\
20	311.462117\\
30	437.321723\\
40	472.546064\\
};

\addplot 
  table[row sep=crcr]{5	1559.71204399997\\
10	1280.34135\\
20	1494.37148300069\\
30	2654.5082\\
40	4508.013781\\
};

\end{axis}
\end{tikzpicture}%
  \end{figure}
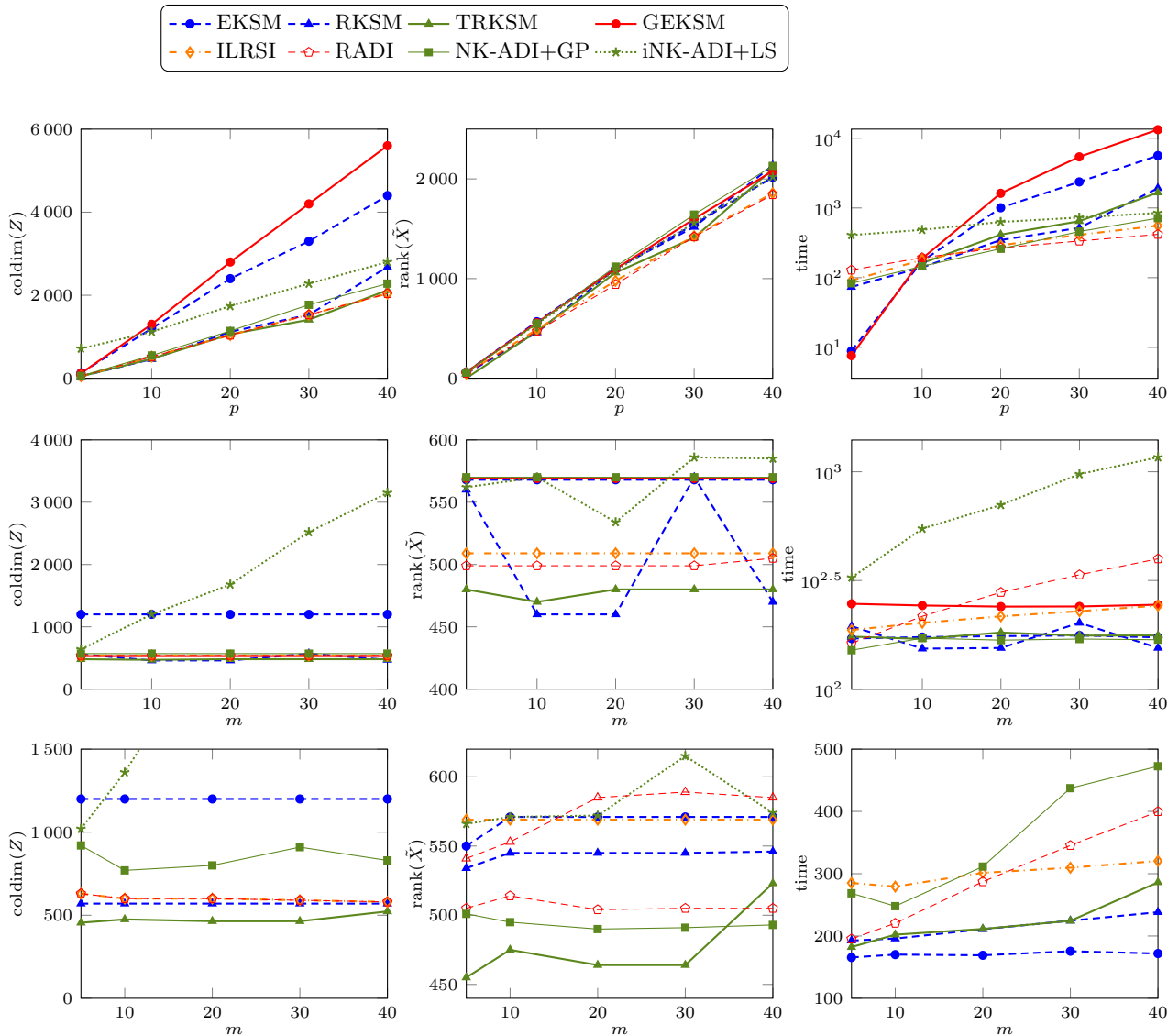
  Table~\ref{tab:cube_scale_pm} and Figure~\ref{fig:cubefd_scale} show the
  results.  We first comment on the top row in Figure~\ref{fig:cubefd_scale};
  the number $m$ of columns of $B$ is fixed to $10$, while the number $p$
  of rows in $C$ increases.
  The
  generated subspace dimensions as well as the obtained solution ranks increase almost linearly with increasing $p$
  for all methods.  This is expected, as each step in each method (except TRKSM)
  expands this subspace by a multiple of $p$
  vectors.
  The slope in top-left graph of Figure~\ref{fig:cubefd_scale} is the
  largest for {(G)}EKSM, since these two methods add $2p$
  vectors in each step, while the rest add only $p$
  vectors. The iNK-ADI+LS method appears to build subspaces larger than the
  rest, with the exception of {(G)}EKSM.\@ The linear increase in
  $t_\text{LS}$
  for all methods is obvious, and for the methods not based on projection, the
  total time of computation increases linearly as well. However, in the
  projection methods, the costs of solving larger projected GCAREs are getting
  increasingly expensive as $p$
  gets larger, resulting in weaker performance. This is, again, most obvious for
  {(G)}EKSM.\@

  Next, we study the middle row of Figure~\ref{fig:cubefd_scale}: now $p=10$,
  while $m$
  increases. Apparently, this increase has no effect on the
  generated subspace dimension at all, which is clear from the way the subspaces
  are expanded, and, thus, indicates that the iteration numbers also remain approximately constant. The exception is iNK-ADI+LS, where after the first Newton
step,
  $p+m$
  vectors are added in each inner iteration.  However, the obtained solution ranks appear to be more sensitive regarding different values of $m$, especially
for RKSM, iNK-ADI+LS. The effect on the total
  time is also different: the methods that use the SMW formula for the solution of
  linear systems (RADI, iNK-ADI+LS) see a linear increase in total time, which
  is a consequence of the increase in $t_\text{LS}$.
  The NK-ADI+GP is for this setting not affected by larger $B$,
  because in all cases the Galerkin acceleration led to satisfaction of the stopping criterion after the
  first Newton step. In consequence, each encountered linear system had only $p$
  right hand sides and the SMW formula was not necessary. The timings for ILRSI
  increase slightly since, for similar reasons as outlined above, the shift
  generation became more expensive with increasing $m$.

  Finally, we analyze the case of unstable Riccati equations for which we fix
  $p=10$ and vary $m=5,10,20,\ldots,40$.  The bottom row in
  Figure~\ref{fig:cubefd_scale} shows the results. For this setting, GEKSM
does not produce reasonable results (stagnation at large residual norms or problems while solving the projected CARE) and is hence
    omitted. ILRSI
  encounters problems in the shift generation routine, and thus we use it with
  the shifts generated by RADI.\@ The related timings $t_\text{shift}$ of ILRSI are copied from those of RADI. Despite the nonzero initial feedback, larger
values of $m$ do
  not appear to severely affect the final subspace dimensions for all methods.  The
  strong exception is the iNK-ADI+LS iteration, which again performs
  comparatively bad for this example, similar to the cases described above.
  Only for this method, the subspace dimensions are much larger than for the
  other methods, surpassing even EKSM, and increasing $m$ clearly seems to lead
  to larger dimensions, as well. We cut off the associated bottom left plot in
  Figure~\ref{fig:cubefd_scale} at $1\,500$ because otherwise the very large
  values of iNK-ADI+LS distort the whole plot.  Some variations in the
  solution ranks can also be seen for varying $m$, similar to the
  previous experiment with stable matrix pair $(A,M)$.
  The total computation times
  $t_{\text{total}}$ increase at various rates for different methods as $m$
  increases, which is due to the increase in the time $t_\text{LS}$ needed for
  solving linear systems, since all of the methods use the SMW formula now. As
  explained in Section~\ref{ssec:linsys}, EKSM has an advantage here and its
  computation times only increase marginally.  RADI and NK-ADI+GP exhibit the
  strongest increase of $t_{\text{total}}$. The same holds for iNK-ADI+LS whose
  curve is omitted in the bottom right plot in Figure~\ref{fig:cubefd_scale}.

  For $m=20$, the results are listed in the bottom section of
  Table~\ref{tab:cube_scale_pm}. We see that EKSM is the fastest, followed by
  {(T)}RKSM;\@ these two methods achieve the smallest subspace dimensions.
  It is also evident that iNK-ADI+LS cannot keep up with any of the other
  algorithms.  Note that, the residual norm estimation via Lanczos in ILRSI
  takes up a significant portion of the overall time, here.

  In respect to the memory consumptions for scenario~\textbf{S2}, RADI clearly requires the smallest number of $n$-vectors to be stored.
  The NK-ADI+GP achieves, similar to the \textit{rail} example before, smaller final residual norms which is again a result of the larger generated subspace
dimensions.
  \qed{}

\subsection{Comparison of theoretically equivalent methods}
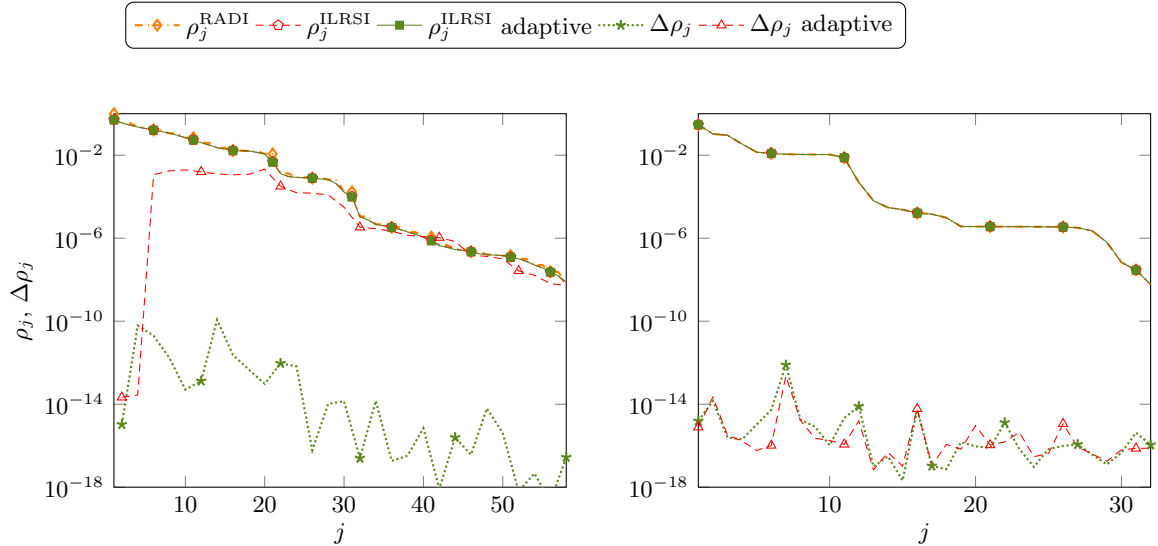
\begin{figure}[!t]
  \caption{Comparison of RADI and ILRSI with and without adaptive shifts for the
    \textit{CUBE-FD} (left) and \textit{rail} (right) examples: scaled residual
    norms $\rho_j$ and discrepancies $\Delta\rho_j$ against iteration number
    $j$.}\label{fig:ILRSIvsRADI}
  \centering
%
%
\begin{tikzpicture}

\begin{axis}[%
name = cubeRI,
width=0.46\linewidth,
xmin=1,
xmax=58,
ymode=log,
ymin=1e-18,
ymax=1,
yminorticks=false,
cycle list name=resN,
cycle list shift=4,
xlabel={$j$},
ylabel={$\rho_j$, $\Delta\rho_j$},
legend entries={$\rho^{\text{RADI}}_j$, $\rho^{\text{ILRSI}}_j$, $\rho_j^{\text{ILRSI}}$  adaptive ,
$\Delta\rho_j$, $\Delta\rho_j$ adaptive},
legend style={draw=black,fill=white,legend pos=north west,below right,
at={(0.025,1.3)},
legend columns=5,
font=\small,
legend cell
align=left,rounded
corners=3pt},
every axis x label/.style={at={(ticklabel cs:0.5)},font=\small,rotate=0,anchor=center,yshift=-.5em},
every axis y label/.style={at={(ticklabel cs:0.5)},font=\small,rotate=90,anchor=center,yshift=.5em},
every  y tick label/.style={font=\small},
every  x tick label/.style={font=\small},
]
\addplot 
  table[row sep=crcr]{1	1\\
2	0.355348873734358\\
3	0.355348873734358\\
4	0.226205640052024\\
5	0.226205640052024\\
6	0.163277296340095\\
7	0.163277296340095\\
8	0.113347461285675\\
9	0.113347461285675\\
10	0.0699649315643731\\
11	0.0699649315643731\\
12	0.03977099978916\\
13	0.03977099978916\\
14	0.0228731100343177\\
15	0.0228731100343177\\
16	0.0168098098417858\\
17	0.0168098098417858\\
18	0.0159976137785775\\
19	0.0159976137785775\\
20	0.0115696417919226\\
21	0.0115696417919226\\
22	0.00133038815677127\\
23	0.00133038815677127\\
24	0.000859241446404829\\
25	0.000859241446404829\\
26	0.000791816793310555\\
27	0.000791816793310555\\
28	0.000671269250951168\\
29	0.000671269250951168\\
30	0.000173818727985971\\
31	0.000173818727985971\\
32	1.11477570916188e-05\\
33	1.11477570916188e-05\\
34	4.64763079834962e-06\\
35	4.64763079834962e-06\\
36	3.34522056311623e-06\\
37	3.34522056311623e-06\\
38	1.97326808276775e-06\\
39	1.97326808276775e-06\\
40	1.14167293400719e-06\\
41	1.14167293400719e-06\\
42	4.70902853600987e-07\\
43	4.70902853600987e-07\\
44	2.9813911702983e-07\\
45	2.9813911702983e-07\\
46	2.26266192944107e-07\\
47	2.26266192944107e-07\\
48	1.62531607375472e-07\\
49	1.62531607375472e-07\\
50	1.45880457246685e-07\\
51	1.45880457246685e-07\\
52	1.00586433314122e-07\\
53	1.00586433314122e-07\\
54	5.10199576011494e-08\\
55	5.10199576011494e-08\\
56	2.37268661173036e-08\\
57	2.37268661173036e-08\\
58	6.26957780600882e-09\\
};

\addplot 
  table[row sep=crcr]{1	0.508883396431817\\
2	0.355348873734359\\
3	0.281712868407725\\
4	0.226205639990526\\
5	0.192839979660071\\
6	0.163277296320169\\
7	0.137120415357788\\
8	0.113347461283998\\
9	0.0909104028126601\\
10	0.0699649315643227\\
11	0.0544585564184362\\
12	0.0397709997890284\\
13	0.0311057730518157\\
14	0.0228731099173672\\
15	0.0197706186025708\\
16	0.0168098098394051\\
17	0.0163997009183692\\
18	0.0159976137781142\\
19	0.0138271757212637\\
20	0.0115696417918258\\
21	0.00480975357609115\\
22	0.00133038815583858\\
23	0.000939084867076153\\
24	0.000859241445735893\\
25	0.000823114072887129\\
26	0.000791816793310495\\
27	0.000730298223493544\\
28	0.000671269250940831\\
29	0.000434321750021684\\
30	0.000173818727972087\\
31	0.000101758362872758\\
32	1.11477570916438e-05\\
33	7.72408859676661e-06\\
34	4.64763078366888e-06\\
35	3.98529452367739e-06\\
36	3.34522056309774e-06\\
37	2.67250281387226e-06\\
38	1.97326808273488e-06\\
39	1.55558899953444e-06\\
40	1.14167293330775e-06\\
41	7.86373094340656e-07\\
42	4.7090285360015e-07\\
43	3.78162673330015e-07\\
44	2.98139116778897e-07\\
45	2.5904514881742e-07\\
46	2.26266192907258e-07\\
47	1.74457751200281e-07\\
48	1.62531601097389e-07\\
49	1.53817527472213e-07\\
50	1.45880456841218e-07\\
51	1.23134013376494e-07\\
52	1.00586433314757e-07\\
53	7.57470145940111e-08\\
54	5.10199576057238e-08\\
55	3.74504631108441e-08\\
56	2.37268661171049e-08\\
57	1.52935278867251e-08\\
58	6.26957777840288e-09\\
};

\addplot 
  table[row sep=crcr]{1	0.510002964997232\\
2	0.357062840282178\\
3	0.283079783435622\\
4	0.227396798919374\\
5	0.193818744377633\\
6	0.164077782465313\\
7	0.137731010109085\\
8	0.113789983227668\\
9	0.0912017154890883\\
10	0.0701194609350689\\
11	0.0545199714138148\\
12	0.0397443315798637\\
13	0.0310381252775375\\
14	0.0227662212347177\\
15	0.0196462534583833\\
16	0.0166689757299301\\
17	0.0162542317826915\\
18	0.0158475844006301\\
19	0.013690370348346\\
20	0.0114536764306894\\
21	0.00468995514247742\\
22	0.00129102510408361\\
23	0.000927744570894964\\
24	0.000852919861283645\\
25	0.000816626939491935\\
26	0.000785085277941438\\
27	0.00072229102839961\\
28	0.000661972198433341\\
29	0.00042616269315825\\
30	0.000166039916877417\\
31	9.7302548006578e-05\\
32	1.08233811293266e-05\\
33	7.56735550428782e-06\\
34	4.61640431864653e-06\\
35	3.9626085856082e-06\\
36	3.3323516361024e-06\\
37	2.67700477112937e-06\\
38	2.01466599699385e-06\\
39	1.55118162413949e-06\\
40	1.08847089041938e-06\\
41	7.47854813998003e-07\\
42	4.39262290727365e-07\\
43	3.56076778898183e-07\\
44	2.83575891658224e-07\\
45	2.51579461893206e-07\\
46	2.23122112553664e-07\\
47	1.77195140978287e-07\\
48	1.6250615260497e-07\\
49	1.53498912441218e-07\\
50	1.45326932681942e-07\\
51	1.22995974736719e-07\\
52	1.00912274688171e-07\\
53	7.56806988318392e-08\\
54	5.05155519837944e-08\\
55	3.70355539690048e-08\\
56	2.34123828560298e-08\\
57	1.52890070387948e-08\\
58	7.61628450490967e-09\\
};

\addplot 
  table[row sep=crcr]{2	1.0547118733939e-15\\
4	6.14977235802172e-11\\
6	1.99256167121575e-11\\
8	1.67640901160837e-12\\
10	5.04596364692134e-14\\
12	1.31512856160754e-13\\
14	1.16950428508122e-10\\
16	2.38072755953667e-12\\
18	4.63271782047414e-13\\
20	9.67195074030869e-14\\
22	9.32686003082828e-13\\
24	6.68936152130317e-13\\
26	6.01732205729455e-17\\
28	1.0336891932694e-14\\
30	1.38832359237312e-14\\
32	2.49891660098964e-17\\
34	1.46807386193948e-14\\
36	1.84903057220877e-17\\
38	3.28644548369932e-17\\
40	6.99443808942343e-16\\
42	8.36815612328045e-19\\
44	2.50932769470258e-16\\
46	3.68487125324202e-17\\
48	6.27808302946114e-15\\
50	4.05466757386533e-16\\
52	6.34546791501679e-19\\
54	4.57432202230804e-18\\
56	1.98639152298484e-19\\
58	2.76059427787211e-17\\
};

\addplot 
  table[row sep=crcr]{1	0\\
2	2.14841810691588e-14\\
3	0\\
4	2.78411922133488e-14\\
5	0\\
6	0.00115087021406837\\
7	0\\
8	0.00175915305975514\\
9	0\\
10	0.00195401151811336\\
11	0\\
12	0.00157639319792073\\
13	0\\
14	0.00125572693984344\\
15	0\\
16	0.0011441930482701\\
17	0\\
18	0.00119614755128791\\
19	0\\
20	0.00211522866575453\\
21	0\\
22	0.000315994435647878\\
23	0\\
24	0.000155698456904195\\
25	0\\
26	0.000148167979320008\\
27	0\\
28	0.000122783275084518\\
29	0\\
30	3.11339922776895e-05\\
31	0\\
32	3.38008031536389e-06\\
33	0\\
34	2.85819875229349e-06\\
35	0\\
36	2.18392576660854e-06\\
37	0\\
38	1.38523721322468e-06\\
39	0\\
40	1.19972042498105e-06\\
41	0\\
42	1.06875451245838e-06\\
43	0\\
44	6.94721720306515e-07\\
45	0\\
46	1.51649432457809e-07\\
47	0\\
48	1.32763268770142e-07\\
49	0\\
50	1.01613535362487e-07\\
51	0\\
52	2.65806393499771e-08\\
53	0\\
54	1.61012389649514e-08\\
55	0\\
56	6.53706847222442e-09\\
57	0\\
58	5.09492032046335e-09\\
};

\end{axis}
~
%
%

\begin{axis}[%
width=0.46\linewidth,
at={(cubeRI.north east)},
xshift=5em,
anchor=north west,
xmin=1,
xmax=32,
ymode=log,
ymin=1e-18,
ymax=1,
yminorticks=false,
cycle list name=resN,
cycle list shift=4,
xlabel={$j$},
every axis x label/.style={at={(ticklabel cs:0.5)},font=\small,rotate=0,anchor=center,yshift=-.5em},
every axis y label/.style={at={(ticklabel cs:0.5)},font=\small,rotate=90,anchor=center,yshift=.5em},
every  y tick label/.style={font=\small},
every  x tick label/.style={font=\small},
]
\addplot 
  table[row sep=crcr]{1	0.295280437231571\\
2	0.105500140012762\\
3	0.0915999775527808\\
4	0.0344723079843699\\
5	0.0139680735852687\\
6	0.012260766669331\\
7	0.011144106953418\\
8	0.0109081330980024\\
9	0.0107950846777802\\
10	0.0107758348454274\\
11	0.00765367565494039\\
12	0.000486363835434287\\
13	6.44345157205118e-05\\
14	2.99912866761558e-05\\
15	2.37113293613849e-05\\
16	1.64846341506818e-05\\
17	1.44088178437758e-05\\
18	9.70452214056634e-06\\
19	3.66142924911146e-06\\
20	3.65070344388593e-06\\
21	3.63204441547835e-06\\
22	3.61110305756079e-06\\
23	3.60982534606393e-06\\
24	3.56885608264532e-06\\
25	3.5322021843806e-06\\
26	3.45836946963871e-06\\
27	3.21084912006159e-06\\
28	2.25345146897407e-06\\
29	6.36591567244485e-07\\
30	6.85119804657415e-08\\
31	2.92142349350283e-08\\
32	5.87757982590785e-09\\
};

\addplot 
  table[row sep=crcr]{1	0.295280437231573\\
2	0.105500140012776\\
3	0.0915999775527811\\
4	0.0344723079843701\\
5	0.0139680735852687\\
6	0.0122607666693256\\
7	0.0111441069526405\\
8	0.0109081330980005\\
9	0.0107950846777794\\
10	0.0107758348454273\\
11	0.00765367565493822\\
12	0.000486363835426441\\
13	6.44345157205018e-05\\
14	2.99912866761862e-05\\
15	2.3711329361387e-05\\
16	1.64846341451153e-05\\
17	1.44088178437863e-05\\
18	9.70452214057375e-06\\
19	3.66142924926114e-06\\
20	3.65070344379311e-06\\
21	3.63204441539959e-06\\
22	3.61110305626566e-06\\
23	3.60982534613329e-06\\
24	3.56885608263584e-06\\
25	3.53220218444927e-06\\
26	3.45836946954125e-06\\
27	3.21084911994712e-06\\
28	2.25345146893675e-06\\
29	6.36591567232046e-07\\
30	6.85119804196987e-08\\
31	2.92142353579079e-08\\
32	5.87757993540234e-09\\
};

\addplot 
  table[row sep=crcr]{1	0.29528043723157\\
2	0.105500140012784\\
3	0.0915999775527812\\
4	0.0344723079843701\\
5	0.0139680735852686\\
6	0.0122607666693311\\
7	0.0111441069531865\\
8	0.0109081330980007\\
9	0.01079508467778\\
10	0.0107758348454272\\
11	0.00765367565494028\\
12	0.000486363835432729\\
13	6.44345157205048e-05\\
14	2.99912866761091e-05\\
15	2.37113293613746e-05\\
16	1.648463414471e-05\\
17	1.44088178437868e-05\\
18	9.7045221404511e-06\\
19	3.66142924904271e-06\\
20	3.65070344295389e-06\\
21	3.63204441537074e-06\\
22	3.6111030574079e-06\\
23	3.60982534562667e-06\\
24	3.56885608261598e-06\\
25	3.53220218433828e-06\\
26	3.45836946850314e-06\\
27	3.21084912014783e-06\\
28	2.25345146893415e-06\\
29	6.3659156726222e-07\\
30	6.8511980402257e-08\\
31	2.92142350085804e-08\\
32	5.87757990201506e-09\\
};
\addplot 
  table[row sep=crcr]{1	1.55431223447522e-15\\
2	1.44190215323192e-14\\
3	2.4980018054066e-16\\
4	2.0122792321331e-16\\
5	0\\
6	5.36549971119626e-15\\
7	7.77497857762377e-13\\
8	1.81972492629967e-15\\
9	8.11850586757146e-16\\
10	1.05818132034585e-16\\
11	2.16753698323302e-15\\
12	7.84620849195172e-15\\
13	9.93400240539843e-18\\
14	3.03915421474843e-17\\
15	2.08031291845656e-18\\
16	5.56651418210687e-15\\
17	1.04981263482698e-17\\
18	7.40815015668611e-18\\
19	1.49681733208149e-16\\
20	9.28187173930735e-17\\
21	7.87537353039158e-17\\
22	1.29513370514256e-15\\
23	6.93575988200239e-17\\
24	9.48549845982728e-18\\
25	6.86698080668534e-17\\
26	9.74532581639754e-17\\
27	1.14468879524893e-16\\
28	3.73202716560245e-17\\
29	1.24388905886662e-17\\
30	4.60428051886124e-17\\
31	4.22879564407819e-16\\
32	1.09494487463421e-16\\
};
\addplot 
  table[row sep=crcr]{1	7.7715611723761e-16\\
2	2.22599716437344e-14\\
3	3.19189119579733e-16\\
4	1.66533453693773e-16\\
5	5.89805981832114e-17\\
6	1.02348685082632e-16\\
7	2.31538746509052e-13\\
8	1.6306400674182e-15\\
9	2.27248775352962e-16\\
10	1.78676518025611e-16\\
11	1.15359111152458e-16\\
12	1.55854062294791e-15\\
13	6.91178884959509e-18\\
14	4.67799356109605e-17\\
15	1.03778476697597e-17\\
16	5.97172796780975e-15\\
17	1.09453597444201e-17\\
18	1.15245608737526e-16\\
19	6.87464645485799e-17\\
20	9.32040326807795e-16\\
21	1.07604101003871e-16\\
22	1.52893682144137e-16\\
23	4.37265512426982e-16\\
24	2.93391037105208e-17\\
25	4.23215776930875e-17\\
26	1.13556328579169e-15\\
27	8.62398124917493e-17\\
28	3.99274390676732e-17\\
29	1.77346464540185e-17\\
30	6.34844841219994e-17\\
31	7.35520761963251e-17\\
32	7.61072131202637e-17\\
};
\end{axis}
\end{tikzpicture}%
\end{figure}

As a final experiment, we test to what extent the predicted
equivalence~\cite[Theorem~2]{BenBKetal15} of ILRSI~\cite{LinSim15,MasOR16} and
RADI~\cite{BenBKetal15} holds in practice.  At first, both methods are executed
with the same set of predetermined shift parameters.  We mimic this by running
RADI with the adaptive residual-Hamiltonian shift strategy and then use the
generated shifts within ILRSI.\@
In the second experiment, we let ILRSI
compute its own shifts using the residual-Hamiltonian approach, too. The
\textit{Rail} example with $n=79\,841$ and the \textit{CUBE-FD} example with
$n=10\,648$, $m=p=10$ are used for this study.  Figure~\ref{fig:ILRSIvsRADI} shows
the scaled residual norms $\rho_j:=\|\cR_j\|/\|CC^{\trp}\|$ for RADI and ILRSI
with precomputed and adaptive shifts, as well as the norm differences
$\Delta
\rho_j:=\vert\|\cR^{\text{RADI}}_j\|-\|\cR^{\text{ILRSI}}_j\|\vert/\|CC^{\trp}\|$
against the iteration number $j$.

In both cases, we indeed observe that RADI and ILRSI produce visually nearly
indistinguishable residual curves when using exactly the same shifts.  When
ILRSI generates its own shifts for the \textit{CUBE-FD} example, the
discrepancies $\Delta\rho_j$ are larger, but still small. The reason is that for
this example, some of the generated shifts come in complex conjugated pairs, but
in its current implementation, ILRSI cannot handle those complex shifts
similarly well
as the RADI method. At some point, the low-rank solution factors generated by
ILRSI will be complex, which will lead to small differences in the computed
shifts since parts of the low-rank factor are used to project the Hamiltonian
matrices. Similar observations can be made for the \textit{Rail} example, but the
norm differences $\Delta\rho_j$ in case of adaptive shifts are much smaller compared to the \textit{CUBE-FD}
example as the majority of generated shifts for this example is real.
 In all
cases, the differences with respect to the computed approximate Riccati
solutions $\|X_j^{\text{RADI}}-X_j^{\text{ILRSI}}\|$ show the same behavior as
the residual norms differences.

\begin{figure}[t]
  \caption{\change{%
A summary of the performance of all methods over all
    examples. For each example, a method gains a penalty equal to its
    running time divided by running time of the fastest method for
    that example. If the quotient is larger than $10$, or the method
    failed in that example, the penalty is set to $10$. The smaller
    the total penalty, the faster the method.}
        }\label{fig:totalTimeComparison}
  \centering

\pgfplotstableread[row sep=\\,col sep=&]{
Method     & Chip        & Filter3D    & Rail        & Toeplitz    & Lung        & Stokes      & CUBE-10-10  & CUBE-10-20  & CUBE-20-10  & CUBE-20-10-unstable \\
EKSM       & 1           & 8.410546139 & 4.288196959 & 6.591836735 & 1.153846154 & 1.366515837 & 1.227920228 & 3.872496148 & 1.13372093  & 1 \\
RKSM       & 4.252747253 & 10          & 1.68211441  & 4.142857143 & 6.653846154 & 1.443438914 & 1           & 1.335130971 & 1           & 1.247928994 \\
TRKSM      & 4.505494505 & 2.15913371  & 1.686459088 & 4           & 10          & 10          & 1.1502849   & 1.596687211 & 1.177648579 & 1.250295858 \\
GEKSM      & 1.351648352 & 10          & 10          & 5.020408163 & 1           & 10          & 1.338319088 & 6.247303544 & 1.549741602 & 10 \\
ILRSI      & 5.857142857 & 1.490583804 & 1.304851557 & 5.979591837 & 10          & 4.117647059 & 1.286324786 & 1.12211094  & 1.398578811 & 1.783431953 \\
RADI       & 5.769230769 & 1           & 1           & 1           & 3.692307692 & 1           & 1.378917379 & 1.031201849 & 1.802971576 & 1.699408284 \\
NK-ADI+GP  & 5.208791209 & 10          & 1.451122375 & 2.244897959 & 10          & 1.022624434 & 1.038461538 & 1           & 1.086563307 & 1.843195266 \\
iNK-ADI+LS & 10          & 6.649717514 & 1.865314989 & 2.244897959 & 10          & 3.411764706 & 3.453703704 & 2.42449923  & 4.547157623 & 8.84260355 \\
    }\mydata

\definecolor{cChip}{RGB}{146,123,103}
\definecolor{cFilter}{RGB}{119,189,152}
\definecolor{cRail}{RGB}{224,119,97}
\definecolor{cToeplitz}{RGB}{249,199,140}
\definecolor{cLung}{RGB}{57,134,132}
\definecolor{cStokes}{RGB}{246,232,161}
\definecolor{cCUBEa}{RGB}{158,213,254}
\definecolor{cCUBEb}{RGB}{9,72,105}
\definecolor{cCUBEc}{RGB}{230,186,97}
\definecolor{cCUBEd}{RGB}{150,96,60}

\begin{tikzpicture}
  \begin{axis}[
    y=-0.55cm,
    xbar stacked,
    width=0.95\linewidth,
    symbolic y coords={EKSM,RKSM,TRKSM,GEKSM,ILRSI,RADI,NK-ADI+GP,iNK-ADI+LS},
    cycle list={%
        {draw=none,fill=cChip,postaction={pattern=dots,pattern color=black}},%
        {draw=none,fill=cFilter,postaction={pattern=grid,pattern color=black}},%
        {draw=none,fill=cRail,postaction={pattern=vertical lines,pattern color=black}},%
        {draw=none,fill=cToeplitz,postaction={pattern=horizontal lines,pattern color=black}},%
        {draw=none,fill=cLung,postaction={pattern=north west lines,pattern color=black}},%
        {draw=none,fill=cStokes,postaction={pattern=dots,pattern color=black}},%
        {draw=none,fill=cCUBEa,postaction={pattern=crosshatch,pattern color=black}},%
        {draw=none,fill=cCUBEb},%
        {draw=none,fill=cCUBEc},%
        {draw=none,fill=cCUBEd,postaction={pattern=north east lines,pattern color=black}}
      },
    yticklabel style={font=\small},
    xticklabel style={font=\small},
    ytick=data,
    legend style={at={(0.425,-0.1)}, anchor=north,legend columns=5,
      font=\small, rounded corners=3pt},
    ]

\addplot table[x=Chip,y=Method]{\mydata};
\addplot table[x=Filter3D,y=Method]{\mydata};
\addplot table[x=Rail,y=Method]{\mydata};
\addplot table[x=Toeplitz,y=Method]{\mydata};
\addplot table[x=Lung,y=Method]{\mydata};
\addplot table[x=Stokes,y=Method]{\mydata};
\addplot table[x=CUBE-10-10,y=Method]{\mydata};
\addplot table[x=CUBE-10-20,y=Method]{\mydata};
\addplot table[x=CUBE-20-10,y=Method]{\mydata};
\addplot table[x=CUBE-20-10-unstable,y=Method]{\mydata};

\legend{Chip,Filter3D,Rail,Toeplitz,Lung,Stokes,CUBE-10-10,CUBE-10-20,CUBE-20-10,CUBE-20-10-unst};

  \end{axis}
\end{tikzpicture}
\end{figure}
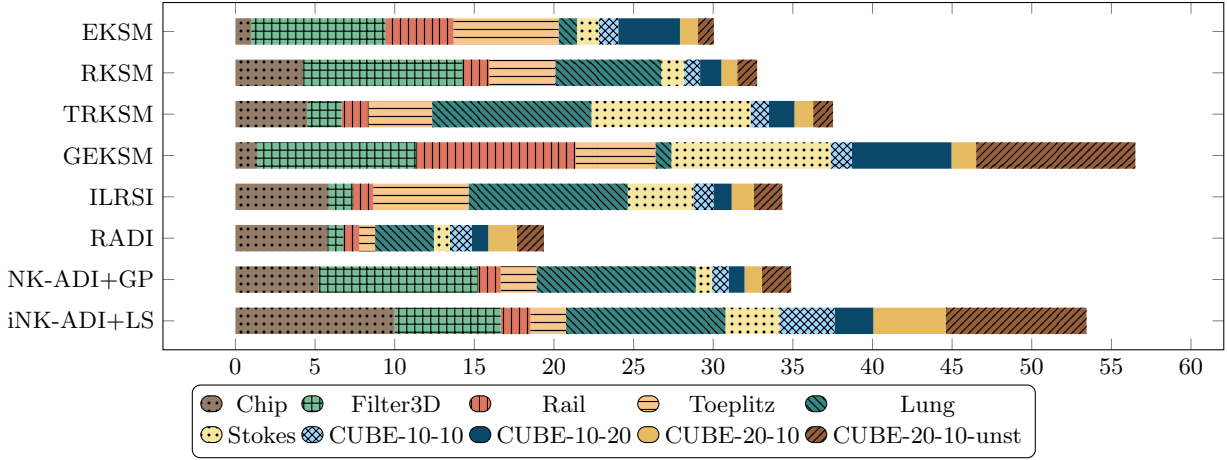

\section{Conclusions}\label{sec:concl}
In this paper, we studied three classes of low-rank algorithms for finding approximate solutions of large-scale GCAREs:
the projection based extended and rational Krylov
subspace methods (EKSM and RKSM), direct iterations not bound to a projection framework (ILRSI and RADI), and Newton-Kleinman methods. Some modifications of the
considered methods were also taken into account. We discussed the methods with respect to the computational complexity
of solving the matrix equation, and the amount of memory required. The theoretical analysis was confirmed by a number of numerical experiments.

We detected the most important subtasks that influence the computational work: solving linear systems with multiple right hand sides,
solving small, dense GCAREs, and generating shift parameters.
The only subtask present in each iteration of every method is the solution of linear systems, where the number of columns in the
right hand sides varies from one method to the other. This leads to different behavior regarding the ranks of the
constant and quadratic term defining the GCARE.\@ In all methods which do not rely on solving small, dense GCAREs (ILRSI, RADI, iNK-ADI+LS) these linear solves
constitute the largest portion of the overall computational work.
For the projection based methods ({(G)}EKSM, {(T)}RKSM, NK-ADI+GP), the solution of the dense
matrix equations could, depending on $p,m$ and the required number of iteration steps, be a second substantial portion of the work load. The generation of shift
parameters and the estimation of the GCARE residual norms are, due to several recent advances, often only minor fractions of the total computational cost.

Projecting to rational Krylov subspaces often fared better than projecting to extended Krylov subspaces, typically keeping the
built up subspace dimension close to the actual rank of the computed low-rank solution.
The extended Krylov subspace methods tend to generate subspaces of dimensions considerably larger than the rank of the approximate solution, and
this results in significant increase of the effort to solve the projected dense GCAREs.
The direct iterative methods, ILRSI and RADI, are theoretically equivalent which can also be observed in numerical tests where only smaller
deviations occur due to round off. In the present form, ILRSI still lacks some efficiency improvements found in the other methods, like reducing the occurrence
of complex arithmetic operations and a cheap GCARE residual norm estimation. Hence, ILRSI could sometimes not keep up with RADI or some of the other
approaches in our experiments. If these issues can be solved in future research
efforts, ILRSI can potentially become a very competitive approach.
The low-rank Newton-Kleinman iterations showed a mixed performance. While the Galerkin projection version could compete most of the time (provided the reduced
GCARE could be solved), the inexact version with line search could not keep up in most cases. The reasons were mostly the more expensive linear systems due to a
higher number of right hand sides, and the substantially larger number of iteration steps required.

Another considered performance indicator was the memory consumption, where it is important to distinguish between two scenarios. The
first scenario occurs when an
approximate low-rank solution $\tX\in\Rnn$ of the GCAREs is sought. Therefore, the smaller the constructed subspace or low-rank solution factor the better.
In this context, the {(T)}RKSM provided especially good results, followed by
RADI and ILRSI.\@
At the expense of somewhat larger generated subspace dimensions, NK-ADI+GP achieved sometimes several orders of magnitude smaller residual norms.
In the second scenario, only a stabilizing feedback matrix $\tK=M^{\trp}\tX B\in\Rnm$ is required, naturally asking for less data to be computed.
Obviously, projection based methods do not gain an advantage here, since they still have to
build the basis for the entire subspace.
The same is true for ILRSI.\@ Only RADI and iNK-ADI+LS were able to solely operate on approximate feedback matrices without ever
forming the low-rank solution factors, making them more memory efficient in this
situation.

 \change{%
Averaging over all carried out tests, {(T)}RKSM and RADI yielded
smaller subspace dimensions (low-rank solution factors) as well as
smaller total computation times (see
Figure~\ref{fig:totalTimeComparison}) compared to the other approaches. The ultimate choice between {(T)}RKSM
and RADI should then take into account whether the full low-rank solution, or
only a stabilizing feedback is sought.}

We also experimented with unstable GCAREs, which might ask for an initial stabilizing guess. The generation of such an initial guess proved to be demanding
in several tests because the conventional approach used for that purpose often failed. This issue requires further research effort.

In this study we restricted ourselves to sparse direct solution techniques for
the arising linear systems of equations. In general, preconditioned iterative
solvers could also be applied. This points towards the interesting research
topic of investigating the effects of the errors made by solving the linear
systems only inexactly, and to establish rules on the minimal accuracy of the linear solves required to obtain a low-rank GCARE solution of a certain quality.

Similar further comparative studies could address low-rank methods for continuous-time Lyapunov equations, and also discrete-time Lyapunov and
Riccati equations, as well as nonsymmetric variants. A comparison of different low-rank approaches for differential matrix equations is ongoing work.
Especially in the latter, but also in other applications, the inhomogeneities of the arising GCAREs can be indefinite, which requires some
smaller changes in the implementations of the algorithms, see, e.g.,~\cite{LanMS15,KreKM19}.
A more demanding alteration is when the quadratic term in the
GCARE is positive semidefinite (i.e., $M^{\trp}XBB^{\trp}XM$ occurs with a positive
sign in~\eqref{gcare}).
This arises, e.g., in certain model order reduction approaches, for which some numerical methods are proposed in~\cite{BenS14}.

\subsection*{Acknowledgments} We thank Valeria Simoncini, Tatjana
Stykel, Arash Massoudi, Peter Chang-Yi Weng, and Heiko Weichelt for helpful
discussions and sharing their implementations of {(T)}RKSM, GEKSM, ILRSI, SDA,
and iNK-ADI+LS.\@ Especially the latter were a great asset for developing the
implementations used in this study. Additional thanks go to Tony Stillfjord for
his careful revision of our final draft. This work was primarily generated while PK was affiliated with the Max Planck Institute for Dynamics of Complex Technical Systems in Magdeburg, Germany.

\bibliographystyle{siamplain}
\bibliography{BenBKEta_ARE_arxiv}
\end{document}